\documentclass[12pt]{article}
\usepackage[dvips]{epsfig}
\usepackage{amsmath,amssymb,amsfonts,enumerate,bbm}
\usepackage{caption}
\usepackage{url}
\usepackage{ifthen}
\usepackage[pdftex]{hyperref}
\usepackage{graphicx}

\usepackage{subfigure}

\usepackage{color,soul}

\long\def\comment #1\commentend{}

\renewcommand{\baselinestretch}{1.1}
\newtheorem{theorem}{Theorem}[section]
\newtheorem{lemma}[theorem]{Lemma}
\newtheorem{observation}[theorem]{Observation}
\newtheorem{corollary}[theorem]{Corollary}

\newtheorem{claim}[theorem]{Claim}
\newtheorem{fact}[theorem]{Fact}

\def\inline#1:{\par\vskip 7pt\noindent{\bf #1:}\hskip 10pt}
\def\Proof{\par\noindent{\bf Proof:~}}
\def\blackslug{\hbox{\hskip 1pt \vrule width 4pt height 8pt
    depth 1.5pt \hskip 1pt}}
\def\inQED{\quad\quad\blackslug}

\def\dnsitem{\vspace{-7pt}\item}
\def\ddnsitem{\vspace{-4pt}\item}
\def\dnsbibitem{\vspace{-4pt}\bibitem}

\newcommand{\PFT}[0]{PF Theorem}
\newcommand{\PFE}[0]{PF eigenvalue}
\newcommand{\PFR}[0]{PF root}
\newcommand{\PFV}[0]{PF vector}

\newcommand{\R}{\mathbb{R}}

\newcommand{\Sign}[0]{\mathrm{sign}}

\newcommand{\CP}[0]{\mathrm{P}} %characteristic polynomial

\newcommand{\MaxPower}[0]{X_{\max}}

\newcommand{\dist}[1]{\mathrm{d} ({#1})}
\newcommand{\EigenValue}[0]{EigVal}
\newcommand{\SpectralRatio}[0]{\rho}

\newcommand{\Entity}[0]{\mathit{v}}
\newcommand{\EntitySet}[0]{\mathcal{V}}
\newcommand{\Affectors}[0]{\mathcal{A}}
\newcommand{\Supporters}[0]{\mathcal{S}}
\newcommand{\Repressors}[0]{\mathcal{R}}
\newcommand{\SupportersMatrix}[0]{\mathcal{M}^+}
\newcommand{\RepressorsMatrix}[0]{\mathcal{M}^-}

\newcommand{\WeakSupportersMatrix}[0]{\mathcal{M}^{+}_{w}}
\newcommand{\WeakRepressorsMatrix}[0]{\mathcal{M}^{-}_{w}}

\newcommand{\PFEigenValue}[0]{\mathit{r}}
\newcommand{\PFEigenVector}[0]{\mathbf{\overline{P}}}
\newcommand{\Period}[0]{\mathit{h}}
\newcommand{\FilterMatrix}[0]{F}
\newcommand{\FilterMatrixFamily}[0]{\mathcal{F}}

\newcommand{\ConstraintsGraph}[0]{\mathcal{CG}}
\newcommand{\GCGraph}[0]{CG}
\newcommand{\CGLevel}[0]{L}

\newcommand{\Polytope}[0]{\mathcal{P}}
\newcommand{\ZeroStar}[0]{\boldsymbol{0^{*}}}
\newcommand{\WeakZeroStar}[0]{\boldsymbol{{}_w0^{*}}}
\newcommand{\Zero}[0]{\boldsymbol{0^{f}}}
\newcommand{\WeakZero}[0]{\boldsymbol{{}_w0^{f}}}
\newcommand{\MaxGain}[0]{\mathcal{G}_{max}}
\newcommand{\AlgoName}[0]{\tt Compute\PFEigenVector(\System)}
\newcommand{\TestIrred}[0]{Irr\_Test}

\newcommand{\System}[0]{\mathcal{L}}

\newcommand{\WeakSystem}[0]{\mathcal{L}^{W}}
\newcommand{\SystemFamily}[0]{\mathfrak{L}^{NS}}
\newcommand{\SquareSystemFamily}[0]{\mathfrak{L}^{S}}
\newcommand{\WeakSystemFamily}[0]{\mathfrak{L}^{W}}

\newcommand{\SCC}[0]{strongly connected component}
\newcommand{\BFS}[0]{BFS}

\newcommand{\TotR}[0]{T^-}
\newcommand{\TotS}[0]{T^+}

\newcommand{\LPRunTime}[0]{\mbox{T}_{LP}}
\newcommand{\SelectionVec}[0]{\bold{S}}
\newcommand{\fff}[0]{f}

\def\inline#1:{\par\vskip 7pt\noindent{\bf #1:}\hskip 10pt}
\def\Proof{\par\noindent{\bf Proof:~}}
\def\blackslug{\hbox{\hskip 1pt \vrule width 4pt height 8pt
    depth 1.5pt \hskip 1pt}}
\def\QED{\quad\blackslug\lower 8.5pt\null\par}
\def\inQED{\quad\quad\blackslug}
\newcommand{\IndS}[0]{\gamma}
\begin{document}
%%%%%%%%%%%%%%%%%%%%%%

\title{Generalized Perron--Frobenius Theorem for Nonsquare Matrices}
\author{
Chen Avin
\thanks{
%Department of Communication Systems Engineering,
\hbox{Ben Gurion University, Beer-Sheva, Israel. Email:}
{\tt \{avin,}
{\tt borokhom,}
{\tt zvilo\}@cse.bgu.ac.il,}
{\tt yoram.haddad@gmail.com.}}
\and
Michael Borokhovich $^*$
\and
Yoram Haddad
\thanks{Jerusalem College of Technology, Jerusalem, Israel.
}
$^*$
\and
Erez Kantor
\thanks{The Technion, Haifa, Israel.
Email: {\tt erez.kantor@gmail.com}. Supported by Eshkol fellowship,
the Israel Ministry of Science and Technology.}
\and
Zvi Lotker $^*$
\thanks{Supported by a grant of the Israel Science Foundation.}
\and
Merav Parter
\thanks{
%Department of Computer Science and Applied Mathematics,
The Weizmann Institute of Science, Rehovot, Israel.
Email: {\tt \{merav.parter,david.peleg\}@ weizmann.ac.il}.}
\thanks{Recipient of the Google Europe Fellowship in distributed computing;
 research supported in part by this Google Fellowship.}
\and
David Peleg $^\ddag$\thanks{Supported in part by the Israel Science Foundation (grant 894/09),
the United States-Israel Binational Science Foundation
(grant 2008348),
the I-CORE program of the Israel PBC and ISF (grant 4/11),
the Israel Ministry of Science and Technology
(infrastructures grant), and the Citi Foundation.}
}
%\institute{
%Dept. of Commun. Syst. Eng., Ben Gurion University,
%%Department of Communication Systems Engineering, Ben Gurion University,
%Beer-Sheva, Israel.
%E-mail:{\tt \{avin,coasaf,zvilo\}@cse.bgu.ac.il, yoram.haddad@gmail.com}.
%\and
%Faculty of Elect. Eng., The Technion, Haifa, Israel.\\
%E-mail:{\tt erez.kantor@gmail.com}.
%\and
%Dept. of Computer Sci., The Weizmann Institute, Rehovot, Israel\\
%%Department of Computer Science and Applied Mathematics,\\
%%The Weizmann Institute of Science, Rehovot, Israel\\
%\email{\{merav.parter,david.peleg\}@weizmann.ac.il}
%}
\maketitle
\begin{abstract}
The celebrated Perron--Frobenius (PF) theorem is stated for irreducible nonnegative square matrices, and provides a simple characterization of their eigenvectors and eigenvalues.
The importance of this theorem stems from the fact that eigenvalue problems
on such matrices arise in many fields of science and engineering, including dynamical systems theory, economics, statistics and optimization.
However, many real-life scenarios give rise to nonsquare matrices.
Despite the extensive development of spectral theories for nonnegative matrices, the applicability of such theories to non-convex optimization problems is not clear. In particular, a natural question is whether the \PFT~(along with its applications) can be generalized to a nonsquare setting.
Our paper provides a generalization of the \PFT{ }
to nonsquare matrices.
The extension can be interpreted as representing client-server systems with additional
degrees of freedom, where each client may choose between multiple
servers that can cooperate in serving it
(while potentially interfering with other clients).
This formulation is motivated by applications to power control
in wireless networks, economics and others,
all of which extend known examples for the use of the original \PFT.

We show that the option of cooperation between servers does not
improve the situation, in the sense that in the optimal solution
no cooperation is needed, and only one server needs to serve each client.
Hence, the additional power of having several potential servers per client
translates into \emph{choosing} the best single server and not into \emph{sharing} the load between the servers in some way, as one might have expected.

The two main contributions of the paper are
(i) a generalized \PFT { }that characterizes the optimal solution
for a non-convex nonsquare problem, and
(ii) an algorithm for finding the optimal solution in polynomial time.
Towards achieving those goals, we extend the definitions of irreducibility and largest eigenvalue
of square matrices to nonsquare ones in a novel and non-trivial way,
which turns out to be necessary and sufficient for our generalized theorem
to hold.
The analysis performed to characterize the optimal solution uses techniques from a wide range of areas and exploits combinatorial
properties of polytopes, graph-theoretic techniques and analytic tools
such as spectral properties of nonnegative matrices and root characterization
of integer polynomials.
%-------------------------------------

%We call such a solution a $\ZeroStar$ solution, or the
%\emph{generalized \PFT} of the nonsquare system.
%Establishing the existence of a $\ZeroStar$ solution, we then provide
%a polynomial time algorithm to compute it.  In addition the paper defines
%a notion of irreducibility of nonsquare system (in away that is both
%sufficient and necessary) and provide a polynomial time algorithm to test
%whether a system of Nonsquare Matrices is irreducible.
%A key notion in our analysis is the constraint graph, which represents
%the flow of influence between the constraints of the system.
%This graph can be though of the square-essence of the nonsquare essence
%and plays an essential rule in our analysis.

\end{abstract}

%\newpage
%%%%%%%%%%%%%%%%%%%%%%%%%%%%%%%%%%%%%%%%%%%%
\section{Introduction}
%%%%%%%%%%%%%%%%%%%%%%%%%%%%%%%%%%%%%%%%%%%%
\paragraph{Motivation and main results.}
%\subsection{Motivation and main results}
%%%%%%%%%%%%%%%%%%%%%%%%%%%%%%%%%%%%%%%%%%%%
This paper presents a generalization of the well known Perron--Frobenius (PF)
Theorem \cite{PF_Frobenius,PF_Perron}.
As a motivating example, let us consider the \emph{Power control problem}, one of the most fundamental problems in wireless networks.
The input to this problem consists of $n$ receiver-transmitter pairs and their physical locations.
All transmitters are set to transmit at the same time with the same frequency,
thus causing interference to the other receivers.
Therefore, receiving and decoding a message at each receiver depends on the
transmitting power of its paired transmitter as well as the power of the rest
of the transmitters.
%The \emph{physical} model we use to represent this setting
%is known as {\emph{signal-to-interference \& noise ratio}} (SINR) model.
%In the SINR model, the energy of a signal fades with the distance and
If the \emph{signal to interference ratio} at a receiver, namely, the signal strength received by a receiver divided by the interfering
strength of other simultaneous transmissions,
%(plus the fixed \emph{background noise} \(\Noise\))
is above some \emph{reception threshold} $\beta$, then the
receiver successfully receives the message, otherwise it does not \cite{R96}.
The power control problem is then to find an optimal power assignment
for the transmitters, so as to make the reception threshold $\beta$
as high as possible and ease the decoding process.

As it turns out, this power control problem can be solved elegantly by casting it as an optimization program and using the Perron--Frobenius (PF) Theorem \cite{Zander92b}.
%This paper extends the human knowledge in two ways.
%First we show how to use the ideas from ellipsoid algorithm and
%convex programming algorithms to solve a family of non-convex problems.
%Traditionally,  convexity plays a central role in developing algorithms,
%and there are few examples of problems that are not convex and can be solved
%in polynomial time accurately.
%Second we extend the well known theorem of Perron--Frobenius (PF)
%from square to nonsquare matrices.
The theorem can be formulated as dealing with
the following optimization problem (where $A \in \R^{n \times n}$):
\begin{eqnarray}\label{eq:basic}
&&\text{maximize $\beta$ subject to: }\\
&&A \cdot \overline{X} \leq 1/\beta \cdot \overline{X},~~
\displaystyle ||\overline{X}||_{1}=1,~~
\displaystyle \overline{X} \geq \overline{0}.\nonumber
\end{eqnarray}
%\begin{equation}\label{eq:basic}
%\framebox{ $
%\begin{array}{lll}
%  \max & \beta
%  \\
%  \mathrm{s.t.}
%  & \displaystyle A \cdot \overline{X} \leq 1/\beta \cdot \overline{X} &
%  \\
%  & \displaystyle \bold{1}^{T} \overline{X}=1 &
%  \\
%  & \displaystyle \overline{X} \geq \overline{0} &
%\end{array}
%%$ }
%\end{equation}
Let $\beta^*$ denote the optimal solution for Program (\ref{eq:basic}).
The Perron--Frobenius (PF) Theorem characterizes the solution to this
optimization problem and shows the following:
\begin{theorem} {\sc (\PFT, short version, \cite{PF_Frobenius,PF_Perron})}
%\begin{itemize}
Let $A$ be an irreducible nonnegative  matrix. Then $\beta^* = 1/\PFEigenValue$, where $\PFEigenValue \in \R_{>0}$
is the largest eigenvalue of $A$,
%(the \emph{Perron-Frobenius root} of $A$).
called the \emph{Perron--Frobenius (PF) root} of $A$.
There exists a unique (eigen-)vector $\PFEigenVector>0$,
$||\PFEigenVector||_{1}=1$, such that
$A \cdot \PFEigenVector= r \cdot \PFEigenVector$,
%(the \emph{Perron vector} of $A$).
called the \emph{Perron vector} of $A$.
(The pair $(\PFEigenValue,\PFEigenVector)$ is hereafter referred to as an {\em eigenpair}
of $A$.)
\end{theorem}
%Finding $\beta^*$ and $\mathbf{p}$ is very useful and has many applications,
%for examples:
%\begin{itemize}
%\item In power control problems in wireless networks where each receiver
%is paired with a single transmitter and vice versa,
%$\mathbf{p}$ is an optimal power assignment and $\beta^*$ is the best
%reception threshold that could be achieved.
%\item In Leontief’s Input--Output Economic Model where each factory produces
%a single product, $\mathbf{p}$ gives the pricing and $\beta^*$ is the largest
%minimum profit for a factory.
%\item In the Leslie (and others) Population model $\mathbf{p}$ is the
%steady state of the system.
%\item In Markov chains $\mathbf{p}$ is the unique stationary distribution
%and $\beta^* = 1$.
%\end{itemize}
Returning to our motivating example, let us consider a more complicated
variant of the power control problem,
where each receiver has several transmitters that can transmit to it
(and only to it) synchronously.
Since these transmitters are located at different places, it may conceivably
be better to divide the power (or work) among them, to increase
the reception threshold at their common receiver. Again, the question
concerns finding the best power assignment among all transmitters.

%To handle this multiple transmitters scenario,
%
In this paper
we extend Program (\ref{eq:basic}) to \emph{nonsquare matrices} and consider
the following extended optimization problem,
which in particular captures the multiple transmitters scenario.
(Here $A, B \in \R^{n \times m}$, $n \leq m$.)
%maximize $\beta$ subject to:
\begin{eqnarray}\label{eq:extended}
&&\text{maximize $\beta$ subject to: }~~\\
&&A \cdot  \overline{X} \leq 1/\beta \cdot B \cdot \overline{X},
\displaystyle ~~~||\overline{X}||_{1}=1,~~~
\displaystyle \overline{X} \geq  \overline{0}.\nonumber
\end{eqnarray}
%\begin{equation}\label{eq:extended}
%\framebox{ $
%\begin{array}{lll}
%  \max & \beta
%  \\
%  \mathrm{s.t.}
% & \displaystyle A \cdot\overline{X} \leq 1/\beta \cdot B \cdot %\overline{X}&
%  \\  %\label{eq:SR}
%  & \displaystyle \bold{1}^{T} \overline{X}=1 &
%  \\
%  & \displaystyle \overline{X} \geq  \overline{0}&
%  %\label{eq:Greater1}
%\end{array}
%%$ }
%\end{equation}
%
%We interpret the nonsquare matrix as some freedom the system designer have. This freedom appears in many applications of PF. For example, consider the following problem. We are given a sets of $R = \{r_1,...,r_n\}$ of $n$ receivers embedded in the plane. For each receiver there is a set of transmitters. The transmitters can cooperate. We would like to compute the bust way the transmitters can cooperate. I.e. witch transmitters transmit and with what power so that you can decode all broadcasts at the maximum SINR.
%
We interpret the nonsquare matrices $A,B$ as representing some additional
freedom given to the system designer. In this setting,
each \emph{entity} (receiver, in the power control example) has
several \emph{affectors} (transmitters, in the example),
referred to as its \emph{supporters}, which can cooperate in serving it
and share the workload. In such a general setting, we would like
%For example in a wireless network setting this problem appears when we allow multiple (and not only single) transmitters to transmit to the same receivers, at the same time, but from different locations. Or, in the economy model we allow each factory to manufacture multiple (and not a single) product to maximize its profit.
to find the best way to organize the cooperation between the supporters
of each entity.

The original problem was defined for a square matrix, so the appearance of
eigenvalues in the characterization of its solution seems natural. In contrast, in the generalized setting the
situation seems more complex. Our main result is an extension of the
\PFT~to nonsquare matrices and systems that give rise to an optimization problem in the form of (\ref{eq:extended}), with optimal solution $\beta^*$.

\begin{theorem}
{\sc (Nonsquare \PFT, short version)}
Let $\langle A, B \rangle$ be an irreducible nonnegative system
%(in a sense
(to be made formal later).
Then $\beta^* = 1/\PFEigenValue$, where $\PFEigenValue \in \R_{>0}$ is
the smallest \emph{Perron--Frobenius (PF) root} of all ${n \times n}$ square
sub-systems (defined formally later).
There exists a vector $\PFEigenVector \ge 0$ such that
$A \cdot \PFEigenVector = \PFEigenValue \cdot B \cdot \PFEigenVector$ and $\PFEigenVector$ has $n$
entries greater than 0 and $m-n$ zero entries
(referred to as a $\ZeroStar$ solution).
\end{theorem}

In other words, the theorem implies that the option of cooperation
does not improve the situation, in the sense that in the optimum solution,
no cooperation is needed and only one supporter per entity needs to work.
Hence, the additional power of having several potential supporters per entity
translates into \emph{choosing} the best single supporter and not into \emph{sharing} the load between the supporters in some way, as one might have expected.
\par As it turns out, the lion's share of our analysis involves such a characterization
of the optimal solution for (the non-convex) problem
of Program (\ref{eq:extended}).
The main challenge is to  show that at the  optimum, there exists a solution
in which only one supporter per entity is required to work;
we call such a solution a $\ZeroStar$ solution.
Namely, the structure that we establish is that the optimal solution
for our nonsquare
system is in fact the optimal solution of an \emph{embedded} square PF system.
Indeed, to enjoy the benefits of an equivalent square system, one should show that there exists a solution in which only one supporter per entity is required to work.
Interestingly,
%the effort in proving that one supporter per entity is sufficient is not linear. In fact,
it turned out to be relatively easy to show that there exists an optimal
``almost $\ZeroStar$'' solution, in which each entity
\emph{except at most one} has a single active
supporter and the remaining entity has at most \emph{two} active supporters.
Despite the presumably large ``improvement" of decreasing the number of servers from $m$ to $n+1$, this still leaves us in the frustrating situation of a nonsquare $n \times (n+1)$ system, where no spectral characterization for optimal solutions exists. In order to allow us to characterize the optimal solution using
the eigenpair of the best square matrix embedded within the nonsquare system,
one must overcome this last hurdle, and reach the ``phase transition'' point of $n$ servers, in which the system is \emph{square}.
Our main efforts went into showing that the remaining entity, too,
can select just one supporter while maintaining optimality, ending with
a \emph{square} $n\times n$ irreducible system where the traditional \PFT\
can be applied.
Proving the existence of an optimal $\ZeroStar$ solution requires techniques
from a wide range of areas to come into play and provide a rich understanding
of the system on various levels. In particular, the analysis exploits
combinatorial properties of polytopes, graph-theoretic techniques and
analytic tools such as spectral properties of nonnegative matrices and
root characterization of integer polynomials.
%We call such a solution a $\ZeroStar$ solution.

In the context of the above example of power control in wireless network with multiple
transmitters per receiver, a $\ZeroStar$ solution means that the best
reception threshold is achieved when only a single transmitter transmits
to each receiver.
%Fig. \ref{figure:sinr_diag}(a) shows the optimal $\ZeroStar$ solution, where each receiver $s_i$ is the reception zone of one of its transmitters $t_{ij}$, $j=1,2$. Fig. \ref{figure:sinr_diag}(b) shows what would be happen if we choose a non-optimal $\ZeroStar$ solution. We see that now the receivers are not covered by any of their transmitters. Finally,  Fig. \ref{figure:sinr_diag}(c) shows that even when two transmitters collaborate, receivers may not be covered by their reception zones.
%The lion's share of our analysis involves a characterization of the optimal solution for (the non-convex) problem of Program (\ref{eq:extended}).
Other known applications of the \PFT~can also be extended in a similar manner.
An Example for such applications is the {\em input-output economic model}
\cite{pillai2005pft}.
%and the {\em population growth model}
%\cite{meyer2000matrix}.
In this economic model, each industry produces a commodity and buys commodities
(raw materials) from other industries. The percentage
profit margin of an industry is the ratio
of its total income and total expenses (for buying its raw materials). It is required to find a pricing maximizing the ratio of the total income and total expenses of all industries. The extended PF variant of the problem concerns the case
where an industry can produce multiple commodities instead of just one.
%In the population growth model, the population is divided into groups
%(e.g., age groups) and there are migrations or splits between the groups
%(via births, deaths, and aging). The question is to find the distribution
%over ages that maximizes the population growth rate. An extended PF variant
%may consider the case where each group has several different types.
In this example, the same general phenomenon holds:
each industry should charge money only for \emph{one} of the commodities it produces.  That is, in the optimal pricing, one commodity per industry has nonzero price,  therefore the optimum is a $\ZeroStar$ solution.
For a more detailed discussion of applications, see Sec. \ref{short:sec:Applications}. In addition,  in Sec. \ref{sec:limit}, we provide a characterization of systems in which a $\ZeroStar$ solution does not exist.
\par While in the original setting the \PFR\ and \PFV\ can be computed in polynomial time,
this is not clear in the extended case, since the problem is not convex
\cite{Boyd-Conv-Opt-Book} (and not even log-convex) and there are
exponentially many choices in the system even if we know that
the optimal solution is  $\ZeroStar$
and each entity has only two supporters
to choose from. Our second main contribution is providing a polynomial time
algorithm to find  $\beta^*$ and $\PFEigenVector$.
The algorithm uses the fact that fixing $\beta$ yields a relaxed problem
which is convex (actually it becomes a linear program). This allows us to employ
the well known interior point method \cite{Boyd-Conv-Opt-Book},
for testing a specific $\beta$ for feasibility.
Hence, the problem reduces to finding the maximum feasible $\beta$,
and the algorithm does so by applying binary search on $\beta$.
Clearly, the search results in an approximate solution, in fact yielding
a fully polynomial time approximation scheme (FPTAS) for program (\ref{eq:extended}). This, however, leaves open the
intriguing question of whether program (\ref{eq:extended}) is polynomial.
Obtaining an exact optimal $\beta^*$, along with an appropriate vector
$\PFEigenVector$, is thus another challenging aspect of the problem.
%Note that although there exists FPTAS to Program (\ref{eq:extended}),
%(via binary search on $\beta^{*}$), due to the special structure of the optimal solution, it can be computed \emph{exactly}, showing that Program [\ref{eq:extended}] is in P. Essentially, this exact solution attains a spectral characterization despite the fact the the system is nonsquare.
\par A central notion in the generalized PF theorem is the \emph{irreducibility}
of the system. While irreducibility is a well-established concept for square
systems, it is less obvious how to define irreducibility for a nonsquare matrix
or system as in Program \eqref{eq:extended}. We provide a suitable definition
% of irreducibility for nonsquare systems that is
based on the property that every maximal square (legal) subsystem is
irreducible, and show that our definition is necessary and sufficient
for the theorem to hold.
%But, since there could be exponentially many such
%square subsystems, it is not a priori clear if one can check that a
%nonsquare system is irreducible in polynomial time.
A key tool in our analysis is what we call the \emph{constraint graph} of
the system, whose vertex set is the set on $n$ constraints (one per entity)
and whose edges represent direct influence between the constraints.
For a square system, irreducibility is equivalent to the constraint graph
being strongly connected, but for nonsquare systems the situation is more
delicate. Essentially, although the matrices are not square, the notion of
constraint graph is well defined and provides a valuable \emph{square}
representation of the nonsquare system (i.e., the adjacency matrix of
the graph). In \cite{PF_Irred,PF_Archive}, we also present a polynomial-time algorithm for testing the irreducibility of a given system, which exploits the properties
of the constraint graph.
\par\noindent{\bf Related work.}
%%%%%%%%%%%%%%%%%%%%%%%%%%%%%%%%%%%%%%%%%%%%
The \PFT~establishes the following two important ``PF properties" for a nonnegative square matrix $A \in \R^{n \times n}$:
(1) the \emph{Perron--Frobenius property}:
$A$ has a maximal nonnegative eigenpair.
If in addition the matrix $A$ is \emph{irreducible} then its maximal eigenvalue is strictly positive, dominant and with a strictly positive eigenvector. Thus nonnegative irreducible matrix $A$ is said to enjoy
the \emph{strong Perron--Frobenius property} \cite{PF_Frobenius,PF_Perron}.
(2) the \emph{Collatz--Wielandt property} (a.k.a. min-max characterization):
the maximal eigenpair is the optimal solution of Program (\ref{eq:basic}) \cite{PF_Collatz, PF_Wielandt}.
\par Matrices with these properties have played
an important role in a wide variety of applications.
The wide applicability of the \PFT, as well as the fact that the necessary and sufficient properties required of a matrix $A$ for the PF properties to hold are still not fully understood, have led to the emergence of many generalizations. We note that whereas all generalizations concern the Perron--Frobenius property, the Collatz--Wielandt property is not always established.
The long series of existing PF extensions includes \cite{PF_with_some_negative_entries,PF_evNONNEG,PF_complex_matrices,PF_for_non_linear_mapping, PF_nonlinear_more,PF_concave_mappings,PF_for_matrix_polynomials,AXBX}.
We next discuss these extensions in comparison to the current work.
%In addition, in Section \ref{subsec:power_control_app}
%we discuss the existing literature for the wireless power control problem
%with multiple transmitters.

Existing PF extensions can be broadly classified into four classes.
The first concerns matrices that do not satisfy the irreducibility and nonnegativity requirements. For example, \cite{PF_with_some_negative_entries,PF_evNONNEG} establish the Perron-Frobenius property for \emph{almost} nonnegative matrices or \emph{eventually} nonnegative matrices.
%and \cite{PF_evNONNEG} considers matrices that are \emph{eventually} nonnegative or positive.
A second class of generalizations concerns square matrices over different domains. For example, in  \cite{PF_complex_matrices}, the \PFT~was established for complex matrices $A \in  \mathbb{C}^{n \times n}$.
In the third type of generalization, the linear transformation obtained by applying the nonnegative irreducible matrix $A$ is generalized to a nonlinear mapping \cite{PF_for_non_linear_mapping, PF_nonlinear_more}, a concave mapping \cite{PF_concave_mappings} or a matrix polynomial mapping \cite{PF_for_matrix_polynomials}.

Last, a much less well studied generalization deals with nonsquare matrices,
i.e., matrices in $\R^{n \times m}$ for $m \neq n$.
Note that when considering a nonsquare system, the notion of eigenvalues
requires definition. There are several possible definitions for eigenvalues
in nonsquare matrices.
One possible setting for this type of generalizations considers a pair
of nonsquare  ``pencil" matrices $A, B \in \R^{n \times m}$,
where the term ``pencil" refers to the expression $A- \lambda \cdot B$,
for $\lambda \in \mathbb{C}$. Of special interest here are the values
that reduce the pencil rank, namely, the $\lambda$ values satisfying
$(A -\lambda B) \cdot \overline{X}=\overline{0}$
for some nonzero $\overline{X}$.
This problem is known as the \emph{generalized eigenvalue problem} \cite{AXBX,NonSQPencil,boelgomi05,Kres11},
which can be stated as follows:
Given matrices $A, B \in \R^{n \times m}$, find a vector $\overline{X}\neq \overline{0}$, $\lambda \in \mathbb{C}$, so that $A \cdot\overline{X}=\lambda B \cdot \overline{X}$. The complex number $\lambda$ is said to be an \emph{eigenvalue of $A$ relative to $B$} iff $A \overline{X}=\lambda \cdot B
\cdot \overline{X}$ for some nonzero $\overline{X}$ and $\overline{X}$ is called the \emph{eigenvector of $A$ relative to $B$}. The set of all eigenvalues of $A$ relative to $B$ is called the \emph{spectrum of $A$ relative to $B$}, denoted by $sp(A_{B})$.

Using the above definition, \cite{AXBX} considered
pairs of nonsquare matrices $A,B$ and was the first to characterize
the relation between $A$ and $B$
required to establish their PF property,
%in a way that establishes their PF property,
i.e., guarantee that the generalized eigenpair is nonnegative.
Essentially, this is done by generalizing the notions of positivity and
nonnegativity in the following manner. A matrix $A$ is said to be
\emph{positive} (respectively,\emph{nonnegative}) with respect to $B$,
if $B^{T} \cdot \overline{Y} \geq 0$ implies that $A^{T} \cdot \overline{Y}>0$
(resp., $A^{T} \cdot \overline{Y}\geq 0$). Note that for $B=I$,
these definitions coincide with the classical definitions of a positive
(resp., nonnegative) matrix.  Let $A, B \in \R^{n \times m}$, for $n \geq m$,
be such that the rank of $A$ or the rank of $B$ is $n$.
It is shown in \cite{AXBX} that if $A$ is positive (resp., nonnegative)
with respect to $B$, then the generalized eigenvalue problem
$A \cdot\overline{X}=\lambda \cdot B \cdot \overline{X}$ has a discrete and finite spectrum,
the eigenvalue with the largest absolute value is real and positive
(resp., nonnegative), and the corresponding eigenvector is positive
(resp., nonnegative). Observe that under the definition used therein,
the cases where $m > n$ (which is the setting studied here) is uninteresting,
as the columns of $A -\lambda \cdot B$ are linearly dependent for any real
$\lambda$, and hence the spectrum $sp(A_{B})$ is unbounded.
%Overall, the importance of \cite{AXBX} was in providing definitions for eigenvalues and eigenvector of a pair of nonsquare matrices and establishing the connection to \PFT~in the sense that they characterizes pairs of nonsquare matrices for which the PF property hold.

Despite the significance of \cite{AXBX} and its pioneering generalization of
the \PFT~to nonsquare systems, it is not clear what are the applications
of such a generalization, and no specific implications are known for the
traditional applications of the PF theorem.
%, such as the power-control problem
%or the economy model.
Moreover, although \cite{AXBX} established the PF
property for a class of pairs of  nonsquare matrices,
the Collatz--Wielandt property, which provides the algorithmic power for the
\PFT, does not necessarily hold with the spectral definition of \cite{AXBX}.
%\textbf{MP: I have an example for it. It is not stated explicitly in their paper.\\
%Let $A=[-1 ~2, 0 ~1, 0~1]$ and $B=[-1~ 3, 0~ 1, 0 ~1]$. Then for nonnegative $X=[1 ~1,0 ~1]$ it holds that $BX=A$ and thus the condition of \cite{AXBX} holds. Let $\lambda=1$ and $x=(1~ 0)$ it can be verified that $A x = B x$. However for any $\lambda>1$ we get that  $A x \geq \lambda B x$ (where there is a strict inequality in the first coordinate and equalities in the second and third). In other words if the goal is: $\max \lambda$ s.t $A x \geq \lambda Bx$ for $x\geq 0$, $x \neq 0$, this optimization problem is not bounded. Although $sp(A_{B})$ is bounded (by 1).}\\
%Therefore the spectral radius of $A$ with respect to $B$ doesn't correspond to a value of an optimization problem as in the classical PF. \\
%\textbf{MP: I wrote this but I am not sure in this since the Collatz-Wielandt property holds also for nonnegative (reducible) matrices}.\\
%Indeed it seems crucial that for having an optimization problem interpretation for the spectral characterization one needs to define the notion of irreducibility. Essentially, the irreducibility of the system implies that there exists a proper graph representation of the system (the pair of nonsquare matrices) that is strongly connected. The strong connectivity reflects the interplay between the (nonnegative) coordinates of a nonnegative solution. The balance of such interplay is then accomplished by the nonnegative eigenvector of the system.
In addition, since no notion of irreducibility was defined in \cite{AXBX},
the spectral radius of a nonnegative system (in the sense of the definition
of \cite{AXBX}) might be zero, and the corresponding eigenvector might be
nonnegative in the strong sense (with some zero coordinates).
These degenerations can be handled only by considering irreducible
nonnegative matrices, as was done
by Frobenius in \cite{PF_Frobenius}.

In contrast, the goal of the current work is to develop the spectral theory for a pair of
nonnegative matrices in a way that is both necessary and sufficient for both
the %strong
PF property and the Collatz--Wielandt property to hold
(allowing the nonsquare system to be of the
``same power" as the square systems considered by Perron and Frobenius).
%We consider nonsquare matrices of dimension $n \times m$ for $n \leq m$,
%which can be interpreted as describing a system with nonsquare (of columns) per entity (row).
Towards this we define the eigenvalues and eigenvectors of pairs of $n \times m$ matrices $A$ and $B$ in a novel manner. Such eigenpair $(\lambda, \overline{X})$ satisfies $A \cdot \overline{X}=\lambda \cdot B \cdot \overline{X}$. In \cite{AXBX}, alternative spectral definitions for pairs of nonsquare matrices $A$ and $B$ are provided. We note that whereas in \cite{AXBX} formulation, the maximum eigenvalue is not bounded if $n < m$, with our definition it is bounded.
%
%Interestingly, the maximum eigenvalue of the spectrum we define is also the maximum of the spectrum according to the definition of \cite{AXBX} and therefore we can show that the Collatz-Wielandt property is extended as well.
\par Let us note that although the generalized eigenvalue problem has been studied for many years, and multiple approaches for nonsquare spectral theory in general have been developed, the algorithmic aspects of such theories with respect to the the Collatz--Wielandt property have been neglected when concerning nonsquare matrices (and also in other extensions). This paper is the first, to the best of our knowledge, to provide spectral definitions for nonsquare systems that have the same algorithmic power as those made for square systems (in the context of the \PFT).
The extended optimization problem that corresponds to this nonsquare setting is a nonconvex problem (which is also not log-convex), therefore its polynomial solution and characterization are of interest.
%connection between this problem and a complete characterization of ~\PFT. In particular, the Collatz--Wielandt property, which has been neglected when concerning nonsquare matrices (and other applications) is shown arise in nonsquare systems and therefore our generalization is applicable in many areas, similarly to the classical theorem.

%In sum, in this paper we present a novel spectral characterization of nonsquare matrices and show that this characterization is both sufficient and necessary for the strong PF and the Collatz--Wielandt property to hold. The mathematical framework is interesting on its own right, combining technique from graph theory, geometry of polytopes and linear algebra of high degree polynomials. In addition, this formulation turns out to extend many well known applications of the classical \PFT~.

Another way to extend the notion of eigenvalues and eigenvectors of a square
matrix to a nonsquare matrix is via \emph{singular value decompositions (SVD)}
\cite{meyer2000matrix}. Formally, the singular value decomposition of an
$n\times m $ real matrix $M$ is a factorization of the form $M=U\Sigma V^{*}$,
where $U$ is an $m\times m$ real or complex unitary matrix, $\Sigma$ is an
$m \times n$ diagonal matrix with nonnegative reals on the diagonal,
and $V^{*}$ (the conjugate transpose of $V$) is an $n\times n$ real or complex
unitary matrix.
The diagonal entries $\Sigma_{i,i}$ of $\Sigma$ are known as the singular
values of $M$. After performing the product $U \Sigma V^{*}$,
it is clear that the dependence of the singular values of $M$ is linear.
In case all the inputs of $M$ are positive, we can add the absolute
value, and thus the SVD  has a flavor of $L^1$ dependence. In contrast
to the SVD definition,
here we are interested in finding the maximum, so our interpretation has
the flavor of $L^\infty$.
%Since SVD has many application we hope that our extension will be useful to some of the problems.

In a recent paper \cite{Vazirani12}, Vazirani defined the notion of
{\em rational convex programs} as problems that have a rational number as
a solution. Our paper can be considered as an example for
{\em algebraic programming},
since we show that a solution to our problem is an algebraic number.

%%%%%%%%%%%%%%%%%%%%%%%%%%%%%%%%%%%%%%%%%%%%
\section{Preliminaries}
\label{sec:per}
%%%%%%%%%%%%%%%%%%%%%%%%%%%%%%%%%%%%%%%%%%%%
\subsection{Definitions and terminology}
%%%%%%%%%%%%%%%%%%%%%%%%%%%%%%%%%%%%%%%%%%%%
Consider a directed graph $G=(V,E)$. A subset of the vertices $W \subseteq V$ is called a \emph{strongly connected component}
%(or an \emph{SCC} for short)
if $G$ contains a directed path from $v$ to $u$ for every $v, u \in W$. $G$ is said to be \emph{strongly connected} if $V$ is a strongly connected component.
\par Let $A \in \R^{n \times n}$ be a square matrix.
Let $\EigenValue(A)= \{\lambda_1, \ldots, \lambda_k\}$, $k \leq n$,
be the set of real eigenvalues of $A$.
The \emph{characteristic polynomial} of $A$, denoted by $\CP(A,t)$,
is a polynomial whose roots are precisely the eigenvalues of $A$,
$\EigenValue(A)$, and it is given by
\begin{equation}
\label{eq:CP}
\CP(A,t) = \det(t \cdot I -A)
\end{equation}
where $I$ is the $n \times n$ identity matrix.
%It then follows
Note that $\CP(A,t)=0$ iff $t \in \EigenValue(A)$.
The {\em spectral radius} of $A$ is defined as
$\SpectralRatio(A) =
\displaystyle \max\limits_{\lambda \in \EigenValue(A)} |\lambda|.$
The $i^{th}$ element of a vector $\overline{X}$ is given by $X(i)$, and
the $i,j$ entry of a matrix $A$ is denoted $A(i,j)$.
Let $A_{i,0}$ (respectively, $A_{0,i}$) denote the $i$-th row (resp., column)
of $A$. Vector and matrix inequalities are interpreted in the component-wise sense. $A$ is \emph{positive} (respectively, \emph{nonnegative})
%($A>0$)
if all its entries are.
$A$ is \emph{primitive} if there exists a natural number $k$ such that
$A^{k}>0$. $A$ is \emph{irreducible} if for every $i,j$,
there exists a natural $k_{i,j}$ such that $(A^{k_{i,j}})_{i,j} >0.$
An \emph{irreducible} matrix $A$ is \emph{periodic} with period $\Period$ if
$(A^{t})_{ii}=0$ for $t \neq k \cdot \Period$.
%Given $A, B \in R^{n \times n}$, let $C= A \bowtie B^{T}$ correspond to
%$C'= A \cdot B^{T}$ but treating positive values as 1. I.e., $C_{i,j}=1$ iff $C'_{i,j}>0$.

%%%%%%%%%%%%%%%%%%%%%%%%%%%%%%%%%%%%%%%%%%%%
\subsection{Algebraic Preliminaries}
\label{sec:algper}
%%%%%%%%%%%%%%%%%%%%%%%%%%%%%%%%%%%%%%%%%%%%
%%%%%%%%%%%%%%%%%%%%%%%%%%%%%%%%%%%%%%%%%%%%
\paragraph{Generalization of Cramer's rule to homogeneous linear systems.}
%%%%%%%%%%%%%%%%%%%%%%%%%%%%%%%%%%%%%%%%%%%%
Let $A_{i,0}$ (respectively, $A_{0,i}$) denote the $i$-th row (resp., column)
of $A$. Let $A_{-(i,j)}$ denote the matrix that results from $A$ by removing
the $i$-th row and the $j$-th column. Similarly, $A_{-(i,0)}$ and $A_{-(0,j)}$
denote the matrix after removing the $i$-th row (respectively, $j$-th column)
from $A$. Let $\widetilde{A}_{i}=\left(A(1,i), \ldots, A(n-1,i) \right)^T$,
i.e., the $i$-th column of $A$ without the last element $A(n,i)$.
For $\overline{X}=(X(1), \ldots, X(n)) \in \R^{n}$, denote
$\overline{X}_{i}=(X(1), \ldots, X(i)) \in \R^{i}$.

We make use of the following extension of Cramer's rule to homogeneous square linear systems.
\begin{claim}
\label{cl:cramer_square}
Let $A \cdot \overline{X} = \overline{0}$ such that $A_{-(n,n)}$ is invertible.
Then,
\begin{description}
\item{(a)}
%\begin{equation}
%\label{eq:square_cramer_i}
$\displaystyle
X(i) ~=~ (-1)^{n-i} \cdot X(n) \cdot \frac{\det(A_{-(n,i)})}{\det(A_{-(n,n)})}~.$
%\end{equation}
\item{(b)}
%\begin{equation}
%\label{eq:square_cramer_n}
$\displaystyle
X(n) \cdot \frac{\det(A)}{\det(A_{-(n,n)})}=0~.$
%\end{equation}
\end{description}
\end{claim}
\Proof
Since $A \cdot \overline{X} = \overline{0}$, it follows that
$A_{-(n,n)} \cdot \overline{X}_{n-1}=-X(n) \cdot \widetilde{A}_{n}$.
As $A_{-(n,n)}$ is invertible, we can apply Cramer's rule to express $X(i)$.
Let $M_{i}=[\widetilde{A}_{1}, \ldots, \widetilde{A}_{i-1},\widetilde{A}_{n},
\widetilde{A}_{i+1}, \ldots, \widetilde{A}_{n-1}] \in \R^{(n-1) \times (n-1)}$,
for $i>1$ and
$M_{1}=[\widetilde{A}_{n}, \widetilde{A}_{2},\ldots, \widetilde{A}_{n-1}]$.
By Cramer's rule, it then follows that
$X(i)=-X(n) \cdot \det(M_{i}) / \det(A_{-(n,n)})$.
We next claim that $\det(M_{i})=(-1)^{n-1-i} \cdot \det(A_{-(n,i)})$.
To see this, note that $A_{-(n,i)}$ and $M_{i}$ are composed of the same set of
columns up to order. In particular, $M_{i}$ can be transformed to
$A_{-(n,i)}=[\widetilde{A}_{1}, \ldots, \widetilde{A}_{i-1}, \widetilde{A}_{i+1},
\ldots, \widetilde{A}_{n-1},\widetilde{A}_{n}]$
by a sequence of $n-1-i$ swaps of consecutive columns starting from the
$i$-th column of $M_{i}$.  It therefore follows that
$X(i)=(-1)^{n-1-i}  \cdot  -(1) \cdot X(n) \cdot
\frac{\det(A_{-(n,i)})}{\det(A_{-(n,n)})}$
establishing part (a) of the claim.
%Eq. (\ref{eq:square_cramer_i}).
We continue with part (b).
%Eq. (\ref{eq:square_cramer_n}).
Since $A \cdot \overline{X} = \overline{0}$, it follows that
$A_{(n,0)} \cdot \overline{X}=0$ or that
\begin{eqnarray*}
A_{n,0} ~\cdot~ \overline{X} &=&\sum_{i=1}^{n} A(n,i) \cdot X(i)\nonumber
\\&=&
X(n) \cdot \sum_{i=1}^{n-1} \left( (-1)^{n-i} \cdot A(n,i)  \cdot
\frac{\det(A_{-(n,i)})}{\det(A_{-(n,n)})} \right)+A(n,n) \cdot X(n) \nonumber
\\&=&
X(n) \cdot \frac{ \sum_{i=1}^{n-1} \left((-1)^{n-i} \cdot A(n,i) \cdot
\det(A_{-(n,i)}) \right)+ A(n,n) \cdot (-1)^{2n} \det(A_{-(n,n)})}{\det(A_{-(n,n)})}
\\&=&
X(n) \cdot \frac{\det(A)}{\det(A_{-(n,n)})}=0~. ~~~\inQED
\nonumber
\end{eqnarray*}
We now turn to a nonsquare matrix $A \in \R^{n \times (n+1)}$.
The matrix $B=B(A) = [\widetilde{A}_{1}, \ldots, \widetilde{A}_{n-1}] \in
\R^{(n-1) \times (n-1)}$ corresponds to the upper left  $(n-1) \times (n-1)$
square matrix of $A$. Let $C^{1}=[A_{1}, \ldots, A_{n}]$ i.e.,
$C^{1}=A_{-(0,n+1)}$ and $C^{2}=A_{-(0,n)}$. Note that
$C^{1}, C^{2} \subseteq  \R^{n \times n}$, i.e., both are square matrices.

\begin{claim}
\label{cl:cramer_non_square}
Let $A \cdot \overline{X}=\overline{0}$  and $B=B(A)$ is invertible. Then,
\begin{description}
\item{(a)}
%\begin{equation}
%\label{eq:non_square_cramer_i}
$\displaystyle
X(i) ~=~ (-1)^{n-i}  \cdot \left( \frac{\det(C^{1}_{-(n,i)})}{\det(B)} \cdot
X(n) +\frac{\det \left(C^{2}_{-(n,i)} \right)}{\det \left(B \right)} \cdot X(n+1) \right) ~,$
%\end{equation}
\item{(b)}
%\begin{equation}
%\label{eq:non_square_cramer_n}
$\displaystyle
X(n) \cdot \frac{\det \left(C^1 \right)}{\det(B)} ~=~
-X(n+1) \cdot \frac{\det \left(C^2 \right)}{\det(B)}~.$
%\end{equation}
\end{description}
\end{claim}
\Proof
Since $A \cdot \overline{X}= \overline{0}$, it follows that
$B \cdot \overline{X}_{n-1}= - \left(X(n) \cdot \widetilde{A}_{n}   + X(n+1)
\cdot \widetilde{A}_{n+1} \right)$.
As $B$ is invertible we can apply Cramer's rule to express $x_{i}$.
Let $M_{i}=[\widetilde{A}_{1}, \ldots, \widetilde{A}_{i-1},x_{n} \cdot
\widetilde{A}_{n}+x_{n+1} \cdot \widetilde{A}_{n+1}, \widetilde{A}_{i+1},
\ldots, \widetilde{A}_{n-1}] \in \R^{n \times n}$.
Let $M_{i}^{1}=[\widetilde{A}_{1}, \ldots, \widetilde{A}_{i-1},
\widetilde{A}_{n} , \widetilde{A}_{i+1}, \ldots, \widetilde{A}_{n-1}]$ and
\\
$M_{i}^{2}=[\widetilde{A}_{1}, \ldots, \widetilde{A}_{i-1}, \widetilde{A}_{n+1} ,
\widetilde{A}_{i+1}, \ldots, \widetilde{A}_{n-1}]$.
By the properties of the determinant function, it follows, that
$$X(i)=X(n) \cdot \frac{\det\left(M_{i}^{1}\right)}{\det\left(B\right)}+
X(n+1) \cdot \frac{\det\left(M_{i}^{2}\right)}{\det\left(B\right)}.$$
We now turn to see the connection between $\det(M_{i}^{1})$ and
$\det(C_{-(n,i)}^{1})$. Note that $M_{i}^{1}$ and $C_{-(n,i)}^{1}$ correspond to
the same columns up to order. Specifically, we can now employ
the same argument of Claim  \ref{cl:cramer_square} and show that
$\det(M_{i}^{1})=(-1)^{n-1-i} \cdot \det(C_{-(n,i)}^{1})$
(informally, the square matrix of Claim \ref{cl:cramer_square} is replaced by
a ``combination" of $C_{1}$ and $C_{2}$).  In a similar way, one can show that
$\det(M_{i}^{2})=(-1)^{n-1-i} \cdot \det(C_{-(n,i)}^{2})$.
We now turn to prove part (b) of the claim.
%Eq. (\ref{eq:non_square_cramer_n}).
Since $A_{n,0} ~\cdot~ \overline{X}$, by part (a),
%Eq. (\ref{eq:non_square_cramer_i}),
we get that
\begin{eqnarray*}
A_{n,0} ~\cdot~ \overline{X}&=&
\sum_{i=1}^{n-1} A(n,i) \cdot X(i) +A(n,n) \cdot X(n)+ A(n,n+1) \cdot X(n+1)
\nonumber
\\&=&
X(n) \cdot \left(\sum_{i=1}^{n-1}  (-1)^{n-i} \cdot A(n,i)  \cdot
\frac{\det \left(C^{1}_{-(n,i)} \right)}{\det(B)} +A(n,n) \right)
\nonumber
\\&&
+ X(n+1) \cdot  \sum_{i=1}^{n-1} \left( (-1)^{n-i} \cdot A(n,i)  \cdot
\frac{\det \left(C^{2}_{-(n,i)} \right)}{\det(B)} +A(n,n+1)\right)
\nonumber
\\&=&
X(n) \cdot \frac{\sum_{i=1}^{n-1}(-1)^{n-i} \cdot A(n,i)  \cdot
\det\left(C^{1}_{-(n,i)} \right) + (-1)^{2n} \cdot A(n,n) \cdot \det(B)}{\det(B)}
\nonumber
\\&&
+ X(n+1) \cdot \frac{\sum_{i=1}^{n-1} (-1)^{n-i} \cdot A(n,i)  \cdot
\det \left(C^{2}_{-(n,i)} \right) + (-1)^{2n} \cdot A(n,n+1) \cdot \det(B)}{\det(B)}
\nonumber
\\&=&
X(n) \cdot \frac{\det(C^{1})}{\det(B)}+X(n+1) \cdot
\frac{\det(C^{2})}{\det(B)}=0~.
\nonumber
\end{eqnarray*}
The claim follows.
\QED
%%%%%%%%%%%%%%%%%%%%%%%%%%%%%%%%%%%%%%%%%%%%
\paragraph{Separation theorem for nonsymmetric matrices.}
%%%%%%%%%%%%%%%%%%%%%%%%%%%%%%%%%%%%%%%%%%%%
We make use of the following fact due to Hall and Porsching \cite{HallInterlacing}, which is an extension of the Cauchy Interlacing Theorem for symmetric matrices.
\begin{lemma}[\cite{HallInterlacing}]
\label{lem:sep_thm}
Let $A$ be a nonegative matrix with eigenvalues $\EigenValue(A)=\{\lambda_i(A) \mid i \in \{1, \ldots, n\}\}$. Let $A_i$ be the $i^{th}$ principle $(n-1) \times (n-1)$ minor of $A$, with eigenvalues $\lambda_j(A_i)$, $j \in \{1, \ldots, n-1\}$.
If $\lambda_p(A)$ is any real eigenvalue of $A$ different from $\lambda_1[A]$, then
$$\lambda_1(A) \leq \lambda_1(A_i) \leq \lambda_p(A)$$
for every $i \in \{1, \ldots, n\}$, with strict inequality on the left if $A$ is irreducible.
\end{lemma}
%%%%%%%%%%%%%%%%%%%%%%%%%%%%%%%%%%%%%%%%%%%%
\subsection{\PFT~for square nonnegative irreducible matrices}
%%%%%%%%%%%%%%%%%%%%%%%%%%%%%%%%%%%%%%%%%%%%
The \PFT~states the following.
\begin{theorem} [PF Theorem, \cite{PF_Frobenius,PF_Perron}]
\label{thm:pf_full}
Let $A \in \R_{\geq 0}^{n \times n}$ be a nonnegative irreducible matrix with
spectral ratio $\SpectralRatio(A)$. Then $\max \EigenValue(A)>0$.
There exists an eigenvalue $r \in \EigenValue(A)$ such that
$r=\SpectralRatio(A)$, called the
\emph{Perron--Frobenius (PF) root} of $A$.
The algebraic multiplicity of $~\PFEigenValue$ is one.
There exists an eigenvector $\overline{X}>0$ such that
$A \cdot \overline{X}=\PFEigenValue \cdot \overline{X}$.
The unique normalized vector $\PFEigenVector$ defined by
$A \cdot \PFEigenVector=\PFEigenValue \cdot \PFEigenVector$
and $||\PFEigenVector||_{1}=1$
is called the \emph{Perron--Frobenius (PF) vector}.
There are no nonnegative eigenvectors for $A$  with $r$ except for positive multiples
of $\PFEigenVector$. If $A$ is a nonnegative irreducible periodic matrix
with period $\Period$, then $A$ has exactly $\Period$ eigenvalues,
$\lambda_j= \SpectralRatio(A) \cdot \exp^{2 \pi i \cdot j/\Period}$ for
$j =1,2, \ldots, \Period,$
and all other eigenvalues of $A$ are of strictly smaller magnitude
than $\SpectralRatio(A)$.
\end{theorem}

%%%%%%%%%%%%%%%%%%%%%%%%%%%%%%%%%%%%%%%%%%%%
\paragraph{Collatz--Wielandt characterization (the min-max ratio).}
%%%%%%%%%%%%%%%%%%%%%%%%%%%%%%%%%%%%%%%%%%%%
Collatz and Wielandt \cite{PF_Collatz, PF_Wielandt} established the following
formula for the \PFR, also known as the min-max ratio characterization.
\begin{lemma} \cite{PF_Collatz, PF_Wielandt}
[Collatz--Wielandt]
\label{lem:Collatz-Wielandt}
$\PFEigenValue=\min_{\overline{X} \in \mathcal{N}} \{\mathfrak{f}(\overline{X})\}$
~where~
$$\mathfrak{f}(\overline{X})= \max\limits_{1 \leq i \leq n, X(i)\neq \overline{0}}
\left \{\frac{(A \cdot \overline{X})_{i}}{X(i)} \right \} \mbox{~~and~~}
\mathcal{N}=\{\overline{X} \geq \overline{0},||\overline{X}||_{1}=1\}.$$
\end{lemma}
Alternatively, this can  be written as the following optimization problem.
%by taking $r=1/\beta$.
\begin{equation}\label{LP:Stand_Perron}
\mbox{maximize} ~~ \beta \text{~~~subject to:~~~} \displaystyle A \cdot \overline{X} \leq 1/\beta
\cdot \overline{X},~~~ \displaystyle ||\overline{X}||_{1}=1,~~~
\displaystyle \overline{X} \geq \overline{0}.
\end{equation}
Let $\beta^{*}$ be the optimal solution of Program (\ref{LP:Stand_Perron}) and let
$\overline{X}^{*}$ be the corresponding optimal vector.
Using the representation of Program (\ref{LP:Stand_Perron}),
Lemma \ref{lem:Collatz-Wielandt} translates into the following.
\begin{theorem}
\label{thm:pf}
The optimum solution of \eqref{LP:Stand_Perron} satisfies $\beta^{*}=1/\PFEigenValue$, where $\PFEigenValue \in \R_{>0}$ is the maximal
eigenvalue of $A$ and $\overline{X}^{*}$ is given by eigenvector
$\PFEigenVector$ corresponding for $\PFEigenValue$.
Hence for $\beta^{*}$, the $n$ constraints given by
$A \cdot \overline{X}^{*} \leq 1/\beta^{*} \cdot \overline{X}^{*}$ of
Program (\ref{LP:Stand_Perron}) hold with equality.
\end{theorem}
%\begin{corollary}
%\label{cor:pf}
%At the optimum value $\beta^{*}$, the set of $n$ constraints given by $A \cdot %\overline{X} \leq 1/\beta \cdot \overline{X}$ of Eq. (\ref{LP:Stand_Perron}) hold %with equality.
%\end{corollary}
This can be interpreted as follows. Consider the ratio
$Y(i)= (A\cdot  \overline{X})_{i}/X(i)$, viewed as the ``repression factor"
for entity $i$. The task is to find the input vector $\overline{X}$
that minimizes the maximum repression factor over all $i$,
thus achieving balanced growth.
In the same manner, one can characterize the $\max$-$\min$ ratio.
Again, the optimal value (resp., point) corresponds to the PF eigenvalue
(resp., eigenvector) of $A$.
%\par
In summary, when taking
$\overline{X}$ to be
the PF eigenvector, $\PFEigenVector$, and $\beta^{*}=1/\PFEigenValue$,
all repression factors are equal, and optimize the $\max$-$\min$ and
$\min$-$\max$ ratios.
%%%%%%%%%%%%%%%%%%%%%%%%%%%%%%%%%%%%%%%%%%%%
\section{A generalized \PFT~for nonsquare systems}
%%%%%%%%%%%%%%%%%%%%%%%%%%%%%%%%%%%%%%%%%%%%
\subsection{The Problem}
\paragraph{System definitions.}
%%%%%%%%%%%%%%%%%%%%%%%%%%%%%%%%%%%%%%%%%%%%
Our framework consists of a set
$\EntitySet=\{\Entity_{1}, \ldots, \Entity_{n}\}$
of entities whose growth is regulated by a set of \emph{affectors}
$\Affectors=\{\Affectors_1, \Affectors_2, \ldots, \Affectors_m\}$,
for some $m \geq n$.
As part of the solution, each affector is set to be either {\em passive} or
{\em active}.
If an affector $\Affectors_j$ is set to be active, then it affects
each entity $\Entity_i$, by  either increasing or decreasing it by a certain
amount $g(i,j)$, which is specified as part of the input.
If  $g(i,j) >0$ (resp., $g(i,j) < 0$), then $\Affectors_j$ is referred to as a
\emph{supporter} (resp., \emph{repressor}) of $\Entity_i$.
%The amount by which affector $\Affectors_j$ affects $\Entity_i$
%(when it is active), denoted $g(i,j)$ is specified as part of the input.
%If  $g(i,j) >0$, then $\Affectors_j$ is referred to as a \emph{supporter} of  $\Entity_i$. If $g(i,j) < 0$, then $\Affectors_j$ is referred to as a \emph{repressor} of $\Entity_i$.
%
For clarity we may write $g(\Entity_i,\Affectors_j)$ for $g(i,j)$.
%
%For ease of algebraic representation,
The affector-entity relation is described by two matrices,
the \emph{supporters gain}  matrix $\SupportersMatrix \in \R^{n \times m}$ and
the \emph{repressors gain} matrix $\RepressorsMatrix \in \R^{n \times  m}$,
given by
\begin{equation*}
\SupportersMatrix(i,j) =
\begin{cases}
g(\Entity_i,\Affectors_j), & \text{if $g(\Entity_i,\Affectors_j) >0$;}\\
0, & \text{otherwise.}
\end{cases}
\end{equation*}
\begin{equation*}
\RepressorsMatrix(i,j) =
\begin{cases}
-g(\Entity_i,\Affectors_j),  & \text{if $g(\Entity_i,\Affectors_j) <0$;}\\
0, & \text{otherwise.}
\end{cases}
\end{equation*}
Again, for clarity we may write $\RepressorsMatrix(\Entity_i,\Affectors_j)$
for $\RepressorsMatrix(i,j)$, and similarly for $\SupportersMatrix$.
%$\RepressorsMatrix(\Entity_i,\Affectors_j)=|g(\Entity_i,\Affectors_j)|$.

We can now formally define a {\em system} as
$\System=\langle \SupportersMatrix, \RepressorsMatrix \rangle$, where
$\SupportersMatrix, \RepressorsMatrix \in \R^{m \times n}_{\geq 0}$,
$n=|\EntitySet|$ and $m=|\Affectors|$.
We denote the supporter (resp., repressor) set of $\Entity_i$ by
\begin{eqnarray*}
\Supporters_{i}(\System) &=& \{\Affectors_j \mid
\SupportersMatrix(\Entity_i,\Affectors_j)>0\},
\\
\Repressors_{i}(\System) &=& \{\Affectors_j \mid
\RepressorsMatrix(\Entity_i,\Affectors_j)>0\}.
\end{eqnarray*}
When $\System$ is clear from the context, we may omit it and simply write
$\Supporters_{i}$ and $\Repressors_{i}$.
%%% reducing to the standard PF
Throughout, we restrict attention to systems in which $|\Supporters_i|\geq 1$
for every $\Entity_i \in \EntitySet$.
We classify the systems into three types:
\begin{description}
\item{(a)}
$\SquareSystemFamily=\{\System ~\mid~ m \leq n, |\Supporters_i|=1 \text{~for every} ~ \Entity_i \in \EntitySet\}$
is the family of \emph{Square Systems}.
\item{(b)}
$\WeakSystemFamily=\{\System \mid m \leq n+1, \exists j
\text{~s.t~} |\Supporters_j|=2 \text{~and~}|\Supporters_i|=1
\text{~for every~}\Entity_i \in \EntitySet \setminus \{\Entity_j\} \}$
is the family of \emph{Weakly Square Systems}, and
\item{(c)}
$\SystemFamily=\{\System \mid m>n+1\}$ is the family of
\emph{Nonsquare Systems}.
\end{description}

%%%%%%%%%%%%%%%%%%%%%%%%%%%%%%%%%%%%%%%%%%%%
\paragraph{The generalized PF optimization problem.}
%%%%%%%%%%%%%%%%%%%%%%%%%%%%%%%%%%%%%%%%%%%%
Consider a set of $n$ entities and gain matrices
$\SupportersMatrix,\RepressorsMatrix \in \R^{n \times m}$, for $m \geq n$.
The main application of the generalized \PFT~is the following optimization
problem, which is an extension of Program (\ref{LP:Stand_Perron}).
\begin{align}
\label{LP:Ext_Perron}
  \mbox{maximize~~} ~& \beta \mbox{~~subject to:~~}
  \\
   & \displaystyle \RepressorsMatrix \cdot \overline{X} ~\leq~
       1/\beta \cdot \SupportersMatrix \cdot \overline{X} ~,&
  \label{eq:SR}  \\
  & \displaystyle \overline{X} \geq \overline{0}~, &
  \label{eq:Ineq}\\
  & \displaystyle ||\overline{X}||_{1}=1~. &
   \nonumber
\end{align}
We begin with a simple observation.
An affector $\Affectors_j$ is \emph{redundant} if
$\SupportersMatrix(\Entity_i,\Affectors_j)=0$ for every $i$.
\begin{observation}
\label{obs:only_positive}
If $\Affectors_j$  is \emph{redundant}, then $X(j)=0$ in any optimal solution
$\overline{X}$.
\end{observation}

%Consequently, we restrict attention to systems without redundant affectors.
In view of Obs. \ref{obs:only_positive}, we may hereafter restrict attention
to the case where there are no redundant affectors in the system,
as any redundant affector $\Affectors_j$ can be removed and simply assigned
$X(j)=0$.

We now proceed with some definitions.
Let $X(\Affectors)$ denote the value of $\Affectors$ in $\overline{X}$, i.e., $X(\Affectors)=X(k)$ where the $k'$th entry in $\overline{X}$ corresponds to $\Affectors$.
An affector $\Affectors$ is \emph{active} in a solution $\overline{X}$ if $X(\Affectors)>0$.
Denote the set of affectors taken to be active in a solution $\overline{X}$
by $NZ(\overline{X})=\{\Affectors_j \mid X(\Affectors_j)> 0\}$.
Let $\beta^{*}(\System)$ denote the optimal value of Program
(\ref{LP:Ext_Perron}), i.e., the maximal positive value $\beta$ for which there exists
a nonnegative, nonzero vector $\overline{X}$ satisfying the constraints of
Program (\ref{LP:Ext_Perron}).
When the system $\System$ is clear from the context
we may omit it and simply write $\beta^*$.
A vector $\overline{X}_{\widetilde{\beta}}$ is \emph{feasible} for
$\widetilde{\beta} \in (0,\beta^{*}]$ if it satisfies all the constraints of
Program (\ref{LP:Ext_Perron}) with $\beta=\widetilde{\beta}$.
A vector $\overline{X}^{*}$ is \emph{optimal} for $\System$ if it is feasible
for $\beta^{*}(\System)$, i.e., $\overline{X}^{*}=\overline{X}_{\beta^{*}}$.
The system $\System$ is \emph{feasible} for $\beta$ if
$\beta\leq \beta^{*}(\System)$, i.e., there exists a feasible
$\overline{X}_{\beta}$ solution for Program (\ref{LP:Ext_Perron}).
%The $\System$ is \emph{feasible} if is feasible for some $\beta>0$.
\par For vector $\overline{X}$, the \emph{total repression} on $\Entity_{i}$ in
$\System$ for a given $\overline{X}$ is
$\TotR(\overline{X}, \System)_{i}=(\RepressorsMatrix \cdot \overline{X})_{i}$.
Analogously, the \emph{total support} for $\Entity_{i}$ is
$\TotS(\overline{X}, \System)_{i}=(\SupportersMatrix \cdot \overline{X})_{i}$.
We now have the following alternative formulation for the constraints of Eq. (\ref{eq:SR}), stated individually for each entity $\Entity_i$.
\begin{equation}
\label{eq:SR_ind}
\TotR(\overline{X}, \System)_{i} ~\leq~
1/\beta \cdot \TotS(\overline{X}, \System)_{i}~\text{~for every~} i \in \{1, \ldots, n\}~.
\end{equation}
\begin{fact}
\label{fc:feasible_tots_totr}
Eq. (\ref{eq:SR}) holds iff Eq. (\ref{eq:SR_ind}) holds.
\end{fact}
We classify the $m+n$ linear
inequality constraints of Program (\ref{LP:Ext_Perron}) into two types of constraints:
\begin{description}
\item{(1)}
SR (Support-Repression) constraints:
the $n$ constraints of Eq. (\ref{eq:SR}) or alternatively of Eq. (\ref{eq:SR_ind}).
\item{(2)}
Nonnegativity constraints:
the $m$ constraints of Eq. (\ref{eq:Ineq}).
\end{description}
When $\System$ is clear from context, we may omit it and simply write
$\TotR(\overline{X})_{i}$ and $\TotS(\overline{X})_{i}$.
As a direct application of the generalized PF Theorem, there is an exact
polynomial time algorithm for solving Program (\ref{LP:Ext_Perron})
for irreducible systems, as defined next.

%%%%%%%%%%%%%%%%%%%%%%%%%%%%%%%%%%%%%%%%%%%%
\subsection{Irreducibility of PF systems}
%%%%%%%%%%%%%%%%%%%%%%%%%%%%%%%%%%%%%%%%%%%%
\paragraph{Irreducibility of square systems.}
%%%%%%%%%%%%%%%%%%%%%%%%%%%%%%%%%%%%%%%%%%%%
%\begin{definition}
%\label{def:irreducible_square}
A square system $\System=\langle \SupportersMatrix,\RepressorsMatrix\rangle\in \SquareSystemFamily$ is {\em irreducible} iff
(a) $\SupportersMatrix$ is nonsingular and
(b) $\RepressorsMatrix$ is irreducible.
%\end{definition}
Given an irreducible square $\System$,
%the following formulation is convenient for our later arguments.
%transforming from Program (\ref{LP:Ext_Perron}) to Program
%(\ref{LP:Stand_Perron}):
let
\begin{equation*}
Z(\System) ~=~ \left(\SupportersMatrix \right)^{-1} \cdot \RepressorsMatrix~.
\end{equation*}
Note the following two observations.
\begin{observation}
\label{cl:irreducible_supporter}
(a) If $\SupportersMatrix$ is nonsingular, then
$\Supporters_i \cap \Supporters_j = \emptyset$.\\
(b) If $\System$ is an irreducible system,
then $Z(\System)$ is an irreducible matrix as well.
\end{observation}
\Proof
Consider part (a). Since $\System$ is square, $|\Supporters_i|=1$ for every $i$. Combining with the fact that $\SupportersMatrix$ is nonsingular, it holds that $\SupportersMatrix$ is equivalent (up to column alternations) to a diagonal matrix with a fully positive diagonal, hence $\Supporters_i \cap \Supporters_j= \emptyset$. Part (b) follows by definition.
\QED

%Given an irreducible square system $\System$,
%the following formulation is convenient for our later arguments.
%transforming from Program (\ref{LP:Ext_Perron}) to Program
%(\ref{LP:Stand_Perron}):
%let $Z(\System)=\left(\SupportersMatrix \right)^{-1} \cdot \RepressorsMatrix~$.
%\begin{equation}
%\label{eq:z_matrix}
%Z(\System)=\left(\SupportersMatrix \right)^{-1} \cdot \RepressorsMatrix~.
%\end{equation}
%
%Note that
%\begin{lemma}
%\label{obs:reducing_to_stand_PF}
%If $\System$ is an irreducible system,
%then $Z(\System)$ is an irreducible matrix as well.
%\end{lemma}
%
Throughout, when considering square systems, it is convenient to assume that
the entities and affectors are ordered in such a way that $\SupportersMatrix$
is a diagonal matrix, i.e., in $\SupportersMatrix$
(as well as in $\RepressorsMatrix$) the $i^{th}$ column corresponds to
$\Affectors_{k} \in \Supporters_i$, the unique supporter of $\Entity_i$.

%
%%%%%%%%%%%%%%%%%%%%%%%%%%%%%%%%%%%%%%%%%%%%%%%%%%%%%%%%%%%%%%%%%%%%%%
\paragraph{Selection matrices.}
%%%%%%%%%%%%%%%%%%%%%%%%%%%%%%%%%%%%%%%%%%%%
To define a notion of irreducibility for a nonsquare system
$\System \notin \SquareSystemFamily$, we first present the notion of a
{\em selection matrix}.
A selection matrix $\FilterMatrix \in \{0,1\}^{m \times n}$ is \emph{legal}
for $\System$ iff for every entity $\Entity_i \in \EntitySet$ there exists
exactly one supporter $\Affectors_j \in \Supporters_i$ such that
$\FilterMatrix(j,i)=1$.
Such a matrix $\FilterMatrix$ can be thought of as representing a selection
performed on $\Supporters_i$ by each entity $\Entity_i$, picking exactly one
of its supporters. Let $\System(\FilterMatrix)$ be the square system corresponding to the legal
selection matrix $\FilterMatrix$, namely,
$\System(\FilterMatrix)=\displaystyle \langle \SupportersMatrix \cdot
\FilterMatrix, \RepressorsMatrix \cdot \FilterMatrix\rangle.$ In the resulting system there are $m' \leq n$ non-redundant affectors.
Since
%every entity selects exactly one supporter, and as
redundant affectors can be discarded from the system (by Obs. \ref{obs:only_positive}), it follows that
the number of active affectors becomes at most the number of entities,
resulting in a square system.
Denote the family of legal selection matrices,
capturing the ensemble of all square systems hidden in $\System$, by
\begin{equation}
\label{eq:FilterMatrixFamily}
\FilterMatrixFamily(\System) ~=~
\{\FilterMatrix \mid  \FilterMatrix \text{~is legal for~} \System \}.
\end{equation}
When $\System$ is clear from the context, we simply write $\FilterMatrixFamily$.
%Let $\System(\FilterMatrix)=\displaystyle \langle \SupportersMatrix \cdot
%\FilterMatrix, \RepressorsMatrix \cdot \FilterMatrix\rangle$ be the square
%system corresponding to the legal selection matrix $\FilterMatrix$.
Let $\overline{X}_{\beta} \in \R^{n}$ be a solution for the square system $\System(\FilterMatrix)$ for some $\FilterMatrix$. The \emph{natural extension} of $\overline{X}_{\beta} \in \R^{n}$ into a solution $\overline{X}^{m}_{\beta} \in \R^m$ of the original system $\System$ is defined by letting $X^{m}_{\beta}(\Affectors_k)=X_{\beta}(\Affectors_k)$
if $\sum_{\Entity_i \in \EntitySet}\FilterMatrix(\Affectors_k, \Entity_i)>0$
and $X^{m}_{\beta}(\Affectors_k)=0$ otherwise.
\begin{observation}
\label{obs:filter_to_square}
(a) $\System(\FilterMatrix) \in \SquareSystemFamily$
for every $\FilterMatrix \in \FilterMatrixFamily$.\\
(b) For every solution $\overline{X}_{\beta} \in \R^{n}$  for system $\System(\FilterMatrix)$, for some matrix $\FilterMatrix \in \FilterMatrixFamily$,
its natural extension $\overline{X}_{\beta}^{m}$ is a feasible solution for the original $\System$.\\
(c) $\beta^{*}(\System) \geq \beta^{*}(\System(\FilterMatrix))$ for every selection matrix $\FilterMatrix \in \FilterMatrixFamily$.
\end{observation}
%\Proof
%Part (a) follows by definition.
%Let $\SelectionVec$ be a complete
%selection vector, $|\SelectionVec|=n$, corresponding to the square system
%$\System(\FilterMatrix)$.
%The proof of Part (b) and (c) are immediate.
%\QED
%
%%%%%%%%%%%%%%%%%%%%%%%%%%%%%%%%%%%%%%%%%%%%
\paragraph{Irreducibility of nonsquare systems.}
%%%%%%%%%%%%%%%%%%%%%%%%%%%%%%%%%%%%%%%%%%%%
We are now ready to define the notion of irreducibility for nonsquare systems, as follows.
%\begin{definition}
%\label{def:irred_non_square}
A nonsquare system $\System$ is \emph{irreducible} iff $\System(\FilterMatrix)$ is irreducible for every selection matrix $\FilterMatrix \in \FilterMatrixFamily$.
Note that this condition is the ``minimal'' \emph{necessary} condition
for our theorem to hold, as explained next.
%in Section \ref{sec:irredsystems}.
% for the following reason.
%\def\AppendExplainIrred{
Our theorem states that the optimum solution for the nonsquare system is the optimum solution for the best \emph{embedded} square system.  It is easy to see that for any nonsquare system $\System=\langle \SupportersMatrix, \RepressorsMatrix \rangle$,
one can increase or decrease any entry $g(i,j)$ in the matrices, while maintaining the sign of each entry in the matrices, such that a particular
selection matrix $\FilterMatrix^{*} \in \FilterMatrixFamily$ would correspond to the optimal square system.  With an optimal embedded square system at hand, which is also guaranteed to be irreducible (by the definition of irreducible nonsquare systems), our theorem can then apply the traditional \PFT, where a spectral characterization for the solution of Program (\ref{LP:Stand_Perron}) exists. Note that irreducibility is a \emph{structural} property of the system, in the sense that it does not depend on the exact gain values, but rather on the sign of the gains, i.e., to determine irreducibility, it is sufficient to observe the binary matrices $\SupportersMatrix_{B}, \RepressorsMatrix_{B}$, treating $g(i,j) \neq 0$ as $1$. On the other hand, deciding which of the embedded square systems has the maximal eigenvalue (and hence is optimal), depends on the \emph{precise} values of the entries of these matrices. It is therefore necessary that the structural property of irreducibility would hold for any specification of gain values (while maintaining the binary representation of $\SupportersMatrix_{B}, \RepressorsMatrix_{B}$).
Indeed, consider a reducible nonsquare system, for which there exists an embedded square system $\System(\FilterMatrix)$ that is reducible. It is not hard to see that there exists a specification of gain values that would render this square system $\System(\FilterMatrix)$ optimal (i.e., with the maximal eigenvalue among all other embedded square systems). But since $\System(\FilterMatrix)$ is reducible, the \PFT\ cannot be applied, and in particular, the corresponding eigenvector is no longer guaranteed to be \emph{positive}.
%} %\AppendExplainIrred

%\inline Discussion on the notion of  irreducible nonsquare systems:
%We first argue that our condition for the irreducibility for nonsquare systems
%is the ``minimal'' \emph{necessary} condition for our theorem to hold,
%for the following reason.
%
%\AppendExplainIrred
%\end{definition}
%Define
%$$\HiddenSquareFamily(\System)=\{\System(\FilterMatrix) \mid \FilterMatrix \in %\FilterMatrixFamily\}.$$
%It then follows that
%\begin{corollary}
%\label{cor:irreducible_hidden_square}
%system is irreducible iff every hidden square system $\System^{s} \in %\HiddenSquareFamily(\System)$ is irreducible.
%\end{corollary}

\begin{claim}
\label{cor:distinct}
In an irreducible system $\System$,
$\Supporters_i \cap \Supporters_j=\emptyset$ for every $\Entity_i, \Entity_j$.
\end{claim}
\Proof
Assume, toward contradiction, that there exists some affector
$\Affectors_k \in \Supporters_i \cap \Supporters_j$, and consider a selection
matrix $\FilterMatrix$ for which $\FilterMatrix(k,i)=1$ and
$\FilterMatrix(k,j)=1$. It then follows
by Obs. \ref{cl:irreducible_supporter}(a)
that $\SupportersMatrix \cdot \FilterMatrix$ is singular.
But the irreducibility of $\System$ implies that
$\SupportersMatrix \cdot \FilterMatrix$ is nonsingular for every
$\FilterMatrix \in \FilterMatrixFamily$; contradiction.
\QED
%%%%%%%%%%%%%%%%%%%%%%%%%%%%%%%%%%%%%%%%%%%%
\paragraph{Constraint graphs: a graph theoretic representation.}
%\label{subsec:constraints_graph}
%\paragraph{The constraint graph.}
%%%%%%%%%%%%%%%%%%%%%%%%%%%%%%%%%%%%%%%%%%%%
We now provide a graph theoretic characterization of irreducible systems $\System$.
Let $\ConstraintsGraph_{\System}(V,E)$ be the directed \emph{constraint graph} for the system $\System$, defined as follows:
$V= \EntitySet$, and the rule for a directed edge $e_{i,j}$ from $\Entity_i$ to $\Entity_j$ is
\begin{equation}
\label{eq:cg_condition}
e_{i,j} \in E ~~~~\mbox{~iff~}~~~~
\Supporters_{i} \cap \Repressors_j \neq \emptyset.
\end{equation}
Note that it is possible that
$\ConstraintsGraph_{\System} \nsubseteq \ConstraintsGraph_{\System(\FilterMatrix)}$ for some $\FilterMatrix \in \FilterMatrixFamily$.
A graph $\ConstraintsGraph_{\System}(V,E)$ is \emph{robustly strongly connected}
if $\ConstraintsGraph_{\System(\FilterMatrix)}(V,E)$ is strongly connected
for every $\FilterMatrix \in \FilterMatrixFamily$.
\begin{observation}
\label{obs:reducible_graph_connected}
Let $\System$ be an irreducible system.
\begin{description}
\item{(a)}
If $\System$ is square, then
$\ConstraintsGraph_{\System}(V,E)$ is strongly connected.
\item{(b)}
If $\System$ is nonsquare, then $\ConstraintsGraph_{\System}(V,E)$ is
robustly strongly connected.
\end{description}
\end{observation}
\Proof
Starting with part (a), in a square system $|\Supporters_i|=1$ and therefore
by definition, the two graphs coincide.
Next note that for a diagonal $\SupportersMatrix$
(as can be achieved by column reordering), $\ConstraintsGraph_{\System}(V,E)$
corresponds to $(\RepressorsMatrix)^{T}$ (by treating positive entries as $1$).
Since $\RepressorsMatrix$ is irreducible (and hence corresponds to a
strongly connected digraph), it follows that the matrix
$(\RepressorsMatrix)^{T}$ is irreducible, and hence
$\ConstraintsGraph_{\System}(V,E)$ is strongly connected.
To prove part (b), consider an arbitrary
$\FilterMatrix \in \FilterMatrixFamily$.
Since $\System(\FilterMatrix)$ is irreducible, it follows that
$\RepressorsMatrix \cdot \FilterMatrix$ is irreducible, and by
Obs. \ref{obs:reducible_graph_connected}(a),
$\ConstraintsGraph_{\System(\FilterMatrix)}(V,E)$ is strongly connected.
%This holds for every $\FilterMatrix \in \FilterMatrixFamily$.
The claim follows.
\QED
%%%%%%%%%%%%%%%%%%%%%%%%%%%%%%%%%%%%%%%%%%%%%
\paragraph{Partial selection for irreducible systems.}
%%%%%%%%%%%%%%%%%%%%%%%%%%%%%%%%%%%%%%%%%%%%
Let $\SelectionVec' \subseteq \Affectors$ be a subset of affectors in an irreducible system $\System$.
Then $\SelectionVec'$ is a \emph{partial selection}
if there exists a subset of entities $V' \subseteq \EntitySet$ such that (a) $|\SelectionVec'|=|V'|$, and (b) for every $\Entity_i \in V'$, $|\Supporters_i \cap \SelectionVec'|=1$.
\\
That is, every entity in $V'$ has a single representative supporter
in $\SelectionVec'$. We refer to $V'$ as the set of entities \emph{determined} by $\SelectionVec'$.
In the system $\System(\SelectionVec')$, the supporters $\Affectors_k$
of any $\Entity_i \in V'$ that were not selected by $\Entity_i$, i.e.,
$\Affectors_k \notin \SelectionVec' \cap \Supporters_i$, are discarded.
In other words, the system's affectors set consists of the selected supporters
$\SelectionVec'$, and the supporters of entities that have not made up
their selection in $\SelectionVec'$.
We now turn to describe $\System(\SelectionVec')$ formally. The set of affectors in $\System(\SelectionVec')$ is given by
$\Affectors(\System(\SelectionVec')) =
\SelectionVec' \cup
\bigcup_{\Supporters_i \cap \SelectionVec'=\emptyset}  \Supporters_i$.
The number of affectors in $\System(\SelectionVec')$ is denoted by
$m(\SelectionVec')=|\Affectors(\System(\SelectionVec'))|$.
Recall that the $j^{th}$ column of the
matrices $\SupportersMatrix ,\RepressorsMatrix$ corresponds to $\Affectors_j$.
Let $ind(\Affectors_j) =
j-|\{\Affectors_\ell \notin \Affectors(\System(\SelectionVec')), \ell\leq j-1\}|$
be the index of the affector $\Affectors_j$ in the new system,
$\System(\SelectionVec')$ (i.e, the $ind(\Affectors_j)^{th}$ column in the contracted matrices $\SupportersMatrix(\SelectionVec'),\RepressorsMatrix(\SelectionVec')$  corresponds to $\Affectors_j$).
Define the partial selection matrix $\FilterMatrix(\SelectionVec') \in \{0,1\}^{m \times m(\SelectionVec')}$
such that $\FilterMatrix(\SelectionVec')_{i,ind(\Affectors_j)}=1$ for every
$\Affectors_j \in \Affectors(\System(\SelectionVec'))$, and
$\FilterMatrix(\SelectionVec')_{i,j}=0$ otherwise.
Finally, let
$\System(\SelectionVec') = \displaystyle \langle
\SupportersMatrix(\SelectionVec'),\RepressorsMatrix(\SelectionVec')\rangle,$
where $\SupportersMatrix(\SelectionVec') =
\displaystyle \SupportersMatrix \cdot \FilterMatrix(\SelectionVec')
\mbox{~~and~~} \RepressorsMatrix(\SelectionVec') =
\RepressorsMatrix \cdot \FilterMatrix(\SelectionVec')$.
Note that $\SupportersMatrix(\SelectionVec'), \RepressorsMatrix(\SelectionVec')
 \in \R^{n \times m(\SelectionVec')}$.
Observe that if the selection $\SelectionVec'$ is a complete legal selection,
then $|\SelectionVec'|=n$ and the system $\System(\SelectionVec')$ is a square
system. In summary, we have two equivalent representations for square systems
in the nonsquare system $\System$:
\\
(a) by specifying a complete selection $\SelectionVec$, $|\SelectionVec|=n$,
and
\\
(b) by specifying the selection matrix,
$\FilterMatrix \in \FilterMatrixFamily$.
\\
Representations (a) and (b) are equivalent in the sense that the two square systems
$\System(\FilterMatrix(\SelectionVec))$ and $\System(\SelectionVec)$
are the same.
We now show that if the system $\System$ is irreducible, then so must be
any $\System(\SelectionVec')$, for any partial selection $\SelectionVec'$.

\begin{observation}
\label{obs:irreducible_selection}
Let $\System$ be an irreducible system. Then  $\System(\SelectionVec')$ is also
irreducible, for every partial selection $\SelectionVec'$.
\end{observation}
\Proof
Recall that a system is irreducible iff every hidden square system is
irreducible. I.e., the square system $\System(\FilterMatrix)$ is irreducible
for every $\FilterMatrix \in \FilterMatrixFamily(\System)$.
We now show that if
$\FilterMatrix \in \FilterMatrixFamily(\System(\SelectionVec'))$,
then $\FilterMatrix \in \FilterMatrixFamily(\System)$.
This follows immediately by  Eq. (\ref{eq:FilterMatrixFamily}) and the fact that
$\Supporters_i(\System(\SelectionVec')) \subseteq \Supporters_i(\System)$.
\QED

\par\noindent{\bf Agreement of partial selections.}
Let $\SelectionVec_1, \SelectionVec_2 \subseteq \Affectors$ be partial
selections for $V_1, V_2 \subseteq \EntitySet$ respectively.
Then we denote by $\SelectionVec_1 \sim \SelectionVec_2$ the property that the partial selections \emph{agree}, namely,
$\SelectionVec_1 \cap \Supporters_j=\SelectionVec_2 \cap \Supporters_j$
for every $\Entity_j \in V_1 \cap V_2$.
\begin{observation}
\label{obs:chain_sym}
Consider $V_1,V_2, V_3 \subseteq  \EntitySet$ determined by the partial selections
$\SelectionVec_1,\SelectionVec_2,\SelectionVec_3$ respectively, such that
$V_1 \subset V_2$, $\SelectionVec_1 \sim \SelectionVec_2$ and
$\SelectionVec_2 \sim \SelectionVec_3$.
Then also $\SelectionVec_3 \sim \SelectionVec_1$.
\end{observation}
\Proof
$\SelectionVec_{2}$ is more restrictive than $\SelectionVec_1$ since it defines
a selection for a strictly larger set of entities. Therefore every partial
selection  $\SelectionVec_3$ that agrees with $\SelectionVec_2$ agrees also
with $\SelectionVec_1$.
\QED

%%%%%%%%%%%%%%%%%%%%%%%%%%%%%%%%%%%%%%%%%%%%
\paragraph{Generalized \PFT~for nonnegative irreducible systems.}
%%%%%%%%%%%%%%%%%%%%%%%%%%%%%%%%%%%%%%%%%%%%
Recall that the root of a square system  $\System \in \SquareSystemFamily$ is
$\PFEigenValue(\System)=\max \left \{\EigenValue(Z(\System)) \right\}.$
$\PFEigenVector(\System)$ is the eigenvector of $Z(\System)$
corresponding to $\PFEigenValue(\System)$.
We now turn to define the \emph{generalized Perron--Frobenius (PF) root} of
a nonsquare system $\System \notin \SquareSystemFamily$, which is given by
\begin{equation}
\label{eq:general_pf_root}
\PFEigenValue(\System) ~=~ \min_{\FilterMatrix \in \FilterMatrixFamily}
\left \{\PFEigenValue(\System(\FilterMatrix)) \right\}.
\end{equation}
Let $\FilterMatrix^*$ be the selection matrix that achieves the minimum in
Eq. (\ref{eq:general_pf_root}). We now describe the corresponding eigenvector
$\PFEigenVector(\System)$. Note that $\PFEigenVector(\System) \in \R^{m}$,
whereas $\PFEigenVector(\System(\FilterMatrix^*)) \in \R^{n}$.

Consider $\overline{X}'=\PFEigenVector(\System(\FilterMatrix^*))$ and let
$\PFEigenVector(\System)=\overline{X}$, where
\begin{equation}
\label{eq:general_pf_vector}
X(\Affectors_j) ~=~
\begin{cases}
X'(\Affectors_j), & \text{if $\sum_{i=1}^{n}\FilterMatrix^{*}(j,i)>0$;}
\\
0, & \text{otherwise.}
\end{cases}
\end{equation}
We next state our main result, which is a generalized variant of the \PFT\ for every nonnegative nonsquare irreducible system.
\begin{theorem}
\label{thm:pf_ext}
Let $\System$ be an  irreducible and nonnegative nonsquare system. Then
\begin{description}
\item{(Q1)}
$\PFEigenValue(\System)>0$,
\item{(Q2)}
$\PFEigenVector(\System) \geq 0$,
\item{(Q3)}
$|NZ(\PFEigenVector(\System))|=n$,
\item{(Q4)}
$\PFEigenVector(\System)$ is not unique.
\item{(Q5)}
The generalized Perron root of $\System$ satisfies
$\displaystyle \PFEigenValue = \min\limits_{\overline{X} \in \mathcal{N}}
\left\{ \mathfrak{f}(\overline{X}) \right\}$, where
$$\mathfrak{f}(\overline{X}) ~=~ \max\limits_{1 \leq i \leq n, \left(\SupportersMatrix \cdot \overline{X} \right)_{i}\neq 0}
\{ \frac{\left(\RepressorsMatrix \cdot \overline{X} \right)_{i}}
{ \left(\SupportersMatrix \cdot \overline{X} \right)_{i}} \}$$
and $\mathcal{N}=\{\overline{X} \geq 0,||\overline{X}||_{1}=1,
\SupportersMatrix \cdot \overline{X}\neq 0\}.$
I.e., the Perron-Frobenius (PF) eigenvalue is $1/\beta^{*}$ where $\beta^{*}$
is the optimal value of Program (\ref{LP:Ext_Perron}),
and the \PFE~ is the corresponding optimal point. Hence for $\beta^{*}$, the $n$ constraints of
Eq. (\ref{eq:SR}) hold with equality.
\end{description}
\end{theorem}
%\begin{corollary}
%\label{cor:pf_ext}
%At the optimum value $\beta^{*}$, the set of $n$ constraints of
%Eq. (\ref{eq:SR}) hold with equality.
%\end{corollary}
%To ease of analysis we restrict attention from now on to a strong irreducible systems. The treatment in (not necessarily strong) irreducible system is defer to Section \ref{subsec:irreducible}.

%%%%%%%%%%%%%%%%%%%%%%%%%%%%%%%%%%%
\paragraph{The difficulty: Lack of log-convexity.}
%%%%%%%%%%%%%%%%%%%%%%%%%%%%%%%%%%%
Before plunging into a description of our proof, we first discuss
a natural approach one may consider for proving Thm. \ref{thm:pf_ext}
in general and solving Program (\ref{LP:Ext_Perron}) in particular,
and explain why this approach fails in this case.

A non-convex program can often be turned into an equivalent convex one
by performing a standard variable exchange.
This allows the program to be solved by convex optimization techniques (see \cite{TanFL11} for more information). An example for a program that's amenable to this technique is
Program (\ref{LP:Stand_Perron}), which is \emph{log-convex}
(see Claim \ref{cl:non_convex}(a)), namely, it becomes convex
after certain term replacements.
Unfortunately, in contrast with Program (\ref{LP:Stand_Perron}), the generalized
Program (\ref{LP:Ext_Perron}) is not log-convex
(see Claim \ref{cl:non_convex}(b)),
and hence cannot be handled in this manner.

More formally, for vector $\overline{X}=(X(1), \ldots, X(m))$ and $\alpha\in \R$,  denote the component-wise $\alpha$-power of $\overline{X}$ by $\overline{X}^{\alpha}=(X(1)^{\alpha}, \ldots, X(m)^{\alpha})$. An optimization program $\Pi$ is \emph{log-convex} if given two feasible solutions $\overline{X}_1, \overline{X}_2$ for $\Pi$, their log-convex combination
$\overline{X}_{\delta}= \overline{X}_1^{\delta} \cdot \overline{X}_2^{(1-\delta)}$ (where  ``$\cdot$" represents component-wise multiplication)
is also a solution for $\Pi$, for every $\delta \in [0,1]$.
In the following we ignore the constraint $||\overline{X}||_{1}=1$, since we only
validate the feasibility of nonzero nonnegative vectors; this constraint
can be established afterwards by normalization.
\begin{claim}
\label{cl:non_convex}
(a) Program (\ref{LP:Stand_Perron}) is log-convex (without the $||\overline{X}||_{1}=1$ constraint). \\
(b) Program (\ref{LP:Ext_Perron}) is not log-convex (even without the $||\overline{X}||_{1}=1$ constraint).
\end{claim}
\Proof
We start with (a). In \cite{LogConvex} it is shown that the power-control problem is log-convex.
The log-convexity of Perron-Frobenius eigenvalue is also discussed in \cite{Boyd-Conv-Opt-Book}, for completeness we prove it here.
We use the same technique of \cite{LogConvex} and show it directly for
Program (\ref{LP:Stand_Perron}). Let $A$ be a non-negative irreducible matrix and let $\overline{X}_1,\overline{X}_2$ be two feasible solutions for Program (\ref{LP:Stand_Perron}) with $\beta_1$, resp. $\beta_2$. We now show that $\overline{X}_3=\overline{X}_1^{\alpha} \cdot \overline{X}_2^{(1-\alpha)}$ (where  ``$\cdot$" represents entry-wise multiplication). is a feasible solution for $\beta_3=\beta_1^{\alpha} \cdot \beta_2^{1-\alpha}$, for any $\alpha \in [0,1]$. I.e., we show that $A \cdot \overline{X}_3 \leq 1/\beta_3 \cdot \overline{X}_3$. Let $\eta_i=X_1(i)/(A \cdot \overline{X}_1)_{i}$, $\gamma_i=X_2(i)/(A \cdot \overline{X}_2)_{i}$, $\delta_i=X_3(i)/(A \cdot \overline{X}_3)_{i}$. By the feasibility of $X_1$ (resp., $X_2$) it follows that $\eta_i \geq \beta_1$ (resp., $\gamma_i \geq \beta_2$) for every $i \in \{1, \ldots, n\}$.
It then follows that
\begin{equation}
\label{eq:log_con}
\frac{\delta_i}{\eta_i^{\alpha} \cdot \gamma_i^{1-\alpha}}= \frac{\left(\sum_{j}A(i,j) \cdot X_1(j)\right)^{\alpha} \cdot \left(\sum_{j}A(i,j) \cdot X_2(j)\right)^{1-\alpha} }{\sum_{j}A(i,j) \cdot X_1(j)^{\alpha} \cdot X_2(j)^{1-\alpha}}~.
\end{equation}
Let $p_j=\left(A(i,j)X_1(j)\right)^{\alpha}$ and $q_j=\left(A(i,j)X_2(j)\right)^{1-\alpha}$. Then Eq. (\ref{eq:log_con}) becomes
\begin{eqnarray*}
\label{eq:log_con2}
\frac{\delta_i}{\eta_i^{\alpha} \cdot \gamma_i^{1-\alpha}}&=& \frac{\left(\sum_{j} p_{j}^{1/\alpha}\right)^{\alpha} \cdot \left( \sum_{j} q_{j}^{1/(1-\alpha)}\right)^{1-\alpha} }{\sum_{j}p_j \cdot q_j}\geq 1
\end{eqnarray*}
where the last inequality follows by Holder Inequality which can be safely applied since $p_j,q_j \geq 0$ for every $j \in \{1, \ldots,n\}$. We therefore get that
for every $i$, $\delta_i \geq \eta_i^{\alpha} \cdot \gamma_i^{1-\alpha} \geq \beta_3$,
concluding that $X_3(i)/(A \cdot \overline{X}_3)_{i}\geq \beta_3$ and $A \cdot X_3 \leq 1/\beta_3 \cdot X_3$ as required. Part (a) is established. We now consider (b).
For vector $\overline{Y} \in \R^{m}$, $m \geq i$, recall that $\overline{Y}_{i}=(Y(1), \ldots, Y(i))$, the $i$ first coordinates of $\overline{Y}$.
For given repressor and supporter matrices
$\RepressorsMatrix,\SupportersMatrix\in\mathbb{R}^{n\times m}$,
define the following program. For $\overline{Y}\in \mathbb{R}^{m+1}$:
%\begin{align}
%\label{LP:not convex}
%  \max ~& Y({m+1})
%  \\
%  \mathrm{s.t.}~
%  &
%\displaystyle Y(m+1)\cdot \RepressorsMatrix\cdot(\overline{Y}_{m})^T
%     \leq \SupportersMatrix \cdot(\overline{Y}_m)^T & \nonumber
%\\
%  & \displaystyle \overline{Y} \geq \overline{0} & \nonumber
%\\
%  & \displaystyle \overline{Y}_m \neq \overline{0} & \nonumber
%\end{align}
%
\begin{eqnarray}
\label{LP:not convex}
  &&\max ~ Y({m+1}) ~
%  \\
  \mathrm{s.t.}~\\
  &&
\displaystyle Y(m+1)\cdot \RepressorsMatrix\cdot(\overline{Y}_{m})^T
     \leq \SupportersMatrix \cdot(\overline{Y}_m)^T  \nonumber
\\
  && \displaystyle \overline{Y} \geq \overline{0}  \nonumber
\\
  && \displaystyle \overline{Y}_m \neq \overline{0}  \nonumber
\end{eqnarray}
\noindent
This program is equivalent to Program (\ref{LP:Ext_Perron}).
An optimal solution $\overline{Y}$ for Program (\ref{LP:not convex})
``includes" an optimal solution for Program (\ref{LP:Ext_Perron}),
where $\beta=Y(m+1)$ and $\overline{X}=\overline{Y}_m$.
%
%We say that a program is convex, if for any given $0\leq\lambda\leq 1$
%and every two feasible solutions $\overline{Y},\overline{Z}$, it follows that
%$\lambda\cdot\overline{Y}+(1-\lambda)\overline{Z}$ is also feasible.
%
We prove that Program (\ref{LP:not convex}) is not log-convex by showing
the following example. Consider the repressor and supporters matrices
\begin{eqnarray*}
\RepressorsMatrix=
\begin{pmatrix}
0 & 2 & 1\\
1 & 0 & 0
\end{pmatrix}
\mbox{ ~~~and~~~ }
\SupportersMatrix=
\begin{pmatrix}
1/2 & 0 & 0 \\
0  & 4 &  4
\end{pmatrix}.
\end{eqnarray*}
It can be verified that $Y_1=(2, 1/2, 0, 1)$ and $Y_2= (4, 0, \sqrt{2}, \sqrt{2})$ are feasible. However, their log-convex combination $Y=Y_1^{1/2} \cdot Y_2^{1/2}$ is not a feasible solution for this system. Lemma follows.
\QED

\subsubsection{Algorithm for testing irreducibility}
%%%%%%%%%%%%%%%%%%%%%%%%%%%%%%%%%%%%%%%%%%%%

In this subsection, we provide a polynomial-time algorithm for testing the irreducibility
of a given nonnegative system $\System$. Note that if $\System$ is a square
system,  then irreducibility can be tested in a straightforward manner
by checking that $\RepressorsMatrix$ is irreducible and that $\SupportersMatrix$ is nonsingular.

However, recall that a nonsquare system $\System$ is irreducible iff every hidden square system $\System(\FilterMatrix)$,
$\FilterMatrix \in \FilterMatrixFamily$, is irreducible.
Since $\FilterMatrixFamily$ might be exponentially large, a brute-force testing of $\System(\FilterMatrix)$ for every $\FilterMatrix$ is too costly, hence another approach is needed.
Before presenting the algorithm, we provide some notation.
\par Consider a directed graph $G=(V,E)$.
Denote the set of incoming neighbors of a node $\Entity_k$ by $\Gamma^{in}(\Entity_k,D)=\{ \Entity_j \mid e_{j,i} \in E(D)\}$. The incoming neighbors of a set of nodes $V' \in \EntitySet$ is denoted $\Gamma^{in}(V',D)=\bigcup_{\Entity_k \in V'}\Gamma^{in}(\Entity_k,D)$.

\paragraph{Algorithm Description.}
To test irreducibility, Algorithm ~\TestIrred~ (see Fig. \ref{figure:irreducibility_tester}) must verify that the constraint graph $\ConstraintsGraph_{\System(\FilterMatrix)}$ of every $\FilterMatrix \in \FilterMatrixFamily$ is strongly connected.
The algorithm consists of at most $n-1$ rounds.
In round $t$, it is given as input a partition $\mathcal{C}^{t}=\{C^{t}_{1}, \ldots, C^{t}_{k_t}\}$ of $\EntitySet$ into $k_t$ disjoint clusters such that $\bigcup_{i} C^{t}_{i}=\EntitySet$. For round $t=0$, the input is a partition $\mathcal{C}^{0}=\{C^{0}_{1}, \ldots, C^{0}_{n}\}$ of the entity set $\EntitySet$ into $n$ singleton clusters $C^{0}_{i}=\{\Entity_i\}$.
The output at round $t$ is a coarser partition $\mathcal{C}^{t+1}$, in which at least two clusters of $\mathcal{C}^{t}$ were merged into a single cluster in $\mathcal{C}^{t+1}$. The partition $\mathcal{C}^{t+1}$ is formed as follows.
The algorithm first forms a graph $D_t=(\mathcal{C}^{t}, E_t)$ on the clusters of the input partition $\mathcal{C}^{t}$, treating each cluster $C^{t}_i \in \mathcal{C}^{t}$ as a node, and including in $E_t$ a directed edge $(i,j)$ from $C^{t}_i$ to $C^{t}_j$ if and only if there exists an entity node $\Entity_k \in C^{t}_{i}$ such that \emph{each} of its supporters $\Affectors_i \in \Supporters_k$ is a repressor of \emph{some} entity $\Entity_{k'} \in C^{t}_{j}$, i.e., $\Supporters_k \subseteq \bigcup_{\Entity_{k'} \in C^{t}_{j}} \Repressors_{k'}$.
\par  The partition $\mathcal{C}^{t+1}$ is now formed by merging clusters $C^{t}_j$ that belong to the same \SCC~ in $D_{t}$ into a single cluster $C^{t+1}_{k'}$ in $\mathcal{C}^{t+1}$.
Each cluster of $\mathcal{C}^{t+1}$ corresponds to a unique \SCC~ in $D_t$. If $D_t$ contains no \SCC~ except for singletons, which implies that no two cluster nodes of $D_t$ can be merged, then the algorithm declares the system $\System$ as reducible and halts. Otherwise, it proceeds with the new partition $\mathcal{C}^{t+1}$. Importantly, in $\mathcal{C}^{t+1}$ there are at least two entity subsets that belong to distinct clusters in $\mathcal{C}^{t}$ but to the same cluster node in $\mathcal{C}^{t+1}$. If none of the rounds ends with the algorithm declaring the system reducible (due to clusters ``merging" failure), then the procedure proceeds with the cluster merging until at some round $t^* \leq n-1$ the remaining partition  $\mathcal{C}^{t^*}=\{\{\EntitySet\}\}$ consists of a single cluster node that encompasses the entire entity set.

%%%%%%%%%%%%%%%%%%%%%%%%
%\begin{figure}[htb]
\begin{figure*}[h!]
\begin{center}
\framebox{\parbox{6in}{
\noindent{\bf Algorithm ~\TestIrred($\System$)~}
\begin{enumerate}
\dnsitem
$t \gets 0$;
\dnsitem
$k_t \gets n$;
\dnsitem
$C^{0}_i  \gets\{\Entity_i\}$ for every $i \in [1, k_{t}]$;
\dnsitem
$\mathcal{C}^{0} \gets \{C^{0}_1, \ldots, C^{0}_{k_t}\}$;
\dnsitem
While $|\mathcal{C}^{t}|>1$ do:
\begin{enumerate}
\ddnsitem
$\Repressors(C^{t}_i) \gets\bigcup_{\Entity_k \in C^{t}_i}  \Repressors_{k}$,
for every $i \in [1, k_{t}]$;
\ddnsitem
$E_t \gets \{e(i,j) \mid \exists \Entity_k \in C_i^t, \text{~such that~}
\Supporters_{k} \subseteq \Repressors(C^{t}_j)\}$.
\ddnsitem
Let $D_{t}=(\mathcal{C}^{t},E_t)$;
\ddnsitem
$k_{t+1} \gets$ number of \SCC s in $D_t$;
\ddnsitem
If $k_{t+1}=k_{t}$ and $|\mathcal{C}^{t}|\geq 2$, then return ``no".
\ddnsitem
%Else,
Decompose $D_{t}(\mathcal{C}^{t},E_t)$ into
%$k_{t+1}$
\SCC s $\widehat{C}^{1}, \ldots, \widehat{C}^{k_{t+1}}$.
\ddnsitem
$C^{t+1}_{i} \gets \bigcup_{C_{j} \in \widehat{C}^{i}} C_{j}$ for every $i \in [1,k_{t+1}]$.
\ddnsitem
$\mathcal{C}^{t+1} \gets \{C^{t+1}_1, \ldots, C^{t+1}_{k_{t+1}}\}$;
\ddnsitem
$t \gets t+1$;
\end{enumerate}
\dnsitem
Return ``yes";
\end{enumerate}
}}
\end{center}
\caption{\label{figure:irreducibility_tester}
The pseudocode of Algorithm ~\TestIrred.}
\end{figure*}
\paragraph{Analysis.}
We first provide some high level intuition for the correctness of the algorithm.
Recall, that the goal of the algorithm is to test whether the entire entity set $\EntitySet$ resides in a single \SCC~ in the constraint graph $\ConstraintsGraph_{\System(\FilterMatrix)}$ for every selection matrix $\FilterMatrix \in \FilterMatrixFamily$. This test is performed by the algorithm in a gradual manner by monotonically increasing the subsets of nodes that belong to the same \SCC~ in every $\ConstraintsGraph_{\System(\FilterMatrix)}$. In the beginning of the execution, the most one can claim is that every entity $\Entity_k$ is in its own \SCC. Over time, clusters are merged while maintaining the invariant that all entities of the same cluster belong to the same \SCC~ in every $\ConstraintsGraph_{\System(\FilterMatrix)}$.
More formally, the following invariant is maintained in every round $t$: the entities of each cluster $C^t_i \subseteq \EntitySet$ of the graph $D_t$ are guaranteed to be in the same \SCC~ in the constraint graph $\ConstraintsGraph_{\System(\FilterMatrix)}$ for every selection matrix
$\FilterMatrix \in \FilterMatrixFamily$.
We later show that if the system $\System$ is irreducible, then the merging process never fails and therefore the last partition $\mathcal{C}^{t^*}=\{\{\EntitySet\}\}$ consists of a single cluster node that contains all entities, and by the invariant, all entities are guaranteed to be in the same \SCC~ in the constraint graph of any hidden square subsystem.
\par We now provide some high level explanation for the validity of this invariant. Starting with round $t=0$, each cluster node $C^0_i=\{\Entity_i\}$ is a singleton and every singleton entity is trivially in its own \SCC~ in any constraint graph $\ConstraintsGraph_{\System(\FilterMatrix)}$. Assume the invariant holds up to round $t$, and consider round $t+1$.
The key observation in this context is that the new partition $\mathcal{C}^{t+1}$ is defined based on the graph $D_t=(\mathcal{C}^{t}, E_{t})$, whose edges are independent of the specific supporter selection that is made by the entities (and that determines the resulting hidden square subsystem). This holds due to the fact that a directed edge $(i,j) \in E_{t}$ between the clusters $C^{t}_{i}, C^{t}_{j} \in \mathcal{C}^{t}$ exists if and only if there exists an entity node $\Entity_k \in C^{t}_{i}$ such that \emph{each} of its supporter $\Affectors_i \in \Supporters_k$ is a repressor of \emph{some} entity $\Entity_{k'} \in C^{t}_{j}$. Therefore, if the edge $(i,j)$ exists in the $D_{t}$, then it exists also in the cluster graph corresponding to the constraint graph $\ConstraintsGraph_{\System(\FilterMatrix)}$ (i.e., the graph formed by representing every \SCC~ of $\ConstraintsGraph_{\System(\FilterMatrix)}$ by a single node) for \emph{every} hidden square subsystem $\System(\FilterMatrix)$, no matter which supporter $\Affectors_i \in \Supporters_k$ was selected by $\FilterMatrix$ for $\Entity_k$. Hence, under the assumption that the invariant holds for $\mathcal{C}^{t}$, the coarse-grained representation of the clusters of $\mathcal{C}^t$ in $\mathcal{C}^{t+1}$ is based on their membership in the same \SCC~  in the ``selection invariant" graph $D_{t}$, thus the invariant holds also for $t+1$.

We next formalize this argumentation. We say that round $t$ is \emph{successful} if $D_t$ contains a \SCC~ of size greater than 1.  We begin by proving the following.
\begin{claim}
\label{cl:partition_induc}
For every successful round $t$, the partition
$\mathcal{C}^{t+1}$ satisfies the following properties.
\begin{description}
\item{(A1)}
$\mathcal{C}^{t+1}$ is a partition of $\EntitySet$, i.e.,
$C^{t+1}_i \subseteq \EntitySet$, $C^{t+1}_{j} \cap C^{t+1}_{i}=\emptyset$
for every $i,j \in [1,k_{t+1}]$, and $\bigcup_{j\leq k_{t+1}} C^{t+1}_{j}=\EntitySet$.
\item{(A2)}
Every $C^{t+1}_{j} \in \mathcal{C}^{t+1}$ is a \SCC~ in the constraint graph
$\ConstraintsGraph_{\System(\FilterMatrix)}$ for every selection matrix
$\FilterMatrix \in \FilterMatrixFamily$.
\end{description}
\end{claim}
\Proof
By induction on $t$.
Clearly, since $C^{0}_{i}=\{\Entity_i\}$ for every $i$,
Properties (A1) and (A2) trivially hold for $\mathcal{C}^{0}$.
We now show that if round $t=0$ is successful, then (A1) and (A2) hold for $\mathcal{C}^{1}$. Since the edges of $D_0$ exist also in the corresponding cluster graph of $\ConstraintsGraph_{\System(\FilterMatrix)}$ under any selection $\FilterMatrix$ of the entities, the clusters of $\mathcal{C}^{0}$ that are merged into a single \SCC~ in $\mathcal{C}^{1}$, belong also to the same \SCC~ in the constraint graph $\ConstraintsGraph_{\System(\FilterMatrix)}$ of every $\FilterMatrix \in \FilterMatrixFamily$.
Next, assume these properties to hold for every round up to $t-1$ and consider round $t$.
Since round $t$ is successful, any prior round $t' <t$ was successful as well, and thus the induction assumption can be applied on round $t-1$. In particular, since $\mathcal{C}^{t+1}$ corresponds to \SCC s of $D_t$, it represents a partition of the clusters of
$\mathcal{C}^{t}$. By the induction assumption for round $t-1$, Property (A1) holds for $\mathcal{C}^{t}$ and therefore $\mathcal{C}^{t}$ is a partition of the entity set $\EntitySet$. Since $\mathcal{C}^{t+1}$ corresponds to a partition of $\mathcal{C}^t$, it is a partition of $\EntitySet$ as well so (A1) is established. Property (A2) holds for $\mathcal{C}^{t+1}$ by the same argument provided for the induction base.
%Finally, as round $t$ is successful,
%construction as well (recalling that $E_{t+1}$ exists for every selection) and by (A2) for round $t$. It remains to establish (A3). This trivially holds as $\mathcal{C}^{t+1}$ is formed by merging least two components of $\mathcal{C}^{t}$  into a single cluster in $\mathcal{C}^{t+1}$. If no two components of $\mathcal{C}^{t-1}$ can be merged then the algorithm terminates and do not construct $\mathcal{C}^{t}$.
The claim follows.
\QED
We next show that the algorithm return ``yes" for every irreducible system.
Specifically, we show that for an irreducible system, if $|\mathcal{C}^{t}|>1$ then round $t$ is \emph{successful}, i.e., the merging operation of the cluster graph $D_t$ succeeds. Once $\mathcal{C}^t$ contains a single cluster (containing all entities), the algorithm terminates and returns ``yes".
We first provide an auxiliary claim.
\begin{claim}
\label{cl:aux}
If $\System$ is irreducible and $|\mathcal{C}^{t}|>1$, then
$|\Gamma^{in}(C^{t}_{j},D_{t})| \geq 1$ for every $C^{t}_{j} \in \mathcal{C}^{t}$.
\end{claim}
\Proof
First note that if $\mathcal{C}^{t}$ is defined, then round $t-1$ was successful. Therefore, by Property (A1) of Cl. \ref{cl:partition_induc}, $\mathcal{C}^{t}$
is a partition of the entity set $\EntitySet$.
Assume, towards contradiction that the claim does not hold, and let $C^{t}_{j}\in \mathcal{C}^{t}$ be such that
$\Gamma^{in}(C^{t}_{j},D_{t})=\emptyset$. Denote the set of incoming neighbors
of component $C^{t}_{j}$ in the constraint graph $\ConstraintsGraph_{\System}$ by $W=\Gamma^{in}(C^{t}_{j},\ConstraintsGraph_{\System}) \setminus C^{t}_{j}$.
Since $\ConstraintsGraph_{\System}$ is irreducible, the vertices of $C^{t}_{j}$ are reachable from the outside, so $W \neq \emptyset$.
Let the repressors set of $C^{t}_{j}$ be
$\Repressors(C^{t}_{j})=\bigcup_{\Entity_k \in C^{t}_{j}} \Repressors_{k}$.
We now construct a square hidden system $\System(\FilterMatrix^*)$ which is reducible, in contradiction to the irreducibility of $\System$. Specifically, we look for a selection matrix $\FilterMatrix^*$ satisfying that for every entity $\Entity_k \in W$, its selected supporter $\Affectors_k$ in  $\System(\FilterMatrix^*)$ (i.e., the one for which $\FilterMatrix^*( \Affectors_k, \Entity_k)=1$) is not a repressor of any of the entities in $C^{t}_{j}$, i.e., $\Affectors_k
\in \Supporters_k \setminus \Repressors(C^{t}_{j})$.
Recall, that since $\System$ is irreducible, the supporter sets $\Supporters_i, \Supporters_j$ are pairwise disjoint (see Claim \ref{cor:distinct}).
Note that since $\Gamma^{in}(C^{t}_{j},D_{t})=\emptyset$,
such a selection matrix $\FilterMatrix^*$ exists.
To see this, assume, towards contradiction that $\FilterMatrix^*$ does not exist. This implies that there exists an entity $\Entity_{k} \in W$ such that $\Supporters_{k} \setminus \Repressors(C^{t}_{j})=\emptyset$ and therefore an affector in $\Supporters_k \setminus \Repressors(C^{t}_{j})$ could not be selected for $\FilterMatrix^*$. Hence,  $\Supporters_k \subseteq \Repressors(C^{t}_{j})$.
Let $C^{t}_{i} \in \mathcal{C}^{t}$ be the cluster such that $\Entity_k \in C^{t}_{i}$. Since $\mathcal{C}^{t}$ is a partition of the entity set $\EntitySet$, such $ C^{t}_{i}$ exists. Since $\Supporters_k \subseteq \Repressors(C^{t}_{j})$, it implies that the edge $e_{i,j} \in D_{t}$, in contradiction to the fact that $C^{t}_{j}$ has no incoming neighbors in $D_{t}$.
%there exists some $\Entity_k \in W$ in some cluster $C^{t}_{j'}$ (since $\mathcal{C}^{t}$ is a partition of the entity set $\EntitySet$ such $C^{t}_{j'}$ exists) such that all its supporters are in $\Repressors(C^{t}_{j})$, then by the definition of $D_{t}$, the edge $(j', j)$ exists, in contradiction to the fact that the cluster $C^{t}_{j}$ has no incoming edges in $D_{t}$.
We therefore conclude that $\FilterMatrix^*$ exists.
\par We now show that $\System(\FilterMatrix^*)$ is reducible. In particular, we show that
the incoming degree of the component $C^{t}_{j}$
(from entities in other components) in the constraint graph $\System(\FilterMatrix^*)$ of the square system $\System(\FilterMatrix^*)$, is zero, i.e.,
$\Gamma^{in}(C^{t}_{j},\ConstraintsGraph_{\System(\FilterMatrix^*)})=\emptyset$.
Assume, towards contradiction, that there exists a directed edge $e_{x,y}$ from entity $\Entity_x \in \EntitySet \setminus C^{t}_{j}$ to some $\Entity_y \in C^{t}_{j}$ in $\ConstraintsGraph_{\System(\FilterMatrix^*)}$. This implies that $e_{x,y} \in \ConstraintsGraph_{\System}$ exists in the constraint graph of the original (nonsquare) system $\System$ and thus $\Entity_x$ is in $W$. Let $\Affectors_{x'} \in \Supporters_x$ be the selected supporter
of $\Entity_x$ in $\FilterMatrix^*$. By construction of $\FilterMatrix^*$,
$\Affectors_{x'}  \notin \Repressors(C^{t}_{j})$, in contradiction to the fact that the edge $e_{x,y} \in \ConstraintsGraph_{\System(\FilterMatrix^*)}$ exists.

Since there exists a node in $\ConstraintsGraph_{\System(\FilterMatrix^*)}$ with no incoming neihbors, this graph is not strongly connected, implying that $\System(\FilterMatrix^*)$ is reducible.

Finally, as $\System$ is irreducible, it holds that every hidden square system is irreducible, in particular $\System(\FilterMatrix^*)$, hence, contradiction. The claim follows.
\QED
\begin{lemma}
\label{cl:pos}
If $\System$  is irreducible then
Algorithm ~\TestIrred($\System$) returns ``yes".
\end{lemma}
\Proof
By Cl. \ref{cl:aux}, we have that if $\System$ is irreducible and
$|\mathcal{C}^{t}|>1$, then every node in $D_t$ has an incoming edge, which necessitates that there exists a (directed) cycle
$C=(C_{i_1}, \ldots, C_{i_k})$, for $k \geq 2$ in $D^{t}$. Since the nodes in such cycle $C$ are strongly connected, they can be merged in $\mathcal{C}^{t+1}$, and therefore round $t$ is successful.
Moreover, since at least two clusters of $\mathcal{C}^t$ are merged into a single cluster in $\mathcal{C}^{t+1}$, we have that
$|\mathcal{C}^{t+1}|<|\mathcal{C}^{t}|$.
This means that the merging never fails as long as $|\mathcal{C}^{t}|>1$, so $k_{t}=|\mathcal{C}^{t}|$ is monotonically decreasing.
It follows that the algorithm terminates within at most $n-1$ rounds with a ``yes". The Lemma follows.
\QED
We now consider a reducible system $\System$ and show
that ~\TestIrred($\System$) returns ``no".
\begin{lemma}
\label{lem:neg}
If $\System$  is reducible, then
Algorithm ~\TestIrred($\System$) returns ``no".
\end{lemma}
\Proof
Towards contradiction, assume otherwise, i.e., suppose that the algorithm accepts $\System$.
This implies that every round $t \in [1, t^*]$ in which $|\mathcal{C}^{t}|>1$ is successful.
\par The reducibility of  $\System$ implies that there exists (at least one) hidden square system  $\System(\FilterMatrix)$ which is reducible, namely, its constraint graph $\widehat{D}=\ConstraintsGraph_{\System(\FilterMatrix)}$ is not strongly connected. Thus $\widehat{D}$ contains at least two nodes $\Entity_i$ and $\Entity_j$ that belong to distinct \SCC s in $\widehat{D}$.
Note that $\Entity_i$ and $\Entity_j$ are in distinct clusters in $\mathcal{C}^{0}$, but belong to the same cluster in the partition of the final $\mathcal{C}^{t^*}$. Therefore, there must exists a round $t' \in (0, t^*)$
in which the cluster $C^{t'}_{i'}$ that contains $\Entity_i$ and the cluster $C^{t'}_{j'}$ that contains $\Entity_j$ appeared in the same \SCC~ in $D_{t'}$ and were merged into a single \SCC~ in $\mathcal{C}^{t'+1}$.
(Note that since $t'-1$ is a successful round,
$\mathcal{C}^{t'}$ is a partition of the entity set (Prop. (A1) of Cl. \ref{cl:partition_induc}) and therefore $C^{t'}_{i'}$ and $C^{t'}_{j'}$ exist.)
Since round $t'$ is successful (otherwise the algorithm would terminates with ``no"), by to Property (A2) of Cl. \ref{cl:partition_induc}, it follows that the entity subset of the unified cluster
$\mathcal{C} \in \mathcal{C}^{t'+1}$ is in the same connected component in the constraint graph $\ConstraintsGraph_{\System(\FilterMatrix')}$ for every $\FilterMatrix' \in \FilterMatrixFamily$. Since $\FilterMatrix \in \FilterMatrixFamily$ as well it holds that $\Entity_i$ and $\Entity_j$ are in the same connected component in $\widehat{D}$. Hence, contradiction. The lemma follows.
\QED
By Lemmas \ref{cl:pos} and \ref{lem:neg} it follows that
Algorithm ~\TestIrred($\System$) returns ``yes" iff the system $\System$  is irreducible, which establish the correctness of the algorithm.

%\paragraph{Running time}
\begin{claim}
\label{cl:runtime}
Algorithm ~\TestIrred~ terminates in $O(m \cdot n^2)$ rounds.
\end{claim}
\Proof
The algorithm consists of at most $n-1$ rounds
%$t=1, \ldots, n-1$.
In each round $t$, it constructs the cluster graph
$D_t=(\mathcal{C}^{t-1}, E_t)$ in time $O(n \cdot m)$.
The decomposition into \SCC s can be done
in $O(|D_t|)=O(n^2)$. The claim follows.
\QED

%We have the following.

\begin{theorem}
\label{lem:alg_irred}
There exists a polynomial time algorithm %~\TestIrred($\System$)
for deciding irreducibility on nonnegative systems.
\end{theorem}
%%%%%%%%%%%%%%%%%%%%%%%%%%%%%%%%%%%%%%%%%%%%
\section{Proof of the generalized \PFT}
%%%%%%%%%%%%%%%%%%%%%%%%%%%%%%%%%%%%%%%%%%%%

%%%%%%%%%%%%%%%%%%%%%%%%%%%%%%%%%%%%%%%%%%%%
\subsection{Proof overview and roadmap}
%\paragraph{Proof overview.}
%%%%%%%%%%%%%%%%%%%%%%%%%%%%%%%%%%%%%%%%%%%%
Our main challenge is to show that the optimal value of Program
(\ref{LP:Ext_Perron}) is related to an \emph{eigenvalue} of some hidden
square system  $\System^{*}$ in $\System$ (where ``hidden" implies that
there is a selection on $\System$ that yields $\System^{*}$).
The flow of the analysis is as follows.
%In Subsection \ref{subsec:constraints_graph},
In Subsec. \ref{sec:geometry_n_1}, we consider a convex relaxation of Program (\ref{LP:Ext_Perron}) and show that the set of feasible solutions
of Program (\ref{LP:Ext_Perron}), for every $\beta \in (0, \beta^{*}]$, corresponds to a bounded polytope.
By dimension considerations, we then show that the vertices of such polytope correspond to feasible solutions with at most $n+1$ nonzero entries.
In Subsec. \ref{sec:weak}, we show that for irreducible systems, each vertex of such a polytope corresponds to a hidden \emph{weakly square} system
$\System^{*} \in \WeakSystemFamily$. That is, there exists a hidden weakly square system in $\System$ that achieves $\beta^{*}$. Note that a solution for such a hidden system can be extended to a solution for the original $\System$ (see Obs. \ref{obs:filter_to_square}).

Next, in Subsec. \ref{sec:zerostar},
we exploit the generalization of Cramer's rule for homogeneous linear systems (Cl. \ref{cl:cramer_non_square})
as well as a separation theorem for nonnegative matrices to show that there is
a hidden optimal \emph{square} system in $\System$ that achieves $\beta^{*}$,
which establishes the lion's share of the theorem.

Arguably, the most surprising conclusion of our generalized theorem is that
although the given system of matrices is not square, and eigenvalues cannot
be straightforwardly defined for it, the nonsquare system contains
a \emph{hidden optimal} square system, optimal in the sense that a solution $\overline{X}$ for this system can be translated into a solution $\overline{X}^m$ to the original system
(see Obs. \ref{obs:filter_to_square}) that satisfies
Program (\ref{LP:Ext_Perron}) with the optimal value $\beta^{*}$.
The power of a nonsquare system is thus not in the ability to create a solution
better than \emph{any} of its hidden square systems, but rather in
the \emph{option} to \emph{select} the best hidden square system
out of the possibly exponentially many ones.

%%%%%%%%%%%%%%%%%%%%%%%%%%%%%%%%%%%%%%%%%%%%
\subsection{Existence of a solution with $n+1$ affectors}
\label{sec:geometry_n_1}
%\paragraph{The Geometry of the \PFT~with selection.}
%\par\noindent{\bf The Geometry of the \PFT~with selection.}
%%%%%%%%%%%%%%%%%%%%%%%%%%%%%%%%%%%%%%%%%%%%

We now turn to characterize the feasible solutions of
Program (\ref{LP:Ext_Perron}).
The following is a convex variant of Program (\ref{LP:Ext_Perron}).
\begin{align}
\label{LP:Ext_Perron_convex}
  \mbox{maximize~~} ~& 1 \mbox{~~subject to:~~}
  \\
   & \displaystyle \RepressorsMatrix \cdot \overline{X} ~\leq~
       1/\beta \cdot \SupportersMatrix \cdot \overline{X} ~,&
  \label{eq:SR-convex}  \\
  & \displaystyle \overline{X} \geq \overline{0}~, &
  \label{eq:Ineq-convex}\\
  & \displaystyle ||\overline{X}||_{1}=1~. &
   \label{eq:eq-one-convex}
\end{align}

Note that Program (\ref{LP:Ext_Perron_convex}) has the same set of constraints
as those of Program (\ref{LP:Ext_Perron}). However, due to the fact that
$\beta$ is no longer a variable, we get the following.
\begin{claim}
\label{cl:convex}
Program (\ref{LP:Ext_Perron_convex}) is convex.
\end{claim}
%\textbf{MP:I add this here but maybe it should be placed elsewhere.}
To characterize the set of feasible solutions $(\overline{X}, \beta)$, $\beta>0$ of Program (\ref{LP:Ext_Perron}), we fix some $\beta>0$, and characterize the solution set of Program (\ref{LP:Ext_Perron_convex}) with this $\beta$.
It is worth noting at this point that using the above convex relaxation,
one may apply a binary search for finding a {\em near-optimal} solution
for Program (\ref{LP:Ext_Perron_convex}), up to any predefined accuracy.
%Yet, this results in an approximate solution and not in exact one.
In contrast, our approach, which is based on exploiting the special
geometric characteristics of the optimal solution,
enjoys the theoretically pleasing (and mathematically interesting) advantage
of leading to
%demonstrating the existence of
an efficient algorithm for computing the optimal solution precisely, and thus establishing the polynomiality of the problem.

Throughout, we restrict attention to values of $\beta \in (0, \beta^{*}]$.
Let $\Polytope(\beta)$ be the polyhedron corresponding to Program
(\ref{LP:Ext_Perron_convex}) and denote by $V(\Polytope(\beta))$
the set of vertices of $\Polytope(\beta)$.
\begin{claim}
\label{cl:n_zero_polytope}
(a) $\Polytope(\beta)$ is bounded (or a polytope).
(b) For every $\overline{X} \in V(\Polytope(\beta))$,
$|NZ(\overline{X})| \leq n+1$. This holds even for reducible systems.
\end{claim}
\Proof
Part (a) holds by the Equality constraint (\ref{eq:eq-one-convex}) which enforces
$||\overline{X}||_{1}~=~1$. We now prove Part (b).
Every vertex $\overline{X} \in \R^{m}$ is defined by a set of $m$ linearly
independent equalities. Recall that one equality is imposed by the constraint
$||\overline{X}||_{1}~=~1$ (Eq. (\ref{eq:eq-one-convex})).
Therefore it remains to assign $m-1$ linearly independent equalities out of
the $n+m$ (possibly dependent) inequalities of Program
(\ref{LP:Ext_Perron_convex}). Hence even if all
%To achieve a lower bound on the number of zeros in $\overline{v}$, we distribute the $m-1$ equalities as follows:
the (at most $n$) linearly independent SR constraints (\ref{eq:SR-convex})
become equalities, we are still left with at least $m-1-n$ unassigned
equalities, which must be taken from the remaining $m$ nonnegativity constraints (\ref{eq:Ineq-convex}). Hence, at most $n+1$ nonnegativity inequalities
were not fixed to zero, which establishes the proof.
\QED

%%%%%%%%%%%%%%%%%%%%%%%%%%%%%%%%%%%
\subsection{Existence of a weak $\ZeroStar$-solution}
\label{sec:weak}
%%%%%%%%%%%%%%%%%%%%%%%%%%%%%%%%%%%
We now consider the case where the system $\System$ is irreducible and
a more delicate characterization of $V(\Polytope(\beta))$ can be deduced.

We begin with some definitions.
A solution $\overline{X}$ is called a {\em $\Zero$ solution}
(for Program (\ref{LP:Ext_Perron}))
if it is a feasible solution
$\overline{X}_{\widetilde{\beta}}$, $\widetilde{\beta} \in (0,\beta^{*}]$,
in which for each $\Entity_i \in \EntitySet$ only one affector has a non-zero
assignment, i.e., $NZ(\overline{X}) \cap \Supporters_i=1$ for every $i$.
A solution $\overline{X}$ is called a {\em $\WeakZero$ solution},
or a {\em ``weak'' $\Zero$ solution}, if it is a feasible vector
$\overline{X}_{\widetilde{\beta}}$, $\widetilde{\beta} \in (0,\beta^{*}]$,
in which for each $\Entity_i$, {\em except at most one}, say
$\Entity_\ell \in \EntitySet$, $|NZ(\overline{X}) \cap \Supporters_i| = 1$,
$\Entity_i \in \EntitySet \setminus \{\Entity_\ell\}$ and
$|NZ(\overline{X}) \cap \Supporters_\ell| = 2$.
A solution $\overline{X}$ is called a {\em $\ZeroStar$ solution} if it is
an optimal $\Zero$ solution.
Let $\WeakZeroStar$ be an optimal $\WeakZero$ solution.

%\subsubsection{$\WeakZero$ solutions}
%\label{sec:weak}
%\paragraph{$\WeakZero$ solutions.}

For a feasible vector $\overline{X}$, we say that $\Affectors_k$ is
\emph{active} in $\overline{X}$ iff $X(\Affectors_k)>0$.
A subgraph $\GCGraph$ of a constraint graph $\ConstraintsGraph_{\System}$ is \emph{active} in $\overline{X}$ iff every edge
in $\GCGraph$ can be associated with (or ``explained by") an active affector, namely,
$$e(i, j) \in E(\GCGraph) ~~~\mbox{~iff~}~~~
\Supporters_{i} \cap \Repressors_{j} \cap NZ(\overline{X}) \neq \emptyset.$$
\par Towards the end of this section, we prove the following lemma which holds for every feasible solution of
Program (\ref{LP:Ext_Perron_convex}).
\begin{lemma}
\label{cl:entity_one_nonzero}
Let $\System$ be an irreducible system with a feasible solution $\overline{X}_{\beta}$ of Program (\ref{LP:Ext_Perron}). For every entity $\Entity_i$ there exists
an active affector $\Affectors_{\IndS(i)} \in \Supporters_i$, such that
$X_{\beta}(\Affectors_{\IndS(i)}) >0$, or in other words,
$\Supporters_{i} \cap NZ(\overline{X}_{\beta})\neq \emptyset$.
\end{lemma}
Let $\SelectionVec'$ be a partial selection determining $V' \subseteq \EntitySet$.
Define the collection of constraint graphs agreeing with $\SelectionVec'$ as
\begin{equation}
\label{eq:graph_family_selection}
\mathfrak{G}(\SelectionVec') ~=~ \{\ConstraintsGraph_{\System(\SelectionVec)} \mid
\text{a complete selection~} \SelectionVec \text{~satisfying~}\SelectionVec \sim \SelectionVec'\}.
\end{equation}
Note that by Obs. \ref{obs:reducible_graph_connected}(b), every
constraint graph $\GCGraph \in \mathfrak{G}(\SelectionVec')$ for every partial selection
$\SelectionVec'$ is strongly connected.
I.e., $\mathfrak{G}(\SelectionVec')$ contains the constraint graphs
for all square systems restricted to the partial selection dictated by
$\SelectionVec'$ for $V'$.
Note that when $|\SelectionVec'|=n$, $\SelectionVec'$ is a complete selection,
i.e., $\FilterMatrix(\SelectionVec') \in \FilterMatrixFamily$, and
$\mathfrak{G}(\SelectionVec')$ contains a single graph
$\ConstraintsGraph_{\System(\SelectionVec')}$ corresponding to the square system
$\System(\SelectionVec')$.

Given a feasible vector $\overline{X}$ and an irreducible system $\System$, the main challenge is to find an active (in $\overline{X}$) irreducible spanning
subgraph of $\ConstraintsGraph_{\System}$.  Finding such a subgraph is crucial for both
Lemma \ref{cl:entity_one_nonzero} and Lemma \ref{lem:strict_equality} later on.

We begin by showing that given just one active affector $\Affectors_{p_1}$ in $\overline{X}$, it is possible to ``bootstrap" it and construct an active irreducible
spanning subgraph of $\ConstraintsGraph_{\System}$ (in $\overline{X}$).

Let $\Entity_{i_1}$ be an entity satisfying that $\Affectors_{p_1} \in \Supporters_{i_1}$.
(Such entity $\Entity_{i_1}$  must exist, since there are no redundant affectors).
In what follows, we build an ``influence tree" starting at $\Entity_{i_1}$ and
spanning the entire set of entities $\EntitySet$.

For a directed graph $G$ and vertex $v \in G$ let
$\BFS(G, v)$ be the \emph{breadth-first search} tree of $G$ rooted at $v$, obtained by placing vertex $w$ at level $i$ of the tree if the shortest directed path from $v$ to $w$ is of length $i$.
Given a constraint graph $\GCGraph$, let $L_{i}(\GCGraph)$ be the $i^{th}$ level of $\BFS(\GCGraph,\Entity_{i_1})$.

We now describe an iterative process for constructing a complete selection
$\SelectionVec^{*}$ of $n$ supporters with positive entries in $\overline{X}_{\beta}$,
i.e., such that $\SelectionVec^{*} \subseteq NZ(\overline{X}_{\beta})$ and
$|\Supporters_i \cap \SelectionVec^{*}| =1$ for every $\Entity_i$.
At step $t$, we start from the partial selection $\SelectionVec_{t-1}$ constructed in the previous step, and extend it to $\SelectionVec_{t}$.
The partial selection $\SelectionVec_{t}$ should satisfy the following four properties.
%(A1) $\SelectionVec_{t}$ defines a partial selection. \\
%There exists a bijection function, $\SelectionFunction^{t}(\Entity_i)=\Affectors_j$ such that $\Affectors_j \in \SelectionVec_{t} \cap \Supporters_i$.\\
\begin{description}
\item{(A1)}
$\SelectionVec_{t} \subseteq NZ(\overline{X}_{\beta})$ (i.e., it consists of strictly positive supporters).
\end{description}
Consider the graph family $\mathfrak{G}(\SelectionVec_{t})$
defined in Eq. (\ref{eq:graph_family_selection}), consisting of all constraint
graphs for square systems induced by a selection that agrees with
$\SelectionVec_{t}$.
\begin{description}
\item{(A2)}
For every $i \in \{0, \ldots, t-1\}$ it holds that $L_{i}(\GCGraph_1) = L_{i}(\GCGraph_2), \text{~for every~} \GCGraph_1, \GCGraph_2 \in \bigcup_{j=i}^{t} \mathfrak{G}(\SelectionVec_{j})$, i.e., from step $i$ ahead, the $i$'th first levels coincide.
\item{(A3)}
$L_{t}(\GCGraph_1) = L_{t}(\GCGraph_2), \text{~for every~} \GCGraph_1, \GCGraph_2 \in
\mathfrak{G}(\SelectionVec_{t})$, (i.e., level $t$ coincides as well).
\end{description}
Denote $\CGLevel_{i}=L_{i}(\GCGraph)$, $\GCGraph \in \mathfrak{G}(\SelectionVec_{t})$,
for $i \in \{0, \ldots, t\}$ (by (A2) and (A3) this is well-defined). Let $Q_{-1}=\emptyset$, and $Q_{t}=\bigcup_{i=0}^{t} \CGLevel_{i}$ for $t \geq 0$, be set of entities in the first $t$ levels of $\mathfrak{G}(\SelectionVec_{t})$ graphs.
\begin{description}
\item{(A4)}
$\SelectionVec_{t}$ is a partial selection determining the entities in $Q_{t-1}$,
(i.e.,  $|\SelectionVec_{t}|=|Q_{t-1}|$  and $|\SelectionVec_t \cap \Supporters_i|=1$ for every $\Entity_i \in Q_{t-1}$).
\end{description}
Let us now describe the construction process of $\SelectionVec^{*}$ in more detail.
At step $t=0$, let $\SelectionVec_{0}=\emptyset$. Note that in this case
$$\mathfrak{G}(\SelectionVec_{0}) ~=~
\{\ConstraintsGraph_{\System(\FilterMatrix)} ~\mid~ \FilterMatrix \in
\FilterMatrixFamily\}.$$
It is easy to see that Properties (A1)-(A4) are satisfied.
For $t=1$, let $\SelectionVec_{1}=\{\Affectors_{p_1}\}$. As $L_0(\GCGraph)=\{\Entity_{i_1}\}$ and $L_1(\GCGraph)=\{\Entity_{i_2}~\mid~ \Affectors_{p_1} \in \Repressors_{i_2}\}$ for every $\GCGraph \in \mathfrak{G}(\SelectionVec_{1})$, Properties (A2) and (A3) holds. Property (A4) holds as well since $\SelectionVec_{1}$ determines $Q_0=\{\Entity_{i_1}\}$.

\par Now assume that Properties (A1)-(A4) hold after step $t$ (for $t \geq 1$), and
consider step $t+1$. We show how to construct $\SelectionVec_{t+1}$ given
$\SelectionVec_{t}$, and then show that it satisfies Properties (A1)-(A4).
Note that by definition
$\CGLevel_{t} \subseteq \EntitySet \setminus Q_{t-1}$.
Our goal is to find a partial selection $\Delta_t$ determining $\CGLevel_{t}$ such that $\Delta_t \subseteq NZ(\overline{X}_{\beta})$

Once finding such a set $\Delta_t$, the partial selection $\SelectionVec_{t+1}$ is taken to be
$\SelectionVec_{t+1}=\SelectionVec_{t} \cup \Delta_t$, where $\SelectionVec_{t}$ is the partial selection determining nodes in $Q_{t-1}$ by Property (A4) for step $t$. Note that since $Q_{t-1} \cap \CGLevel_{t}=\emptyset$, the corresponding
selections $\SelectionVec_{t}$ and $\SelectionVec_{t+1}$ agree.
%Following notation can be removed (defines the relation of similarity
%for selection vectors instead for selection matrices):
%
%($\FilterMatrix(\SelectionVec_{t}) \sim \FilterMatrix(\SelectionVec_{t-1}))$.

We now show that such $\Delta_t$ exists. This follows by the next claim.

\begin{claim}
\label{cl:aux1}
For every $t>1$, every entity $\Entity_j \in \CGLevel_t$ has an active repressor in $\overline{X}_{\beta}$, i.e.,
$\Repressors_j \cap NZ(\overline{X}_{\beta}) \neq \emptyset$.
\end{claim}
\Proof
We prove the claim by showing a slightly stronger statement, namely, that for every
$\Entity_j \in \CGLevel_t$ there exists an affector
$\Affectors_k \in \Repressors_j \cap \SelectionVec_{t}$.

For ease of analysis, let's focus on one specific $\GCGraph \in \mathfrak{G}(\SelectionVec_{t})$.
Since $\Entity_j \in \CGLevel_{t}$, it follows that there exists some $\Entity_i \in \CGLevel_{t-1}$
such that $(\Entity_i, \Entity_j) \in E(\GCGraph)$. Since $\SelectionVec_{t}$ determines $Q_{t-1}$ and $\Entity_i \in Q_{t-1}$, there exists a unique affector
$\Affectors_{\IndS(i)}=\SelectionVec_{t} \cap \Supporters_{i}$.
In addition, by Property (A1) for step $t$, $X_{\beta}(\Affectors_{\IndS(i)})>0$.
Therefore, since $\Entity_j$ is an immediate outgoing neighbor of $\Entity_i$, it holds by Eq. (\ref{eq:cg_condition}) that $\Affectors_{\IndS(i)} \in \Repressors_j$,
which establishes the claim.
\QED

We now complete the proof for the existence of $\Delta_t$.
By Claim \ref{cl:aux1}, each entity $\Entity_i \in \CGLevel_{t}$ has a strictly
positive repression, or, $\TotR(\overline{X}_{\beta}, \System)_{i}>0$.
Since $\overline{X}_{\beta}$ is feasible,
it follows by Fact \ref{fc:feasible_tots_totr} that also $\TotS(\overline{X}_{\beta}, \System)_{i}>0$. Therefore we get that for every
$\Entity_i \in \CGLevel_{t}$, there exists an affector
$\Affectors_{\IndS(i)} \in \Supporters_{i} \cap NZ(\overline{X}_{\beta})$.
Consequently, set $\Delta_t=\{\Affectors_{\IndS(i)} ~\mid~ \Entity_i \in \CGLevel_t \}$ and let $\SelectionVec_{t+1}= \SelectionVec_{t} \cup \Delta_t$.
%Recall that $\FilterMatrix(\SelectionVec_{t})$ is the matrix representation
%of $\SelectionVec_{t}$.
%\begin{lemma}
%\label{obs:sim_relation}
%$\FilterMatrix(\SelectionVec_{t}) \sim \FilterMatrix(\SelectionVec_{t+1})$.
%\end{lemma}
\begin{observation}
\label{obs:sim_relation}
$\SelectionVec_{t} \sim \SelectionVec_{t+1}$.
\end{observation}
\Proof
By definition, $\SelectionVec_{t}$ determines
$Q_{t-1}=\bigcup_{j=0}^{t-1} L_{j}(\GCGraph)$, for every $\GCGraph \in \mathfrak{G}(\SelectionVec_{t})$.
The selection $\SelectionVec_{t+1}$ consists of  $\SelectionVec_{t}$ and a new selection for the new layer $\CGLevel_{t}$ such that
$\CGLevel_{t} \cap  Q_{t-1}=\emptyset$ and therefore $\SelectionVec_{t}$ and
$\SelectionVec_{t+1}$ agree on their common part.
% Following notation can be removed:
%, i.e.,
%$\FilterMatrix(\SelectionVec_{t}) \sim \FilterMatrix(\SelectionVec_{t+1})$.
\QED

We now turn to prove Properties (A1)-(A4) for step $t+1$. Property (A1) follows immediately
by the construction of $\SelectionVec_{t+1}$. We next consider (A2).
\begin{claim}
\label{cl:aux4}
$\mathfrak{G}(\SelectionVec_{t+1}) \subseteq \mathfrak{G}(\SelectionVec_{t})$.
\end{claim}
\Proof
Consider some $\GCGraph \in \mathfrak{G}(\SelectionVec_{t+1})$. By Eq. (\ref{eq:graph_family_selection}), there exists
a complete selection $\SelectionVec^{*}$, where
$\GCGraph=\ConstraintsGraph_{\System(\SelectionVec^{*})}$, such that
$\SelectionVec^{*} \sim \SelectionVec_{t+1}$.
Recall that $\CGLevel_{i}=L_{i}(\GCGraph')$ for every $\GCGraph' \in \mathfrak{G}(\SelectionVec_{t})$ and
for every $i \in \{0, \ldots, t\}$
and that $Q_{t-1}=\bigcup_{i=0}^{t-1} \CGLevel_{i}$ and $Q_{t}=Q_{t-1} \cup \CGLevel_{t}$
where $Q_{t-1} \cap \CGLevel_{t}=\emptyset$. Therefore $Q_{t-1} \subset Q_{t}$. By the inductive assumption,
$\SelectionVec_{t}$ determines $Q_{t-1}$ and by construction
$\SelectionVec_{t+1}$ determines $Q_{t}$.
Combining all the above, Obs. \ref{obs:sim_relation},
$\SelectionVec_{t+1} \sim \SelectionVec_{t}$. Obs. \ref{obs:chain_sym} implies that
$\SelectionVec^{*} \sim \SelectionVec_{t}$.
Therefore, by Eq. (\ref{eq:graph_family_selection}) again, $\GCGraph \in \mathfrak{G}(\SelectionVec_{t})$.
\QED

Due to Claim \ref{cl:aux4}, and Properties (A2) and (A3) for step $t$ , Property (A2) follows for step $t+1$.
It is therefore possible to fix some $\GCGraph \in  \mathfrak{G}(\SelectionVec_{t+1})$ and define $\CGLevel_{i}=L_{i}(\GCGraph)$ for every
$i \in \{0,\ldots, t\}$ (by (A2) for $t+1$ this is well-defined)

We consider now Property (A3) and show that $L_{t+1}(\GCGraph_1)=L_{t+1}(\GCGraph_2)$ for every $\GCGraph_1, \GCGraph_2 \in \mathfrak{G}(\SelectionVec_{t+1})$.

For every graph $\GCGraph \in \mathfrak{G}(\SelectionVec_{t+1})$,
define $W(\GCGraph)$ as the set of all
immediate outgoing neighbors of $\CGLevel_{t}$ in $\GCGraph$, $W(\GCGraph) = \{\Entity_k \mid \exists \Entity_i \in \CGLevel_{t} \text{~such that~}
(\Entity_i,\Entity_k) \in E(\GCGraph)\}$.
\begin{observation}
\label{obs:w}
$W(\GCGraph_1)=W(\GCGraph_2)$ for every $\GCGraph_1, \GCGraph_2 \in \mathfrak{G}(\SelectionVec_{t+1})$.
\end{observation}
\Proof
Let $\GCGraph_1=\ConstraintsGraph_{\System(\SelectionVec_1)}$ and
$\GCGraph_2=\ConstraintsGraph_{\System(\SelectionVec_2)}$, where
$\SelectionVec_1, \SelectionVec_2$ correspond to complete legal selections.
Since $\GCGraph_1, \GCGraph_2 \in \mathfrak{G}(\SelectionVec_{t+1})$, it follows that
$\SelectionVec_1, \SelectionVec_2 \sim \SelectionVec_{t+1}$.
Since $\Delta_t$ determines $\CGLevel_t$, every entity $\Entity_i \in \CGLevel_t$
has the same unique supporter
$\Affectors_{\IndS(i)} \in \SelectionVec_{t+1} \cap \Supporters_i$ in both
$\SelectionVec_1, \SelectionVec_2$. By the definition of the constraint graph in Eq. (\ref{eq:cg_condition}),
it then follows that for graph $\GCGraph \in \mathfrak{G}(\SelectionVec_{t+1})$, the immediate outgoing neighbors of $\CGLevel_t$, $W(\GCGraph)$ are fully determined by the partial selection $\Delta_t$. The observation follows.
\QED

Hereafter, let $W=W(\GCGraph)$, $\GCGraph\in \mathfrak{G}(\SelectionVec_{t+1})$, be the set
of immediate neighbors of $\CGLevel_t$ in $\GCGraph$ (by Obs. \ref{obs:w}, this is well-defined).
Finally, note that $L_{t+1}(\GCGraph)=W \setminus \left(\bigcup_{i=1}^{t} L_{i}(\GCGraph) \right)$, for every
$\GCGraph \in \mathfrak{G}(\SelectionVec_{t+1})$. By Property (A2), $\CGLevel_{i}=L_{i}(\GCGraph)$ for every
$\GCGraph \in \mathfrak{G}(\SelectionVec_{t+1})$ and $i \in \{0, \ldots,t\}$. Hence, $L_{t+1}(\GCGraph)=W \setminus Q_t$ and by Obs. \ref{obs:w}, Property (A3) is established.

Finally, it remains to consider Property (A4). First, note that by Property (A2) and (A3) for step $t+1$, we get that $Q_{t}=Q_{t-1} \cup L_{t}(\GCGraph)$ for every $\GCGraph \in \mathfrak{G}(\SelectionVec_{t+1})$. By Property (A4) for step $t$ and Properties (A2) and (A3) for step $t+1$, it follows that the selection $\SelectionVec_{t+1}$ determines $Q_{t}$.

We now turn to discuss the stopping criterion. Let $t^{*}$ be the first time step $t$
where $\SelectionVec_{t^{*}}=\SelectionVec_{t^{*}-1}$. (Since $\SelectionVec_{t} \subseteq \SelectionVec_{t+1}$ for every $t\geq 0$, such $t^*$ exists). We then have the following.

\begin{lemma}
$|\SelectionVec_{t^{*}}|=n$ hence
$\System(\SelectionVec_{t^{*}})$ is a square system, and $\mathfrak{G}(\SelectionVec_{t^{*}})=\{\ConstraintsGraph_{\System(\SelectionVec_{t^{*}})}\}$,
\end{lemma}
\Proof
Recall that for every $i \in \{0, \ldots, t^{*}\}$, by Eq. (\ref{eq:graph_family_selection}), $\GCGraph' \in \mathfrak{G}(\SelectionVec_{i})$ represents a square system, and therefore by Obs. \ref{obs:reducible_graph_connected} it is strongly connected.
Fix some arbitrary $\GCGraph \in \mathfrak{G}(\SelectionVec_{t^{*}})$ and let $\CGLevel_i=L_i(\GCGraph)$ for every $i \in \{0, \ldots, t^{*}\}$ (By Property (A2) and (A3) this is well defined).
By Property (A4) it holds that the partial selection $\SelectionVec_{t^{*}-1}$ (resp., $\SelectionVec_{t^{*}}$) determines $Q_{t^*-2}$ (resp., $Q_{t^*-1}$).
As $\SelectionVec_{t^{*}-1}=\SelectionVec_{t^{*}}$, we have that $Q_{t^*-2}=Q_{t^*-1}$. Hence, $Q_{t^*-1} \setminus Q_{t^*-2}=\CGLevel_{t^*-1}=\emptyset$.
This implies that the BFS graph $BFS(\GCGraph, \Entity_{i_1})$ consists of $t^*-1$ levels $Q_{t^*-2}$. In addition, since $\GCGraph$ is strongly connected it follows that $Q_{t^*-2}=\EntitySet$.
By Property (A4), $\SelectionVec_{t^{*}}$ determines $Q_{t^*}$, hence $|\SelectionVec_{t^{*}}|=n$ meaning that $\SelectionVec_{t^{*}}$ is a complete selection, so $\System(\SelectionVec_{t^{*}})$
corresponds to a unique square system. Finally, since the $t^*-1$ layers of every $\GCGraph \in \mathfrak{G}(\SelectionVec_{t+1})$ are the same (Property (A2) and (A3)) and span all the entities it follows that $\mathfrak{G}(\SelectionVec_{t+1})$ consists of a single constraint graph, the lemma follows.
\QED
In summary, we end with a complete selection $\SelectionVec_{t^{*}}$ that spans the $n$ entities.
Every affector $\Affectors_k \in \SelectionVec_{t^{*}}$ is active and therefore
the constraint graph $\ConstraintsGraph_{\System(\SelectionVec_{t^{*}}}$ is
active in $\overline{X}_{\beta}$.
This establishes the following lemma.
\begin{lemma}
\label{lem:one_active}
For every feasible point $\overline{X}_{\beta}$ for Program (\ref{LP:Ext_Perron_convex}) and every active affector $\Affectors_{p_1}$ in $\overline{X}_{\beta}$,
there exists a complete selection $\SelectionVec^{*}$ for $\EntitySet$
such that $\SelectionVec^{*} \subseteq NZ(\overline{X}_{\beta})$,
hence the corresponding constraint subgraph $\ConstraintsGraph_{\System(\SelectionVec^{*})}$
is active in $\overline{X}_{\beta}$.
\end{lemma}

The following is an interesting implication.
\begin{corollary}
\label{cor:active}
For every feasible vector there exists an active spanning irreducible graph.
\end{corollary}
\Proof
Since every feasible vector is non-negative, there exists at least one active
affector in it, from which an active spanning irreducible graph can be constructed by Lemma \ref{lem:one_active}.
\QED

Finally, we are ready to complete the proof of Lemma \ref{cl:entity_one_nonzero}
for any irreducible system $\System$.
\Proof[Lemma \ref{cl:entity_one_nonzero}]
%%%%%%%%%%%%%%%%%%%%%%%%%%%%%%%% proof for the general case
Since $\sum_{i} X_{\beta}(i)>0$, it follows that there exists at least one
affector $\Affectors_{p_1}$ such that  $X_{\beta}(\Affectors_{p_1})>0$.
By Lemma \ref{lem:one_active}, there is a complete selection vector
$\SelectionVec^{*} \subseteq NZ( X_{\beta})$. The lemma follows.
\QED

We end this subsection by showing that every vertex
$\overline{X} \in V(\Polytope(\beta))$ is a $\WeakZero$ solution.
\begin{lemma}
\label{cl:weak_zero_polytope}
If the system of Program (\ref{LP:Ext_Perron_convex}) is irreducible, then
every $\overline{X} \in V(\Polytope(\beta))$ is a $\WeakZero$ solution for it, and in particular every optimal solution $\overline{X}^* \in V(\Polytope(\beta^*))$ is a $\WeakZeroStar$ solution.
\end{lemma}
\Proof
By  Claim \ref{cl:n_zero_polytope}, for every
$\overline{X} \in V(\Polytope(\beta))$, $|NZ(\overline{X})| \leq n+1$.
By Lemma \ref{cl:entity_one_nonzero}, for every $1 \leq i \leq n$,
$|NZ(\overline{X}) \cap \Supporters_i| \geq 1$. Therefore there exists
at most one entity $\Entity_{i}$ such that
$|NZ(\overline{X}) \cap \Supporters_i|=2$, and $|NZ(\overline{X}) \cap \Supporters_j|=1$ for every $j \neq i$, i.e., the solution is $\WeakZero$. The above holds for every $\beta \in (0, \beta^*]$. In particular, for the optimal $\beta$ value, $\beta^*$, it holds that $\overline{X}^* \in V(\Polytope(\beta^*))$ is a $\WeakZeroStar$ solution. \QED

%\subsubsection{$\ZeroStar$ solutions}
%\label{sec:zerostar}
\subsection{Existence of a $\ZeroStar$ solution}
\label{sec:zerostar}
In the previous section we established the fact that when $\System$ is irreducible, every vertex
$\overline{X} \in V(\Polytope(\beta))$ corresponds to an $\WeakZero$ solution for Program (\ref{LP:Ext_Perron_convex}).
In particular, this statement holds for $\beta=\beta^{*}(\System)$, the optimal $\beta$ for $\System$.
By the feasibility of the system for $\beta^{*}$, the corresponding polytope
is non-empty and bounded
(and each of its vertices is a $\WeakZeroStar$ solution),
hence there exist $\WeakZeroStar$ solutions for the problem.
The goal of this subsection is to establish the existence of a $\ZeroStar$
solution for the problem and thus complete the proof of Thm. \ref{thm:pf_ext}. In particular, we consider Program (\ref{LP:Ext_Perron_convex}) for an irreducible system $\System$ and $\beta=\beta^{*}$, i.e., the optimal value of Program (\ref{LP:Ext_Perron}) for $\System$, and show that {\em every}
optimal $\overline{X} \in V(\Polytope(\beta^{*}))$ solution is in fact a $\ZeroStar$ solution.

We begin by showing that for $\beta^{*}$,
the set of $n$ SR Inequalities (Eq. (\ref{eq:SR-convex}))
hold with equality for every optimal solution $\overline{X}^{*}$, including one that is not a $\WeakZeroStar$ solution.
\begin{lemma}
\label{lem:strict_equality}
If $\System=\langle \SupportersMatrix,\RepressorsMatrix \rangle$ is irreducible, then
$\RepressorsMatrix \cdot \overline{X}^{*} =
1/\beta^{*}(\System) \cdot \SupportersMatrix \cdot \overline{X}^{*}$ for every optimal solution $\overline{X}^{*}$ of Program (\ref{LP:Ext_Perron_convex}).
\end{lemma}
\Proof
%For definitions and notation see Appendix \ref{subsec:partial}.
Consider an irreducible system $\System$.
By Lemma \ref{cl:entity_one_nonzero}, every entity $\Entity_i$ has at least
one active supporter in $NZ(\overline{X}^{*})$. Select, for every $i$, one such
supporter $\Affectors_{\IndS(i)} \in \Supporters_i \cap NZ(\overline{X}^{*})$.
Let $\SelectionVec^{*}= \{\Affectors_{\IndS(i)} \mid 1\le i\le n \}$.
% $\SelectionVec^{*}= \bigcup_i \{\Affectors_{\IndS(i)}\}$.
By definition, $\SelectionVec^{*} \subseteq NZ(\overline{X}^{*})$. Also, by Claim \ref{cor:distinct}
the sets $\Supporters_i$ are disjoint. Therefore $\SelectionVec^{*}$ is a complete
selection
(i.e, for every $\Entity_i$, $|\Supporters_i \cap \SelectionVec^{*}|=1$), and hence $\System^{*}=\System(\SelectionVec^{*})$ is a square irreducible system.
Let $\GCGraph^*=\ConstraintsGraph_{\System^{*}}$ be the constraint graph of $\System^{*}$.
By Obs. \ref{obs:reducible_graph_connected}(a), $\GCGraph^*$ is strongly connected. In addition, since $\System^*$ has exactly one affector $\Affectors_{\IndS(i)}$ for every $\Entity_i$, and this affector is active, it follows that every edge
$e(\Entity_i,\Entity_j) \in E(\GCGraph^*)$ corresponds to an \emph{active} affector
in $\overline{X}^{*}$, i.e.,
$\Supporters_i \cap \Repressors_j \cap NZ(\overline{X}^{*}) \neq \emptyset$, and hence $\GCGraph^*$ is active.
%Lemma \ref{lem:one_active}  established the existence of a spanning directed subgraph $\GCGraph \subseteq \ConstraintsGraph_{\System}$ of the constraint graph, that has the following properties:
%(a) it is irreducible (strongly connected),
%(b) every directed edge $(i,j)$ in this subgraph is ``explained" by an active supporter
%in $\overline{X}^{*}$, i.e., $\left(\Supporters_i \cap \Repressors_j \right)\subseteq NZ(\overline{X}^{*})$.
%In other words, even considering only the active %affectors in
%$\overline{X}^{*}$ (and discarding the others), the constraint graph $\GCGraph^*$
%is guaranteed to be strongly connected
%(due to Lemma \ref{cl:entity_one_nonzero}), and %moreover,
%every directed edge is due to some affector with %positive entry in
%$\overline{X}^{*}$.
\par Therefore, for an edge $(v_i,v_j)$ in $\GCGraph^*$,
if we reduce the power of the active supporter of $v_i$ which,
by the definition of $\GCGraph^*$ (see Eq. (\ref{eq:cg_condition})) is a repressor of $v_j$, then $v_j's$ inequality
can be made strict. Such reduction makes sense
only because we consider active affectors. This intuition
is next used in order to prove the lemma.
For a feasible solution $\overline{X}$ of Program (\ref{LP:Ext_Perron_convex}) and vaule $\beta$, let us formulate the SR constraints in terms of total support and total repression as in (Eq. (\ref{eq:SR_ind})) , and let
\begin{equation}
\label{eq:residual}
R_{i}(\overline{X}) = 1/\beta \cdot
\TotS(\overline{X})_{i}-\TotR(\overline{X})_{i}
\end{equation}
be the residual amount of the $i'th$ SR constraint
of (\ref{eq:SR_ind})(hence $R_{i}(\overline{X})>0$ implies strict inequality on the $i$th constraint with $\overline{X}$).
Then the lemma claims that for the optimal solution $\overline{X}^*$ and $\beta^*$, $R_{i}(\overline{X}^*)=0$ for every $i$.
\par Assume, toward contradiction, that there exists at least one entity, w.l.o.g. $\Entity_{0}$, for which $R_{0}(\overline{X}^*)>0$.
In what follows, we gradually construct a new assignment $\overline{X}^{**}$
that achieves a strictly positive residue $R_{i}(\overline{X}^{**})>0$ , or, a strict inequality in the SR constraint of Eq. (\ref{eq:SR_ind}),
for all $\Entity_i \in \EntitySet$.
Clearly, if all SR constraints are satisfied with strict inequality,
then there exists some larger $\beta^{**}>\beta^{*}(\System)$ that still satisfies all the constraints, in contradiction to the optimality of $\beta^{*}(\System)$.

To construct $\overline{X}^{**}$, we trace paths of
influence in the strongly connected (and active) constraint graph $\GCGraph^*$.
Think of $\Entity_0$ as the root, and let $L_{j}(\GCGraph^*)$ be the $j^{th}$ level of
$\BFS(\GCGraph^*, \Entity_{0})$ (with $L_0=\{\Entity_{0}\}$).
Let $Q_{-1}=\emptyset$, and $Q_{t}=\bigcup_{i=0}^{t} L_{i}(\GCGraph^*)$ for $t \geq 0$.
Let $\SelectionVec_{t}=\{\Affectors_{\IndS(i)} ~\mid~ \Entity_{i} \in Q_{t-1}\} \subseteq \SelectionVec^{*}$ be the partial selection
determining the entities in $Q_{t-1}$. I.e., $|\SelectionVec_{t}|=|Q_{t-1}|$ and
for every $\Entity_i \in Q_{t-1}$, $|\SelectionVec_{t} \cap \Supporters_{i}|=1$.

The process of constructing $\overline{X}^{**}$ consists of $d$ steps,
where $d$ is the depth of $\BFS(\GCGraph^*, \Entity_{0})$.
At step $t$, we are given $\overline{X}_{t-1}$ and use it to construct
$\overline{X}_{t}$.
Essentially, $\overline{X}_{t}$ should satisfy the following properties.
\begin{description}
\item{(B1)}
The set of SR inequalities corresponding to $Q_{t-1}$ entities hold with strict
inequality with $\overline{X}_{t}$. That is, for every $\Entity_i \in Q_{t-1}$, $R_{i}(\overline{X}_{t})>0$, i.e.,
$$ 1/\beta^{*} \cdot \TotS(\overline{X}_t)_{i} ~>~
\TotR(\overline{X}_t)_{i} ~. $$
\item{(B2)}
$\overline{X}_{t}$ is an optimal solution, i.e., it satisfies Program
(\ref{LP:Ext_Perron}) with $\beta^{*}(\System)$.
\item{(B3)}
$X_{t}(\Affectors)=X^{*}(\Affectors)$
for every $\Affectors \notin \SelectionVec_{t}$ and
$X_{t}(\Affectors)<X^{*}(\Affectors)$
for every $\Affectors \in \SelectionVec_{t}$.
\end{description}

Let us now describe the construction process in more detail.
Let $\overline{X}_0=\overline{X}^{*}$.
Consider step $t=1$ and recall that $R_{0}(\overline{X}_0)>0$.
Let $\Affectors_{k_0}$ be the active supporter of
$\Entity_{0}$, i.e.,
$\Affectors_{k_0} \in \Supporters_{0} \cap \SelectionVec^{*}$.
Then it is possible to slightly reduce the value of
$\Affectors_{k_0}$ in $\overline{X}_0$ while still maintaining feasibility, yielding $\overline{X}_1$.
Formally, let
$X_{1}(\Affectors_{k_0}) = X_0(\Affectors_{k_0})-
\min\{X_0(\Affectors_{k_0}),R_{0}(\overline{X}_0)\}/2$ and leave the rest
of the entries unchanged, i.e.,
$X_{1}(\Affectors_{k})=X^{*}(\Affectors_{k})$ for every other $k \neq k_0$.
We now show that Properties (B1)-(B3) are satisfied for $t\in \{0,1\}$ and then
proceed to consider the construction of $\overline{X}_{t}$ for $t>1$.
Since $L_{0}(\GCGraph^*)=\{\Entity_{0}\}$, and $Q_{-1}=\emptyset$, also $\SelectionVec_{0}=\emptyset$, so (B1) holds vacuously, and (B2) and (B3) follow by the fact that
$\overline{X}_{0}=\overline{X}^*$.
Next, consider $\overline{X}_{1}$. By the irreducibility of the system
(in particular, see Cl. \ref{cor:distinct}), since  only $\Affectors_{k_0}$
was reduced in $\overline{X}_1$ (compared to $\overline{X}^{*}$), only the
constraint of $\Entity_0$ could have been damaged (i.e., become unsatisfied).
Yet, it is easy to verify that the constraint of $\Entity_{0}$ still holds
with strict inequality for $\overline{X}_{1}$, so Property (B2) holds. As $Q_{0}=\{\Entity_0\}$, Property (B1) needs to be verified only for $\Entity_0$, and indeed the new value of $X_1(\Affectors_{k_0})$ ensures $R_0(\overline{X}_{1})>0$, so (B1) is satisfied. Finally, $\SelectionVec_{1}=\{\Affectors_{k_0}\}$, and Property (B3) checks out as well.

Next, we describe the general construction step. Assume that we are given solution
$\overline{X}_{r}$ satisfying Properties (B1)-(B3)
for each $r\leq t$. We now describe the construction of $\overline{X}_{t+1}$
and then show that it satisfies the desired properties.
We begin by showing that the set of SR inequalities of Eq. (\ref{eq:SR_ind}) on the entities $\Entity_i$ in $L_{t}(\GCGraph^*)$
hold with strict inequality with $\overline{X}_{t}$.

\begin{claim}
\label{cl:strict_inequality_induc}
$R_j(\overline{X}_{t})>0$, or, $\TotR(\overline{X}_{t})_{j} < 1/\beta^{*} \cdot \TotS(\overline{X}_{t})_{j}$,
for every entity $\Entity_j \in L_{t}(\GCGraph^*)$.
\end{claim}
\Proof
Consider some $\Entity_j \in L_{t}(\GCGraph^*)$. By definition of $L_{t}(\GCGraph^*)$,
there exists an entity $\Entity_i \in L_{t-1}(\GCGraph^*)$ such that
$e(i,j) \in E(\GCGraph^*)$.
Since $\Entity_{i} \in Q_{t-1}$ and $\SelectionVec_{t}$ is a partial selection determining
$Q_{t-1}$, a (unique) supporter
$\Affectors_{\IndS(i)} \in \SelectionVec_{t} \cap \Supporters_{i}$
is guaranteed to exist.
By the definition of $\GCGraph^*$, $e(\Entity_i,\Entity_j) \in E(\GCGraph^*)$
implies that $\Affectors_{\IndS(i)} \in \Repressors_{j}$.
Finally, note that by Property (B3),
$X_{t}(\Affectors_{\IndS(i)})<X^*(\Affectors_{\IndS(i)})$ and
$X_{t}(\Affectors)=X^{*}(\Affectors)=X_{t-1}(\Affectors)$ for every
 $\Affectors \in \Supporters_j$ (since $\SelectionVec_{t} \cap \Supporters_j =\emptyset$). I.e.,
\begin{equation}
\label{eq:tot_sup_rep_ineq}
\TotS(\overline{X}_{t})_{j}=\TotS(\overline{X}_{t-1})_{j} \text{~and~}
\TotR(\overline{X}_{t})_{j}<\TotR(\overline{X}_{t-1})_{j},
\end{equation}
which implies by Eq. (\ref{eq:SR_ind}) that
\begin{equation}
\label{eq:residual_step}
R_j(\overline{X}_{t-1})<R_j(\overline{X}_{t})~.
\end{equation}

By the optimality of $\overline{X}_{t-1}$
(Property (B2) for step $t-1$), we have that
$R_j(\overline{X}_{t-1}) \geq 0$.
Combining this with Eq. (\ref{eq:residual_step}),
$0 \leq R_j(\overline{X}_{t-1})<R_j(\overline{X}_{t})$,
which establishes the claim for $\Entity_j$. The same argument can be applied
for every $\Entity_j \in L_{t}(\GCGraph^*)$, thus the claim is established.
\QED

Let $\Delta_t \subseteq \SelectionVec^{*}$ be the partial selection that determines $L_{t}(\GCGraph^*)$. In the solution $\overline{X}_{t+1}$, only the entries of $\Delta_t$ have been reduced and the other entries remain as in $\overline{X}_{t}$.
Recall that by construction, $\SelectionVec^{*} \subseteq NZ(\overline{X}^{*})$
and therefore also $\SelectionVec^{*} \subseteq NZ(\overline{X}_{t})$.
By Claim \ref{cl:strict_inequality_induc}, the constraints of $L_{t}(\GCGraph^*)$ nodes
hold with strict inequality, and therefore it is possible to slightly reduce
the value of their positive supporters while still maintaining the strict
inequality (although with a lower residue). Formally, for every
$\Entity_k \in L_{t}(\GCGraph^*)$, consider its unique supporter in $\Delta_t$,
$\Affectors_{i_k}\in \Delta_t \cap \Supporters_k$.
By Claim \ref{cl:strict_inequality_induc}, $R_{k}(\overline{X}_{t})>0$.
Set $X_{t+1}(\Affectors_{i_k}) =
X_{t}(\Affectors_{i_k})-\min(X_{t}(\Affectors_{i_k}),R_{k}(\overline{X}_{t}))/2$.
In addition, $X_{t+1}(\Affectors_{i_k})=X_{t}(\Affectors_{i_k})$ for every other
supporter $\Affectors_{i_k} \notin \Delta_t$.

It remains to show that $\overline{X}_{t+1}$ satisfies the Properties (B1)-(B3).
(B1) follows by construction. To see (B2), note that since
$\Supporters_{i} \cap \Supporters_{j} = \emptyset$ for every
$\Entity_i, \Entity_j \in \EntitySet$, only the constraints of $L_{t}(\GCGraph^*)$ nodes
might have been violated by the new solution $\overline{X}_{t+1}$.
Formally, $\TotS(\overline{X}_{t+1})_{i}=\TotS(\overline{X}_{t})_{i}$ and
$\TotR(\overline{X}_{t+1})_{i} \leq \TotR(\overline{X}_{t})_{i}$ for every
$\Entity_i \notin L_{t}(\GCGraph^*)$. Although, for $\Entity_i \in L_{t}(\GCGraph^*)$,
we get that $\TotS(\overline{X}_{t+1})_{i}<\TotS(\overline{X}_{t})_{i}$
(yet $\TotR(\overline{X}_{t+1})_{i} = \TotR(\overline{X}_{t})_{i}$),
this reduction in the total support of $L_{t}(\GCGraph^*)$ nodes was performed
in a controlled manner,
guaranteeing that the corresponding $L_{t}(\GCGraph^*)$ inequalities
hold with \emph{strict} inequality. Finally, (B3) follows immediately.
After $d+1$ steps, by Property (B1) all inequalities hold with strict inequality
(as $Q_{d}=\EntitySet$) with the solution $\overline{X}_{d+1}$.
Thus, it is possible to find some $\beta^{**}>\beta^{*}(\System)$ that would
contradict the optimally of $\beta^{*}$.
Formally, let $R^{*}=\min R_{i}(\overline{X}_{d+1})$. Since $R^{*}>0$,
we get that $\overline{X}_{d+1}$ is feasible with
$\beta^{**}=\beta^{*}(\System)+R^{*}>\beta^{*}(\System)$,
contradicting the optimally of $\beta^{*}(\System)$.
Lemma \ref{lem:strict_equality} follows.
\QED

We proceed by considering a vertex of
$\overline{X}^{*} \in V(\Polytope(\beta^{*}))$.
By Lemma \ref{cl:weak_zero_polytope}, $\overline{X}^{*}$ is a $\WeakZeroStar$ solution.
%In Sec. \ref{sec:geoproofs}, we establish the following and complete the proof for Thm. \ref{thm:pf_ext}.
To complete the proof of Thm. \ref{thm:pf_ext}, we have to prove that it is a $\ZeroStar$ solution.
To do that, we first transform $\System$ into a weakly square system $\WeakSystem$. First, if $m=n+1$, then the system is already weak.
Otherwise, without loss of generality, let the $i^{th}$ entry in
$\overline{X}^{*}$ correspond to $\Affectors_i$ where
$\Affectors_i=NZ(\overline{X}^{*}) \cap \Supporters_i$ for $i \in \{1,\ldots, n-1\}$ and
the $n^{th}$ and $(n+1)^{st}$ entries correspond to $\Affectors_n$ and $\Affectors_{n+1}$
respectively such that
$\{\Affectors_n, \Affectors_{n+1}\}=NZ(\overline{X}^{*}) \cap \Supporters_n$.
It then follows that $X^{*}(i) \neq 0$ for every $i \in \{1,\ldots,n+1\}$ and
$X^{*}(i) = 0$ for every $i \in \{n+2, \ldots, m\}$.
Let $\overline{X}^{**}=\left(X^{*}(1),\ldots, X^{*}(n+1)\right)$.
Let $\WeakSupportersMatrix \in \R^{n \times (n+1)}$ where
$\WeakSupportersMatrix(i,j)=\SupportersMatrix(i,j)$ for every
$i \in \{1,\ldots, n\}$ and every $j \in \{1,\ldots, n+1\}$, and define $\WeakRepressorsMatrix$ analogously.
From now on, we restrict attention to the weakly square system
$\WeakSystem=\langle \WeakSupportersMatrix,\WeakRepressorsMatrix\rangle$ where
$|\Supporters_{n}|=2$. Note that this system results from $\System$
by discarding the corresponding entries of
$\Affectors \setminus NZ(\overline{X}^{*})$.
Therefore,
$\beta^{*}(\System)=\beta^{*}(\WeakSystem)$.
Let $\SupportersMatrix_{n-1}$ correspond to the upper left
$(n-1) \times (n-1)$ submatrix of $\WeakSupportersMatrix$.
Let $\SupportersMatrix_{n}$ be obtained from $\WeakSupportersMatrix$ by
removing the $(n+1)^{st}$ column. Finally, $\SupportersMatrix_{n+1}$
is obtained from $\WeakSupportersMatrix$ by removing the $n^{th}$ column.
The matrices
$\RepressorsMatrix_{n-1},\RepressorsMatrix_{n},\RepressorsMatrix_{n+1}$
are defined analogously.

To study the weakly square system $\WeakSystem$, we consider the following three
\emph{square} systems:
\begin{eqnarray}
\label{eqn:ssystem_def}
\System_{n-1} &=& \langle \SupportersMatrix_{n-1}, \RepressorsMatrix_{n-1}\rangle~,
\\
\System_{n}   &=& \langle \SupportersMatrix_{n}, \RepressorsMatrix_{n}\rangle~, \nonumber
\\
\System_{n+1} &=& \langle \SupportersMatrix_{n+1}, \RepressorsMatrix_{n+1}\rangle~. \nonumber
\end{eqnarray}
Note that
%a feasible solution $\overline{X}_{n-1}$ for the system
%$\System_{n-1}$ corresponds to a feasible solution $\overline{X}_{w}$ for the
%system $\WeakSystem$ by setting $X_{w}(\Affectors_i)= X_{n-1}(\Affectors_i)$,
%for every $i \neq n, n+1$ and setting $X_{w}(\Affectors_n)=0$ and
%$X_{w}(\Affectors_{n+1})=0$. In the same manner,
a feasible solution
$\overline{X}_{n+b}$ for the system $\System_{n+b}$, for $b \in \{0,1\}$,
corresponds to a feasible solution for  $\WeakSystem$ by setting
$X_{w}(\Affectors_j)=X_{n+b}(\Affectors_j)$ for every $j \neq n+(1-b)$
and $X_{w}(\Affectors_{n+(1-b)})=0$.
For ease of notation, let
$\CP_{n}(\lambda) = \CP(Z(\System_{n}), \lambda)$,
$\CP_{n+1}(\lambda) = \CP(Z(\System_{n+1}), \lambda)$ and
$\CP_{n-1}(\lambda) = \CP(Z(\System_{n-1}), \lambda)$,
where $\CP$ is the characteristic polynomial defined in Eq. (\ref{eq:CP}).
Let $\beta^{*}_{n+b}=\beta^*(\System_{n+b})$ be the optimal
value of Program (\ref{LP:Ext_Perron}) for the system
$\System_{n+b}$.
Let $\beta^*=\beta^*(\System)$ and let
\begin{eqnarray*}
\lambda^{*}&=&1/\beta^{*} ,\\
\lambda^{*}_{n+b}&=&1/\beta^{*}_{n+b},
\mbox{~for~}  b \in \{-1,0,1\}~.
\end{eqnarray*}

\begin{claim}
\label{cl:n_1_beta_optimal}
$\max\{\beta^{*}_{n}, \beta^{*}_{n+1}\} \leq \beta^{*} < \beta^{*}_{n-1}$.
\end{claim}
\Proof
The left inequality follows as any optimal solution $\overline{X}^{*}$ for
$\System_{n}$ (respectively, $\System_{n+1}$) can be achieved in the weakly square system
$\WeakSystem$ by setting $X^{*}(\Affectors_{n+1})=0$
(resp., $X^{*}(\Affectors_{n})=0$).
\par Assume towards contradiction that $\beta^*=\beta^*_{n-1}$ and let $\overline{X}'$ be the optimal solution for $\WeakSystem$.

By Lemma \ref{cl:entity_one_nonzero}, it holds that $X'(\Affectors_{n})+ X'(\Affectors_{n+1}) >0$.
Without loss of generality, assume that  $X'(\Affectors_{n})>0$.
By Obs. \ref{obs:reducible_graph_connected}(a) and
the irreducibility of $\WeakSystem$, $\Entity_{n}$ is strongly connected to the rest of the graph for every selection of one of its two supporters. Thus there exists at least one entity $\Entity_{j}$, $j \in [1,n-1]$
such that $\Affectors_{n} \in \Repressors_{j}$.

Let $\overline{X}'' \in \R^{n-1}$ be obtained by taking the values of the first $n-1$ affectors as in $\overline{X}'$ and discarding the values of $\Affectors_{n}$ and $\Affectors_{n+1}$. We have the following.
\begin{eqnarray}
\label{eqn:totstotr}
\TotS(\overline{X}'', \System_{n-1})_{j}=\TotS(\overline{X}', \WeakSystem)_{j} \mbox{~and~}
\TotR(\overline{X}'', \System_{n-1})_{j}<\TotR(\overline{X}', \WeakSystem)_{j}~,
\end{eqnarray}
where strict inequality follows by the assumption that $X'(\Affectors_{n})>0$ and $\Affectors_{n}$ is a repressor of $\Entity_j$.
Since $\overline{X}'$ is an optimal solution for the system $\WeakSystem$, by Lemma \ref{lem:strict_equality}, it holds that $\TotS(\overline{X}', \WeakSystem)_{j}=\TotR(\overline{X}', \WeakSystem)_{j}$. Combining with Eq. (\ref{eqn:totstotr}), we get that $\TotS(\overline{X}'', \System_{n-1})_{j}<\TotR(\overline{X}'', \System_{n-1})_{j}$. Since $\overline{X}''$ is an optimal solution for $\System_{n-1}$, we end with contradiction to Lemma \ref{lem:strict_equality}, concluding that $\beta^* <\beta^*_{n-1}$. The claim follows.
\QED

Our goal in this section is to show that the optimal $\beta^{*}$ value for
$\WeakSystem$ can be achieved by setting either $X^{*}(\Affectors_{n})=0$ or
$X^{*}(\Affectors_{n+1})=0$, essentially showing that the optimal
$\WeakZeroStar$ solution corresponds to a $\ZeroStar$ solution.
This is formalized in the following lemma.
\begin{lemma}
\label{thm:0_solution}
$\beta^{*}=\max\{\beta^{*}_{n}, \beta^{*}_{n+1}\}$.
\end{lemma}

The following observation holds for every $b \in \{-1,0,1\}$ and follows
immediately by the definitions of feasibility and irreducibility and
the \PFT~\ref{thm:pf_full}.

\begin{observation}
\label{obs:perron_application}
\begin{description}
\item{(1)}
$\lambda^*_{n+b}>0$ is the maximal eigenvalue of $Z(\System_{n+b})$.
\item{(2)}
For an irreducible system $\System$, $\lambda^*_{n+b}=1/\beta^*_{n+b}$.
\item{(3)}
If the system is feasible then $\lambda^*_{n+b}>0$.
\end{description}
\end{observation}
For a square system $\System \in \SquareSystemFamily$, let $W^1$ be a modified
form of the matrix $Z$, defined as follows.
$$W^1(\System, \beta) ~=~ Z(\System)-1/\beta \cdot I ~~~\text{for}~~~ \beta \in (0, \beta^{*}].$$
More explicitly,
$$W^{1}(\System, \beta)_{i,j} ~=~
\begin{cases}
-1/\beta, & \text{if $i=j$;}\\
-g(\Entity_i,\Affectors_j)/g(i,i), & \text{otherwise.}
\end{cases}
$$
Clearly, $W^1(\System, \beta)$ cannot be defined for a nonsquare
system $\System \notin \SquareSystemFamily$. Instead, a generalization $W^2$ of $W^1$
for any (nonsquare) $m \geq n$ system $\System$ is given by
$$W^2(\System, \beta) ~=~
\RepressorsMatrix- 1/\beta \cdot \SupportersMatrix,
~~~\text{for}~~~ \beta \in (0, \beta^{*}],$$
or explicitly,
$$W^{2}(\System, \beta)_{i,j} ~=~
\begin{cases}
-g(i,i)/\beta, & \text{if $i=j$;}\\
-g(\Entity_i,\Affectors_j), & \text{otherwise.}
\end{cases}
$$
Note that if $\overline{X}_{\beta}$ is a feasible solution for $\System$,
then $W^{2}(\System, \beta) \cdot \overline{X}_{\beta} \leq 0$.
If $\System \in \SquareSystemFamily$, it also holds that
$W^{1}(\System, \beta) \cdot \overline{X}_{\beta} \leq 0$.

For $\System \in \SquareSystemFamily$, where both $W^{1}(\System, \beta)$ and
$W^{2}(\System, \beta)$ are well-defined, the following connection
becomes useful in our later argument. Recall that $\CP(Z(\System),t)$
is the characteristic polynomial of $Z(\System)$ (see Eq. (\ref{eq:CP})).

\begin{observation}
\label{obs:x_y_relation}
For a square system $\System$,\\
(a) $\det(-W^{1}(\System, \beta))=\CP(Z(\System),1/\beta)$ and \\
(b) $\det(-W^{2}(\System, \beta)) ~=~
\CP(Z(\System),1/\beta) \cdot \prod_{i=1}^{n} g(i,i)$.
\end{observation}
\Proof
The observation follows immediately by noting that
$W^{1}(\System, \beta)_{i,j}=W^{2}(\System, \beta)_{i,j} \cdot g(i,i)$ for every $i$ and $j$, and by Eq. (\ref{eq:CP}).
\QED
The next equality plays a key role in our analysis.
\begin{lemma}
\label{lem:P_n_n+1}
%\begin{eqnarray*}
$\displaystyle \frac{g(n,n) \cdot X^{*}(n) \cdot \CP_{n}(\lambda^{*})}
{\CP_{n-1}(\lambda^{*})} +
\frac{g(n,n+1) \cdot X^{*}(n+1) \cdot \CP_{n+1}(\lambda^{*})}
{\CP_{n-1}(\lambda^{*})} = 0.$
%\end{eqnarray*}
\end{lemma}
\Proof
%We now complete the proof for Lemma \ref{lem:P_n_n+1}.
By Lemma \ref{lem:strict_equality}, it follows that
$-W^{2}(\WeakSystem,\beta^{*}) \cdot \overline{X}^{*}=0$,
or
$$
\begin{pmatrix}
g(1,1)/\beta^{*} &  g(1,2)  & \ldots &
g(1,n)&
g(1,n+1)\\
g(2,1) &  g(2,2)/\beta^{*}  & \ldots &
g(2,n)&
g(2,n+1)\\
\vdots & \ldots & \ldots & \vdots\\
g(n,1) &  g(n,2)  & \ldots &
g(n,n)/\beta^{*}&
g(n,n+1)/\beta^{*}\\
\end{pmatrix}
\cdot
\begin{pmatrix}
X^{*}(1)\\
X^{*}(2)\\
\vdots\\
X^{*}(n)\\
X^{*}(n+1)
\end{pmatrix}
=
\begin{pmatrix}
0\\
\vdots\\
0\\
0
\end{pmatrix}
$$
Next, we need to apply Claim \ref{cl:cramer_non_square}(b).
%Eq. (\ref{eq:non_square_cramer_n}),
To do that, we first need to verify that
$W^{2}(\System_{n-1},\beta^{*})$, i.e., the $(n-1) \times (n-1)$ upper left
submatrix of $W^{2}(\WeakSystem,\beta^{*})$, is nonsingular.
This follows by noting that $\lambda^{*} \in \R_{>0}$ and
by Claim \ref{cl:n_1_beta_optimal}, $\lambda^{*}>\lambda^{*}_{n-1}$.
Moreover, note that $\lambda^{*}_{n-1}$ is the largest real root of
$\CP_{n-1}(\lambda)$, hence
\begin{equation}
\label{eq:cpnz}
\CP_{n-1}(\lambda^{*}) \neq 0~.
\end{equation}
Combining with Obs. \ref{obs:x_y_relation}(b), it follows that
$\det(-W^{2}(\System_{n-1}, \beta^{*})) \neq 0$ or that
$W^{2}(\System_{n-1}, \beta^{*})$ is nonsingular.

Now we can safely apply Claim \ref{cl:cramer_non_square}(b),
%Eq. (\ref{eq:non_square_cramer_n}),
yielding
\begin{eqnarray*}
\label{eq:Cramer_SINR_mid}
X^{*}(n) \cdot \frac{\det \left(-W^{2}(\System_{n},\beta^{*}) \right)}
{\det \left(-W^{2}(\System_{n-1},\beta^{*}) \right)}+ X^{*}(n+1) \cdot \frac{\det (-W^{2}(\System_{n+1},\beta^{*}))}
{\det \left(-W^{2}(\System_{n-1},\beta^{*}) \right)} = 0~.
\end{eqnarray*}
By plugging Obs. \ref{obs:x_y_relation}(b) and simplifying, the lemma follows.
\QED

Our work plan from this point on is as follows. We first define a range of
`candidate' values for $\beta^{*}$. Essentially, our interest is in
\emph{real} positive $\beta^{*}$.
Recall that $Z(\WeakSystem), Z(\System_{n})$ and $Z(\System_{n+1})$ are nonnegative
irreducible square matrices and therefore Theorem \ref{thm:pf_full} can be applied
throughout the analysis.
Without loss of generality, assume that $\beta^{*}_{n} \geq \beta^{*}_{n+1}$
(and thus $\lambda_{n}^{*}\leq \lambda_{n+1}^{*}$) and let
$Range_{\beta^*}=(\beta^{*}_{n}, \beta^{*}_{n-1}) \subseteq \R_{>0}$.
Let the corresponding range of $\lambda^{*}$ be
\begin{equation}
\label{eq:lambda_range}
Range_{\lambda^{*}}=(\lambda^{*}_{n-1},\lambda^{*}_{n})=(1/\beta^{*}_{n-1}, 1/\beta^{*}_{n}).
\end{equation}
To complete the proof for  Lemma \ref{thm:0_solution} we assume, towards contradiction, that
$\beta^{*} > \beta^{*}_{n}$.
According to Claim \ref{cl:n_1_beta_optimal} and the fact that $\beta^{*} \neq \beta^{*}_{n}$,
it then follows that $\beta^{*} < \beta^{*}_{n}$, $\lambda^*< \lambda^*_{n}, \lambda^*_{n+1}$ and hence $\CP_{n}(\lambda^*),\CP_{n+1}(\lambda^*)\neq 0$.

In addition, $\beta^{*} \in Range_{\beta^*}$.
Note that since $Range_{\beta^{*}} \subseteq \R_{>0}$, also
$Range_{\lambda^{*}}\subseteq \R_{>0}$, namely, the corresponding $\lambda^*$ is real and positive as well.
This is important mainly in the context of nonnegative irreducible matrices
$Z(\System')$ for $\System' \in \SquareSystemFamily$.
In contrast to nonnegative primitive matrices
(where $\Period=1$) for irreducible matrices, such as $Z(\System')$,
by Thm. \ref{thm:pf_full} there are $\Period \geq 1$ eigenvalues,
$\lambda_i \in \EigenValue(\System')$, for which
$|\lambda_i|=\PFEigenValue(\System')$.
However, note that only one of these, namely, $\PFEigenValue(\System')$,
might belong to $Range_{\lambda^{*}}\subseteq \R_{>0}$.
(This follows as by Thm. \ref{thm:pf_full}, every other such $\lambda_i$ is either real but negative
or with a nonzero complex component).

Fix $b \in \{-1,0,1\}$ and let $k_{n+b}$ be the number of real and positive
eigenvalues of $Z(\System_{n+b})$. Let
$0<\lambda_{n+b}^{1} \leq  \lambda_{n+b}^{2} \ldots \leq \lambda_{n+b}^{k_{n+b}}$
be the ordered set of \emph{real and positive} eigenvalues for
$Z(\System_{n+b})$, i.e., real positive roots of $\CP_{n+b}(\lambda)$.
Note that $\lambda_{n+b}^{k_{n+b}}=\lambda_{n+b}^{*}$.
By Theorem \ref{thm:pf_full}, we have that for every $b \in \{-1,0,1\}$
\\
(a) $\lambda_{n+b}^{*} \in \R_{>0}$, and
\\
(b) $\lambda_{n+b}^{*} > |\lambda_{n+b}^{p}|$, $p \in \{1,\ldots, k_{n+b}-1\}$.

We proceed by showing that the potential range for $\lambda^{*}$, namely,
$Range_{\lambda^{*}}$, can contain no root of $\CP_{n}(\lambda)$ and
$\CP_{n+1}(\lambda)$.
Since $Range_{\lambda^{*}}$ is real and positive, it is sufficient to consider
only real and positive roots of $\CP_{n}(\lambda)$ and $\CP_{n+1}(\lambda)$
(or real and positive eigenvalues of $Z(\System_{n})$ and $Z(\System_{n+1})$).

%%%%%%%%%%%%%%%
\begin{figure}[htb]
%\label{fig:zero_star_proof}
\begin{center}
\includegraphics[scale=1]{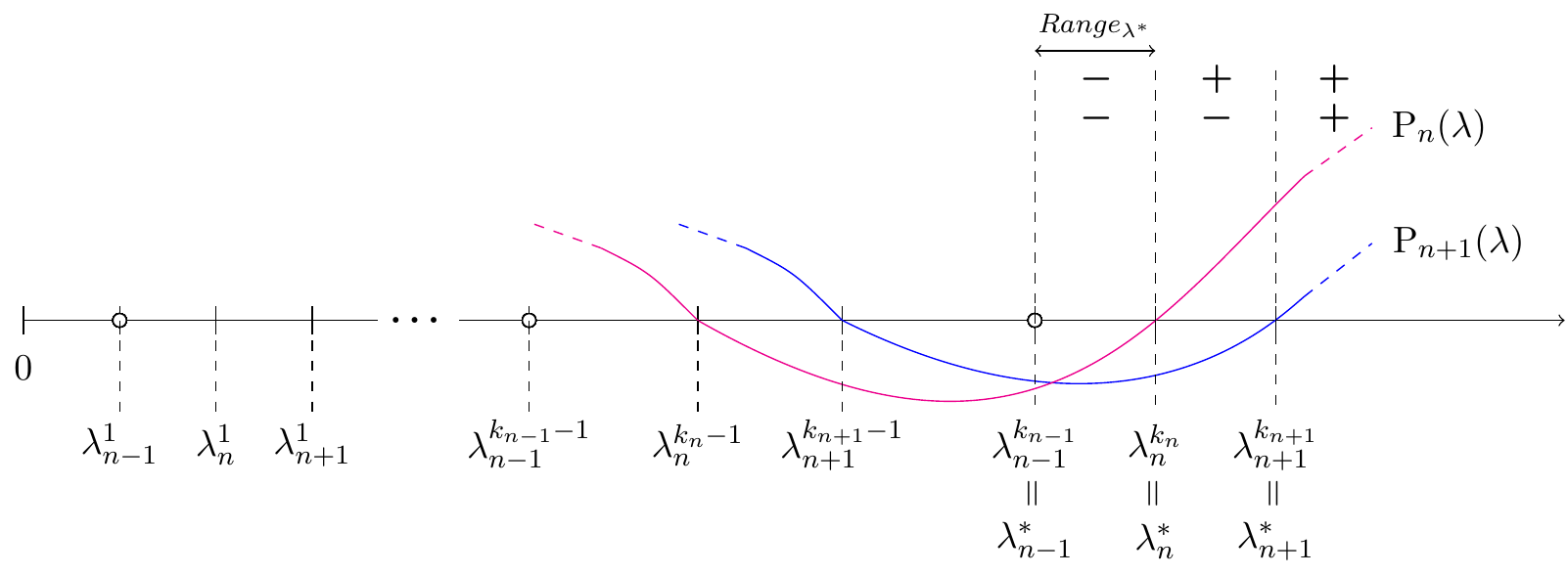}
\caption{Real positive roots of $\CP_{n+1}(\lambda)$, $\CP_{n}(\lambda)$,
and $\CP_{n-1}(\lambda)$. Each eigenvalue sequence $\lambda^p_{n+b}$ is ordered increasingly, but the relative ordering in which the sequences are merged in the figure is arbitrary, except $\lambda^{*}_{n-1}, \lambda^*_{n}$ and $\lambda^*_{n+1}$. Note that in the range $Range_{\lambda^{*}}$ that are no roots of $\CP_{n+b}(\lambda)$ for $b \in \{-1,0,1\}$.
\label{fig:eigenval}}
\end{center}
\end{figure}
%%%%%%%%%%%%%%%

\begin{claim}
\label{cl:no_root_in_range}
$\lambda_{n}^{p_0}, \lambda_{n+1}^{p_1} \notin Range_{\lambda^{*}}$ for every real
$\lambda_{n}^{p_0}, \lambda_{n+1}^{p_1}$, for $p_0 <k_{n}, p_1<k_{n+1}$.
\end{claim}
\Proof
Note that $Z(\System_{n-1})$ is the principal $(n-1)$ minor of both
$Z(\System_{n})$ and $Z(\System_{n+1})$.  By the separation theorem of Hall and Porsching, see Lemma. \ref{lem:sep_thm}, we get that
$\lambda_{n}^{p_0}, \lambda_{n+1}^{p_1}  \leq \lambda_{n-1}^{*}$ for every
$p_0 <k_{n}$ and $p_1 < k_{n+1}$, concluding by Eq. (\ref{eq:lambda_range}) that
$\lambda_{n}^{p_0}, \lambda_{n+1}^{p_1} \notin Range_{\lambda^{*}}$.
\QED

We proceed by showing that $\CP_{n}(\lambda)$ and $\CP_{n+1}(\lambda)$
have the same sign in $Range_{\lambda^{*}}$. See Fig. \ref{fig:eigenval} for a schematic description of the system.
\begin{claim}
\label{cl:same_sign}
$\Sign (\CP_{n}(\lambda)) =\Sign (\CP_{n+1}(\lambda))$
for every $\lambda \in Range_{\lambda^{*}}$.
\end{claim}
\Proof
Fix $b \in \{0,1\}$. By Claim \ref{cl:no_root_in_range}, $\CP_{n+b}$ has no roots in $Range_{\lambda^{*}}$, so $\Sign (\CP_{n+b}(\lambda_1)) =\Sign(\CP_{n+b}(\lambda_2))$ for every
$\lambda_1, \lambda_2 \in Range_{\lambda^{*}}$. Also note that by Thm. \ref{thm:pf_full}, $\Sign(\CP_{n+b}(\lambda_1) ) =\Sign(\CP_{n+b}(\lambda_2))$,
for every $\lambda_1, \lambda_2 > \lambda_{n+b}^{*}$.
We now make two crucial observations. First, as $\CP_{n}(\lambda)$ and
$\CP_{n+1}(\lambda)$ correspond to a characteristic polynomial of an
$n \times n$ matrix, they have the same leading coefficient (any characteristic polynomial is monic, i.e., with leading coefficient 1 and degree $n$) and therefore
$\Sign(\CP_{n}(\lambda))=\Sign(\CP_{n+1}(\lambda))$ for
$\lambda > \lambda_{n+1}^{*}$ (recall that we assume that
$\lambda_{n+1}^{*}\geq\lambda_{n}^{*}$).
Second, due to the \PFT, the maximal roots of
$\CP_{n}(\lambda)$  and $\CP_{n+1}(\lambda)$ are of multiplicity one and therefore
the polynomial $\CP_{n}(\lambda)$ (resp., $\CP_{n+1}(\lambda)$) necessarily changes its sign when $\lambda$ passes through its maximal real positive root $\lambda_{n}^{*}$ (respectively, $\lambda_{n+1}^{*}$).
Using these two observations, we now prove the claim via contradiction.
Assume, toward contradiction, that
$\Sign (\CP_{n}(\lambda)) \neq \Sign (\CP_{n+1}(\lambda))$ for
some $\lambda \in Range_{\lambda^{*}}$. Then
$\Sign (\CP_{n}(\lambda_1)) \neq \Sign (\CP_{n}(\lambda_2))$
for $\lambda_1 > \lambda_{n}^{*}$ and $\lambda_2 \in Range_{\lambda^{*}}$
also $\Sign (\CP_{n+1}(\lambda_1) )\neq \Sign (\CP_{n+1}(\lambda_2))$
for $\lambda_1 > \lambda_{n+1}^{*}$ and $\lambda_2 \in Range_{\lambda^{*}}$.
(This holds since when encountering a root of multiplicity one,
the sign necessarily flips). In particular, this implies that
$\Sign (\CP_{n}(\lambda))\neq \Sign (\CP_{n+1}(\lambda))$
for every $\lambda \geq \lambda_{n+1}^{*}$, in contradiction to the fact that
$\Sign(\CP_{n}(\lambda))=\Sign(\CP_{n+1}(\lambda))$
for every $\lambda>\lambda_{n+1}^{*}$. The claim follows.
\QED
We now complete the proof of Lemma \ref{thm:0_solution}.
\Proof
By Eqs. (\ref{eq:cpnz}) and
(\ref{eq:lambda_range}),
$\CP_{n-1}(\lambda) \neq 0$ for every $\lambda \in Range_{\lambda^{*}}$.
We can safely apply Claim \ref{cl:same_sign} to Lemma \ref{lem:P_n_n+1} and
and get that $\Sign (X^{*}(n)) \neq \Sign (X^{*}(n+1))$.
Since $X^{*}(n),X^{*}(n+1)$ and $g(n,n),g(n,n+1)$ are nonnegative, it follows that
$X^{*}(n)=0$ and $X^{*}(n+1)=0$. In contradiction to Lemma \ref{cl:entity_one_nonzero}.
We conclude that $\beta^{*}=\beta^{*}_{n}$.
\QED

We complete the geometric characterization of the generalized \PFT~by noting the following.
\begin{lemma}
\label{lem:no_weak_vertex}
Every vertex $\overline{X}\in V(\Polytope(\beta^{*}))$ is a $\ZeroStar$ solution.
\end{lemma}
\Proof
By Lemma \ref{cl:weak_zero_polytope}, it is sufficient to show that there
exists no $\overline{X}\in V(\Polytope(\beta^{*}))$ that is weak, namely,
which is a $\WeakZeroStar$ solution but not a $\ZeroStar$ solution.
Assume, towards contradiction, that $\overline{X}\in V(\Polytope(\beta^{*}))$
and that both $X(n) >0$ and $X(n+1)>0$. From now on, we replace
$\overline{X} \in \R^{m}$ by its truncated sub-vector in
$\R^{n+1}$, i.e., we discard the $m-n-1$ zero entries in $\overline{X}$.

Let $\System_{n-1},\System_{n}$ and $\System_{n+1}$ be defined as in Eq. (\ref{eqn:ssystem_def}). Recalling the notation of Sec. \ref{sec:per} where for matrix $A$, we denote $A_{-(i,j)}$ by the matrix that results from $A$ by removing the $i$-th row and the $j$-th column, define
$$a_i=(-1)^{n-i}  \cdot \frac{\det \left(W^{2}(\System_{n},\beta^{*})_{-(n,i)}\right)}{\det\left( W^{2}(\System_{n-1},\beta^{*})\right)}$$
and
$$b_i=(-1)^{n-i} \cdot \frac{\det \left(W^{2}(\System_{n+1},\beta^{*})_{-(n,i)} \right)}{\det\left( W^{2}(\System_{n-1},\beta^{*})\right)}$$
for $i \in \{1, \ldots, n\}$.
By Eq. (\ref{eq:CP}), Claim \ref{cl:cramer_non_square}(a)
%Eq. (\ref{eq:non_square_cramer_i})
and the proof of Lemma \ref{lem:P_n_n+1},
every optimal solution, and in particular every $\overline{X}\in V(\Polytope(\beta^{*}))$, satisfies
\begin{equation}
\label{eq:i_entries_pf}
X(i)=a_{i} \cdot X(n)+ b_{i} \cdot X(n+1)
\end{equation}
for $i \in \{1,\ldots,n-1\}$. This implies that our weak solution $\overline{X}$ is given by
$$
\overline{X}~=~ X(n) \cdot [a_1, \ldots, a_{n-1},1,0]^{T}+X(n+1)
\cdot [b_1, \ldots, b_{n-1},0,1]^{T}.
$$
Let
$$
c_n ~=~ X(n) \cdot \left(1+\sum_{i=1}^{n-1}a_{i} \right)$$ and $$c_{n+1}=X(n+1) \cdot \left(1+\sum_{i=1}^{n-1}b_{i} \right),
$$
where the feasibility of $\overline{X}$ implies $c_{n}+c_{n+1}=1$.
Next, consider Lemma \ref{lem:P_n_n+1}. Since $\overline{X}$ is optimal,
with both $X(n) >0$ and $X(n+1)>0$, it follows that
$\det \left(W^{2}(\System_{n},\beta^{*})\right)=\det \left( W^{2}(\System_{n+1},\beta^{*})\right)=0$.
This means that when constructing an optimal solution $\overline{Y}$,
one has complete freedom to select any $Y(n),Y(n+1) \geq 0$ and the
rest of the coordinates are determined by Eq. (\ref{eq:i_entries_pf}).
In particular, setting $Y(n)=X(n)/c_{n}$ and
$Y(n+1)=X(n+1)/c_{n+1}$ yields the following two optimal solutions:
$\overline{Y}_1=X(n)/c_{n}\cdot [a_1, \ldots, a_{n-1},1,0]^{T}$ and
$\overline{Y}_2=X(n+1)/c_{n+1}\cdot[b_1, \ldots, b_{n-1},0,1]^{T}$.
Note that $\overline{X}$ can be described as a convex combination of
$\overline{Y}_1$ and $\overline{Y}_2$, i.e.,
$\overline{X}=c_{n} \cdot \overline{Y}_1+ c_{n+1}\cdot \overline{Y}_2$
(recall that $c_{n}+c_{n+1}=1$). This is in contradiction to the fact
that $\overline{X}$ is a vertex of a polytope. The lemma follows.
\QED
\begin{lemma}
\label{lem:zero_star}
There exists a selection
$\FilterMatrix^{*} \in \FilterMatrixFamily$
such that $\PFEigenValue(\System(\FilterMatrix^{*}))=1/\beta^{*}$.
\end{lemma}
\Proof
Recall that our $\ZeroStar$ solution,
$\overline{X}^{*}$, is a solution for the weak subsystem $\WeakSystem$, and
therefore $\overline{X}^{*} \in \R^{n+1}$. In addition, $|NZ(\overline{X}^{*})|=n$
and due to Lemma \ref{cl:entity_one_nonzero}, $|NZ(\overline{X}^{*}) \cap \Supporters_i|=1$ for every $\Entity_{i}$, or in other words,
$\SelectionVec'=NZ(\overline{X}^{*})$ is a complete selection for $\EntitySet$
such that $|\SelectionVec'|=n$. Taking $\FilterMatrix^{*}=\FilterMatrix(\SelectionVec')$ yields the desired claim. The lemma follows.
\QED

%\textbf{MP: I add this here, maybe it should be placed elsewhere.}\\
Note that Eq. (\ref{eq:i_entries_pf}) illustrates the additional degrees
of freedom at the optimum point of Program (\ref{LP:Ext_Perron}).
Specifically, to obtain an optimum solution for $\beta^{*}$,
one has the freedom to set $X_{n}\geq 0$ and $X_{n+1} \geq 0$
(as long as at least one of them is positive) and the rest of the coordinates are determined accordingly.

We are now ready to complete the proof of Thm. \ref{thm:pf_ext}.
%The complete Proof for Lemma \ref{lem:zero_star}  and Thm. \ref{thm:pf_ext} are given in Appendix \ref{sec:geoproofs}.
%
\Proof [Thm. \ref{thm:pf_ext}]
Let $\FilterMatrix^{*}$ be the selection such that
$\PFEigenValue(\System)=\PFEigenValue(\System(\FilterMatrix^{*}))$.
Note that by the irreducibility of $\System$, the square system
$\System(\FilterMatrix^{*})$ is irreducible as well and therefore the \PFT~
for irreducible matrices can be applied.
In particular, by Thm. \ref{thm:pf_full}, it follows that
$\PFEigenValue(\System(\FilterMatrix^{*})) \in \R_{>0}$ and that
$\PFEigenVector(\System(\FilterMatrix^{*}))>0$.
Therefore, by Eq. (\ref{eq:general_pf_root}) and (\ref{eq:general_pf_vector}),
Claims (Q1)-(Q3) of Thm. \ref{thm:pf_ext} follow.

We now turn to claim (Q4) of the theorem. Note that for a symmetric system, in which
$g(i,j_1)=g(i,j_2)$ for every
$\Affectors_{j_1},\Affectors_{j_2} \in \Supporters_k$ and every $k,i \in [1,n]$,
the system is invariant to the selection matrix and therefore
$\PFEigenValue(\System(\FilterMatrix_1)) =
\PFEigenValue(\System(\FilterMatrix_2))$
for every $\FilterMatrix_1,\FilterMatrix_2 \in \FilterMatrixFamily$.

Finally, it remains to consider claim (Q5) of the theorem. Note that the optimization problem
specified by Program (\ref{LP:Ext_Perron}) is an alternative formulation
to the generalized Collatz-Wielandt formula given in (Q5).
We now show that $\PFEigenValue(\System)$ (respectively,
$\PFEigenVector(\System)$) is the optimum value (resp., point)
of Program (\ref{LP:Ext_Perron}).
By Lemma \ref{lem:zero_star}, there exists an optimal point
$\overline{X}^{*}$ for Program (\ref{LP:Ext_Perron}) which is a $\ZeroStar$
solution. Note that a $\ZeroStar$  solution corresponds to a unique hidden
square system, given by $\System^{*}=\System(NZ(\overline{X}^{*}))$
($\System^{*}$ is square since $|NZ(\overline{X}^{*})|=n$).
Therefore, by Thm. \ref{thm:pf} and Lemma \ref{lem:zero_star},
we get that
\begin{equation}
\label{eq:put_all_val}
\PFEigenValue(\System^{*}) ~=~ 1/\beta^{*}(\System^{*}) ~=~ 1/\beta^{*}(\System).
\end{equation}
Next, by Observation \ref{obs:filter_to_square}(b), we have that
$\PFEigenValue(\System(\FilterMatrix)) \geq
\PFEigenValue(\System)$. It therefore follows that
\begin{equation}
\label{eq:put_all}
\PFEigenValue(\System^{*}) ~=~
\min_{\FilterMatrix \in \FilterMatrixFamily} \PFEigenValue(\System(\FilterMatrix)).
\end{equation}
Combining Eq. (\ref{eq:put_all_val}), (\ref{eq:put_all}) and
(\ref{eq:general_pf_root}),
we get that the PF eigenvalue of the system $\System$
satisfies $\PFEigenValue(\System)=1/\beta^{*}(\System)$ as required.
Finally, note that
by Thm. \ref{thm:pf}, $\PFEigenVector(\System^{*})$ is the optimal
point for Program (\ref{LP:Ext_Perron}) with the square system $\System^{*}$.
By Eq. (\ref{eq:general_pf_vector}), $\PFEigenVector(\System)$ is an extension
of $\PFEigenVector(\System^{*})$ with zeros (i.e., a $\ZeroStar$ solution).
It can easily be checked that $\PFEigenVector(\System)$ is a feasible solution
for the original system $\System$ with
$\beta=\beta^{*}(\System^{*})=\beta^{*}(\System)$, hence it is optimal.
Note that by Lemma \ref{lem:strict_equality}, it indeed follows that $\RepressorsMatrix \cdot \PFEigenVector(\System) =
1/\beta^{*}(\System) \cdot \SupportersMatrix \cdot \PFEigenVector(\System)$,
for every optimal solution $\overline{X}^{*}$.
Theorem \ref{thm:pf_ext} follows.
\QED

%%%%%%%%%%%%%%%%%%%%%%%%%%%%%%%%%%%%%%
\section{Computing the generalized PF vector}
\label{short:sec:Algorithm}
%%%%%%%%%%%%%%%%%%%%%%%%%%%%%%%%%%%%%%%%%%%%
In this section we present a polynomial time algorithm for computing the generalized Perron eigenvector $\PFEigenVector(\System)$ of an irreducible
system $\System$.

%%%%%%%%%%%%%%%%%%%%%%%%%%%%%%%%%%%%%%%%%%%%
\paragraph{The method.}
%%%%%%%%%%%%%%%%%%%%%%%%%%%%%%%%%%%%%%%%%%%%
By Property (Q5) of Thm. \ref{thm:pf_ext}, computing $\PFEigenVector(\System)$
is equivalent to finding a $\ZeroStar$ solution for
Program (\ref{LP:Ext_Perron}) with $\beta=\beta^{*}(\System)$.
For ease of analysis, we assume throughout  that the gains are integral,
i.e., $g(i,j) \in  \mathbb{Z}^{+}$, for every $i \in \{1, \ldots, n\}$ and $j \in \{1, \ldots, m\}$.
If this does not hold, then the gains can be rounded or scaled to achieve this.
Let
\begin{equation}
\label{eq:gmax}
\MaxGain(\System) ~=~ \max_{i \in \{1,\ldots, n\}, j \in \{1, \ldots, m\}} \left\{|g(i,j)| \right\},
\end{equation}
and define $\LPRunTime$ as the running time of an LP solver such as the
interior point algorithm \cite{Boyd-Conv-Opt-Book} for Program (\ref{LP:Ext_Perron_convex}).
Recall that we look for an exact optimal solution for a non-convex optimization
problem (see Program (\ref{LP:Ext_Perron})). Using the convex relaxation of
Program (\ref{LP:Ext_Perron_convex}), a binary search can be applied
for finding an approximate solution up to a predefined accuracy.
The main challenge is then to find (a) an optimal solution
(and not an approximate one), and (b) among all the optimal solutions,
to find one that is a $\ZeroStar$ solution.
Let $\FilterMatrix_1, \FilterMatrix_2 \in \FilterMatrixFamily$ be two
selection matrices for $\System$. By Thm. \ref{thm:pf_ext}, there exists
a selection matrix $\FilterMatrix^{*}$ such that
$\PFEigenValue(\System)=\PFEigenValue(\System(\FilterMatrix^{*}))$ and
$\PFEigenVector(\System)$ is a $\ZeroStar$ solution corresponding to
$\PFEigenVector(\System(\FilterMatrix^{*}))$
(in addition $\beta^{*}=1/\PFEigenValue(\System(\FilterMatrix^{*}))$).
Our goal then is to find a selection matrix $\FilterMatrix^{*} \in \FilterMatrixFamily$ where
$|\FilterMatrixFamily|$ might be exponentially large.
%In Sec. \ref{sec:alg} we prove the following.
\begin{theorem}
\label{thm:algorithm}
Let $\System$ be an irreducible system. Then $\PFEigenVector(\System)$ can be
computed in time
$O(n^{3} \cdot \LPRunTime \cdot
\left(\log \left(n \cdot \MaxGain \right) +n \right))$.
\end{theorem}
Let
\begin{equation}
\label{eq:delta_beta}
\Delta_{\beta} ~=~ (n\MaxGain)^{-8n^3}.
\end{equation}
The key observation in this context is the following ``minimum gap" observation.
\begin{lemma}
\label{lem:apart_in_range}
Consider a selection matrix $\FilterMatrix \in \FilterMatrixFamily$. If
$\beta^*(\System)-1/\PFEigenValue(\System(\FilterMatrix)) \leq \Delta_{\beta}$,
then $\beta^{*}(\System)=1/\PFEigenValue(\System(\FilterMatrix))$.
\end{lemma}
By performing a polynomial number of steps of binary search for the optimal
$\beta^{*}(\System)$, one can converge to a value $\beta^{-}$  that is at most
$\Delta_{\beta}$ far from $\beta^{*}(\System)$, i.e.,
$\beta^{*}(\System)-\beta^{-}<\Delta_\beta$.
Let $Range_{\beta^{*}}=[\beta^{-}, \beta^{*}]$.
Then by Lemma \ref{lem:zero_star}, we are guaranteed that
$\PFEigenValue(\System(\FilterMatrix))=1/\beta^{*}$
for any selection matrix $\FilterMatrix \in \FilterMatrixFamily$ such that
$1/\PFEigenValue(\System(\FilterMatrix)) \in Range_{\beta^{*}}$
(there could be many such matrices $\FilterMatrix$, but in this case,
they all correspond to systems with PF value $1/\beta^{*}$).
To prove Lemma \ref{lem:apart_in_range}, we first establish a lower
bound on the difference between \emph{any} two different PF eigenvalues of any
two irreducible square systems, i.e., we show that the PF roots
$\PFEigenValue(\System^{s}_1)$ and $\PFEigenValue(\System^{s}_2)$
of any two irreducible square systems
$\System^{s}_1, \System^{s}_2 \in \SquareSystemFamily$
cannot be too close if they are different.
Recall that for an irreducible square system $\System^{s}$,
$Z(\System^{s})=(\SupportersMatrix)^{-1} \cdot \RepressorsMatrix$,
where $\SupportersMatrix$ can be considered to be diagonal
with a strictly positive diagonal.
%
%For ease of analysis,
We begin the analysis by scaling the entries of $Z(\System^{s})$ to obtain
an integer-valued matrix $Z^{\mathrm{int}}$. The scaling is needed in order to
employ a well-known bound due to Bugeaud and Mignotte \cite{OnDistBetwRoots}
on the minimal distance between the roots
of integer polynomials (Lemma \ref{lemma:distance_of_roots}).
The guaranteed distance on
$\PFEigenValue(\System^{s}_1)$ and $\PFEigenValue(\System^{s}_2)$ is later
translated into a minimal bound on distance for their reciprocals
$1/\PFEigenValue(\System^{s}_1)$ and $1/\PFEigenValue(\System^{s}_2)$,
which correspond to $\beta$ values of Program (\ref{LP:Ext_Perron}),
i.e., optimal $\beta$ values of two different irreducible square systems
for Program (\ref{LP:Ext_Perron}).
Specifically, we show that for any given
sufficiently small range of $\beta$ values, $Range_\beta=[\beta_1, \beta_2]$
such that $|\beta_1-\beta_2| \leq \Delta_{\beta}$, there cannot be two
selection matrices
$\FilterMatrix_1,\FilterMatrix_2  \in \FilterMatrixFamily$ such that
$\PFEigenValue(\System(\FilterMatrix_1)) \neq
\PFEigenValue(\System(\FilterMatrix_2))$
and yet both
$1/\PFEigenValue(\System(\FilterMatrix_1)),
1/\PFEigenValue(\System(\FilterMatrix_2)) \in Range_\beta$.
\par The \emph{na\"{\i}ve height} of an integer polynomial $P$, denoted $H(P)$, is the maximum of the absolute
values of its coefficients.
\begin{lemma}[Bugeaud and Mignotte \cite{OnDistBetwRoots}]
\label{lemma:distance_of_roots}
Let $P(X)$ and $Q(X)$ be nonconstant integer polynomials of degree $n$ and $m$,
respectively. Denote by $r_P$ and $r_Q$ a zero of $P(X)$ and $Q(X)$, respectively.
Assuming that $P(r_Q)\ne 0$, we have
\begin{eqnarray*}
|\ r_P - r_Q |\ \ge
2^{1-n}(n+1)^{\frac{1}{2}-m}(m+1)^{-\frac{n}{2}}H(P)^{-m}H(Q)^{-n}.
\end{eqnarray*}
\end{lemma}
We first show the following.
\begin{lemma}
\label{lem:new_height}
$|\PFEigenValue(\System^{s}_1)-\PFEigenValue(\System^{s}_2)| \ge
(n\MaxGain)^{-6n^3}$ for every
$\System^{s}_1, \System^{s}_2\in \SquareSystemFamily$.
\end{lemma}
\Proof
Recall that for an irreducible square system $\System^{s}$,
$Z(\System^{s})=(\SupportersMatrix)^{-1} \cdot \RepressorsMatrix$,
where
$\SupportersMatrix$ can be considered to be diagonal with strictly positive
diagonal. Therefore, $Z(\System^{s})_{i,j}=|g(i,j)|/g(i,i)$ where $g(i,i)$
corresponds to the gain of the unique supporter of $\Entity_i$.

For ease of notation, let $Z_1=Z(\System^{s}_1)$, $Z_2=Z(\System^{s}_2)$,
$r_1=\PFEigenValue(\System^{s}_1)$ and $r_2=\PFEigenValue(\System^{s}_2)$.
Let $i_1$ (resp., $i_2$) be the index of the unique supporter of entity
$\Entity_{i}$ in the square system $\System^{s}_1$ (resp., $\System^{s}_2$).

To employ Lemma \ref{lemma:distance_of_roots}, we first scale $Z_1$ and $Z_2$
to obtain two integer-valued matrices $Z^{\mathrm{int}}_1$ and $Z^{\mathrm{int}}_2$.
The new matrix $Z_b^{\mathrm{int}}$, for $b \in \{1,2\}$, is constructed by multiplying each entry of
$Z_b$ by the common denominator of its entries,
i.e., $Z_{b}^{\mathrm{int}}(i,j)=Z_b(i,j)\cdot\prod_{i} \left(|g(i,i_1)| \cdot |g(i,i_2)| \right)$.
Thus all entries of $Z_b^{\mathrm{int}}$ are integers and bounded by $\MaxGain^{2n}$
(since $|g(i,j)|\le \MaxGain$).
Let $P_1(x)=\CP(Z^{\mathrm{int}}_1,x)$ and $P_2(x)=\CP(Z^{\mathrm{int}}_2,x)$ be the
characteristic polynomials of the matrices $Z^{\mathrm{int}}_1$ and
$Z^{\mathrm{int}}_2$ respectively, see Eq. (\ref{eq:CP}).
Note that $P_1(x)$ and $P_2(x)$ are integer polynomials of degree $n$,
and $H(P_1),H(P_2) \leq \MaxGain^{2n^2}$ (since $|\det(Z)|\le (\MaxGain^{2n})^n$).
Let $r_1^{\mathrm{int}}$ and $r_2^{\mathrm{int}}$ correspond to the PF eigenvalues
of $Z_{1}^{\mathrm{int}}$ and $Z_{2}^{\mathrm{int}}$ respectively.
Lemma \ref{lemma:distance_of_roots} yields
%$$
\begin{eqnarray*}
|r_1^{\mathrm{int}}-r_2^{\mathrm{int}}| &\geq&
2^{1-n}(n+1)^{\frac{1}{2}-n}(n+1)^{-\tfrac{n}{2}}
(\MaxGain^{2n^2})^{-n}(\MaxGain^{2n^2})^{-n}
= 2^{1-n}(n+1)^{\frac{1-3n}{2}}\MaxGain^{-4n^3}.
\end{eqnarray*}
%$$
Finally, by definition of $Z^{\mathrm{int}}_1$ and $Z^{\mathrm{int}}_2$,
%$$
\begin{align*}
|r_1^{\mathrm{int}}-r_2^{\mathrm{int}}| ~=~ |r_1-r_2|\prod_{i} \left(|g(i,i_1)| \cdot |g(i,i_2)|\right) ,
\end{align*}
%$$
and thus
%$$
\begin{eqnarray*}
|r_1-r_2| &\geq&
\frac{2^{1-n}(n+1)^{\frac{1-3n}{2}}\MaxGain^{-4n^3}}{\prod_{i}\left(|g(i,i_1)| \cdot |g(i,i_2)| \right)}
~\geq~
\frac{2^{1-n}(n+1)^{\frac{1-3n}{2}}\MaxGain^{-4n^3}}{\MaxGain^{2n}}
\geq (n\MaxGain)^{-6n^3}.
\end{eqnarray*}
\QED
We now turn to translate the distance between $r_1$ and $r_2$ into a distance
between $1/r_1$ and $1/r_2$ (corresponding to the optimal $\beta$ values of
Program (\ref{LP:Ext_Perron}) with $\System^{s}_1$ and $\System^{s}_2$,
respectively).
The next auxiliary claim gives a bound for $\lambda \in \EigenValue(A)$
as a function of $\MaxGain$.

\begin{lemma}
\label{claim:lambd_up_bound}
Let $\lambda$ be an eigenvalue of an $n\times n$ matrix $Z$ such that
$|Z(i,j)|\le \MaxGain$. Then $|\lambda|\le n \MaxGain$.
\end{lemma}
\Proof
Let $\overline{X}$ be the eigenvector of $Z$ and assume that
$||\overline{X}||_{2}=1$. Since
$\overline{X}^{T} \cdot Z \cdot \overline{X} =
\lambda \overline{X}^{T} \cdot \overline{X}=\lambda$, we have:
\begin{eqnarray*}
|\lambda | &=& |\overline{X}^{T}Z \overline{X}|
~=~ |\sum_i\sum_j X(i) Z(i,j) X(j)|
~\le~ \MaxGain\cdot |\sum_i\sum_j X(i) \cdot X(j)|\\
\\ &=& \MaxGain\cdot |\sum_i X(i)|\cdot |\sum_j X(j)|
~=~ \MaxGain\cdot\|\overline{X} \|_1^2
~\le~ \MaxGain\cdot(\sqrt{n}\|\overline{X}\|_2)^2 ~=~ n\MaxGain~.
\end{eqnarray*}
\QED
We now turn to prove Lemma \ref{lem:apart_in_range}.
\Proof [of Lemma \ref{lem:apart_in_range}]
\\By Lemma \ref{lem:new_height} and \ref{claim:lambd_up_bound},
$$\left|\frac{1}{r_2}-\frac{1}{r_1}\right| =
\left|\frac{r_1-r_2}{r_1 r_2}\right| \ge
\frac{\left|r_1-r_2\right|}{(n\MaxGain)^2} \ge
(n\MaxGain)^{-8n^3}.$$
%\end{align*}
So far, we proved that if
$\PFEigenValue(\System(\FilterMatrix_1)) \not=
\PFEigenValue(\System(\FilterMatrix_2))$,
then
$|1/\PFEigenValue(\System(\FilterMatrix_1))
- 1/\PFEigenValue(\System(\FilterMatrix_2))| \geq \Delta_{\beta}$,
for every $\FilterMatrix_1,\FilterMatrix_2 \in \FilterMatrixFamily$.
By Thm. \ref{thm:pf_ext}, there exists a selection
$\FilterMatrix^{*} \in \FilterMatrixFamily$ such that
$\PFEigenValue(\System(\FilterMatrix^{*}))=1/\beta^{*}(\System)$.
Assume, toward contradiction, that there exists some
$\FilterMatrix^{'} \in \FilterMatrixFamily$ such that
$\PFEigenValue(\System(\FilterMatrix^{'}))\neq 1/\beta^{*}(\System)$ but
$|\beta^{*}(\System)-1/\PFEigenValue(\System(\FilterMatrix^{'}))| \leq
\Delta_{\beta}$. Let $r_1=\PFEigenValue(\System(\FilterMatrix^{*}))$ and
$r_2=\PFEigenValue(\System(\FilterMatrix^{'}))$. In this case, we get that
$|1/r_1-1/r_2|\leq \Delta_{\beta}$, contradiction.
Lemma \ref{lem:apart_in_range} follows.
\QED

%%%%%%%%%%%%%%%%%%%%%%%%%%%%%%%%%%%%%%%%%%%%
\paragraph{Algorithm description.}
%%%%%%%%%%%%%%%%%%%%%%%%%%%%%%%%%%%%%%%%%%%%
We now describe Algorithm $\AlgoName$ for $\PFEigenVector(\System)$ computation.
%\textbf{MP: I do not think this theorem is is needed}\\
%First, we present a well known result of Leonid Khachiyan that proposed the {\em ellipsoid method} to solve linear programming problems in polynomial time.\\
%\begin{theorem}[Khachiyan:79 \cite{Khachiyan:79}]\label{thm:khachiyan}
%The following problems are solvable in polynomial time:
%\begin{itemize}
%\item Given a matrix $A \in \mathbb{Q}^{m\times n}$ and a vector $\overline{B} \in \mathbb{Q}^m$, decide whether $A\overline{X} \le \overline{B}$ has a solution $\overline{X} \in \mathbb{Q}^n$, and if so, find one.
%\item (Linear programming problem) Given a matrix $A \in \mathbb{Q}^{m\times n}$ and vectors $\overline{B} \in \mathbb{Q}^m, \overline{C} \in \mathbb{Q}^n$, decide whether
%$\max\{\overline{C}^T \overline{X} \mid A\overline{X} \le \overline{B},\overline{X} \in \mathbb{Q}^n\}$ is infeasible, finite, or unbounded. If it is finite, find an optimal solution. If it is unbounded, find a feasible solution $\overline{X}_0$, and find a vector $\overline{D} \in \mathbb{Q}^n$ with $A\overline{D} \le 0$ and $\overline{C}^T \overline{D} > 0$.
%\end{itemize}
%\end{theorem}
%
%\textbf{MP: I think we can start from here}.
Consider some partial selection $\SelectionVec'\subseteq\Affectors$ for
$V'\subseteq \EntitySet$.
For ease of notation, let
$\System(\SelectionVec')=\langle \RepressorsMatrix(\SelectionVec'), \SupportersMatrix(\SelectionVec') \rangle$,
where
$\RepressorsMatrix(\SelectionVec')=\RepressorsMatrix \cdot \FilterMatrix(\SelectionVec')$ and
$\SupportersMatrix(\SelectionVec')=\SupportersMatrix \cdot \FilterMatrix(\SelectionVec')$.
Consider the Program
\begin{align}
\label{eq:alg}
  \text{maximize} & ~\beta \text{~subject to:~}
  \\
  &
\displaystyle \RepressorsMatrix(\SelectionVec') \nonumber \cdot\overline{X} \leq
                 1/\beta \cdot \SupportersMatrix(\SelectionVec') \cdot\overline{X} , &
\\
  & \displaystyle \overline{X} \geq \overline{0} ,& \nonumber
\\
  & \displaystyle ||\overline{X}||_{1}~=~1~. & \nonumber
\end{align}

Note that if $\SelectionVec'=\emptyset$, then
Program (\ref{eq:alg}) is equivalent to Program
(\ref{LP:Ext_Perron}), i.e., $\System(\SelectionVec')=\System$.
Define
\begin{eqnarray*}
\fff(\beta,\System(\SelectionVec'))=
\left\{\begin{array}{ll}
1,  &  \mbox{if there exists an~} \overline{X} \mbox{ such that }
  ||\overline{X}||_{1}~=~1, ~\overline{X}\geq \overline{0}, \mbox{ and }
\\ &
\displaystyle \RepressorsMatrix(\SelectionVec')\cdot\overline{X} \leq
        1/\beta \cdot \SupportersMatrix(\SelectionVec')\cdot\overline{X},
\\
0,  &  \hbox{otherwise}.
\end{array}\right.
\end{eqnarray*}

Note that $\fff(\beta,\System(\SelectionVec'))=1$ iff $\System(\SelectionVec')$ is feasible for
$\beta$ and that $\fff$ can be computed in polynomial time using the interior point method.

Algorithm $\AlgoName$ is composed of two main phases.
In the first phase it finds, using binary search, an estimate $\beta^-$ such that
$\beta^{*}(\System)-\beta^- \leq \Delta_{\beta}$.
In the second phase, it finds a hidden square system
$\System(\FilterMatrix^{*})$, $\FilterMatrix^{*} \in \FilterMatrixFamily$,
corresponding to a complete selection vector $\SelectionVec_n$ of size $n$ for $\EntitySet$.
%We note that the system
%$\System(\SelectionVec_n)=\langle %\SupportersMatrix(\SelectionVec_n),\RepressorsMatrix(\SelectionVec_n) \rangle$
%is equivalent to $\System(\FilterMatrix(\SelectionVec_n)) \in %\SquareSystemFamily$.
%This is because
%$\SupportersMatrix(\SelectionVec_n),\RepressorsMatrix(\SelectionVec_n) \in \R^{n \times n}$,
%so by removing the $m-n$ zero columns of
%$\SupportersMatrix(\SelectionVec_n)$ and $\RepressorsMatrix(\SelectionVec_n)$, we get a square system
%$\System(\FilterMatrix(\SelectionVec_n)) \in \SquareSystemFamily$.
By Lemma \ref{lem:apart_in_range}, it follows that
$\PFEigenValue(\System(\FilterMatrix^{*}))=1/\beta^{*}(\System)$.

We now describe the construction of $\SelectionVec_n$ in more detail.
The second phase consists of $n$ iterations. Iteration $t$ obtains a partial
selection $\SelectionVec_{t}$ for $\Entity_1, \ldots, \Entity_t$ such  that
$\fff(\beta^-,\System(\SelectionVec_t))=1$. The final step achieves the desired $\SelectionVec_n$,
where $\System(\SelectionVec_n) \in \SquareSystemFamily$ and
$\fff(\beta^-,\System(\SelectionVec_n))=1$
(therefore also $\fff(\beta^-,\System(\FilterMatrix(\SelectionVec_n)))=1$).
Initially, $\SelectionVec_0$ is empty.
The $t$'th iteration sets $\SelectionVec_t=\SelectionVec_{t-1} \cup \{\Affectors_j\}$ for some supporter
$\Affectors_j \in\Supporters_{t}$
such that $\fff(\beta^-,\System(\SelectionVec_{t-1} \cup \{\Affectors_j\}))=1$.
We later show (in proof of Thm. \ref{thm:algorithm}) that such a supporter
$\Affectors_j$ exists.

Finally, we use $\PFEigenVector(\System(\SelectionVec_n))$ to construct
the Perron vector $\PFEigenVector(\System)$. This vector contains zeros
for the $m-n$ non-selected affectors, and the values of the $n$ selected affectors
are as in  $\PFEigenVector(\System(\SelectionVec_n))$.

The pseudocode is presented formally next.
%in Figure \ref{figure:pseudo code Algorithm ?name?}.

%%%%%%%%%%%%%%%%%%%%%%
%\begin{figure}[htb]
\begin{figure*}[htb]
\begin{center}
\framebox{\parbox{6in}{
\noindent{\bf Algorithm $\AlgoName$}
\begin{enumerate}
\dnsitem[]
/* Binary search phase: finding $\beta^-$ such that
    $\beta^*-\beta^-<\Delta_\beta$ */
\dnsitem
$\beta\leftarrow 1$;
\dnsitem
While  $\fff(\beta,\System)=1$ do: \\
\hskip 30pt
$\beta\leftarrow 2\beta$;
\dnsitem
If $\beta>1$, then $\beta^{-}\leftarrow\beta/2$, else $\beta^-\leftarrow0$;
\dnsitem
$\beta^+\leftarrow \beta$;
\dnsitem
While  $\beta^{+}-\beta^-\geq\Delta_\beta$ do:~~~~
      /* from now on  $\beta^-\leq \beta^*<\beta^+$ */
\label{alg:while loop}
    \begin{enumerate}
    \ddnsitem $\beta\leftarrow \left(\beta^-+\beta^+ \right)/{2}$;
    \ddnsitem If $\fff(\beta,\System)=1$, then $\beta^-\leftarrow\beta$,
           else $\beta^+\leftarrow\beta$;
    \end{enumerate}
\dnsitem[]
~~~~~~~~~~~~~/* Affector elimination phase: Finding a $\ZeroStar$ solution */
\dnsitem
$\SelectionVec_{0}\leftarrow\emptyset$;
\dnsitem
For $t=1$ to $n$ do:
\label{alg:for loop}
  \begin{enumerate}
   \ddnsitem Select some supporter $\Affectors_j\in \Supporters_{t}$
   such that $\fff(\beta^-,\System(\SelectionVec_{t}\cup \{\Affectors_j\}))=1$;
   \ddnsitem Set $\SelectionVec_{t+1} \gets \SelectionVec_{t}\cup \{\Affectors_j\}$;
  \end{enumerate}
\dnsitem[]
/* $|\SelectionVec_n|=n$ and  $\System(\SelectionVec_n)\in \SquareSystemFamily$ ~*/
\dnsitem
Set Perron value $\beta^*=1/r(\System(\SelectionVec_n))$;
\dnsitem
Set $\overline{X}^{*} \gets \PFEigenVector(\System(\SelectionVec_n))$;
\dnsitem
Set $X^{**}(\Affectors_k) \gets
\begin{cases}
X^{*}(\Affectors_k), & \Affectors_k \in \SelectionVec_n,\\
0, & \text{otherwise;}
\end{cases}
$
\dnsitem
Let $\PFEigenVector(\System) \gets \overline{X}^{**}$;
\end{enumerate}
%\caption{\label{figure:pseudo code Algorithm ?name?}
%Pseudocode of Algorithm $\AlgoName$.}
}}
\end{center}
\end{figure*}
%%%%%%%%%%%%%%%%%%%%%%%%

%\begin{lemma}
%\label{correctness}
%Let $\System$ be an irreducible system. Then Algorithm $\AlgoName$ finds
%$\PFEigenVector(\System)$ in time
%$(\log \left(\MaxGain \Delta_{\beta}\right) +n) \cdot \LPRunTime)$.
%\end{lemma}
To establish Theorem \ref{thm:algorithm}, we prove the correctness of
Algorithm $\AlgoName$ and bound its runtime.
We begin with two auxiliary claims.

\begin{claim}
\label{cl:max_beta}
$\beta^{*}(\System) \leq \MaxGain$.
\end{claim}
\Proof
Let $\overline{X}^{*}=\PFEigenVector(\System)$ and let
$\SelectionVec^{*}=NZ(\overline{X}^{*})$.
Then by claims (Q3) and (Q5) of Thm. \ref{thm:pf_ext} we have that $|\SelectionVec^{*}|=n$.
Define $\FilterMatrix^{*}= \FilterMatrix(\SelectionVec^{*})$.
Since $\SelectionVec^{*}$ is a complete selection vector
(see Claim \ref{cl:entity_one_nonzero}), we have that
$\FilterMatrix^{*}\in \FilterMatrixFamily$.
Let $\Affectors_{\IndS(i)}$ be the supporter of entity
$\Entity_i$ in $\SelectionVec^{*}$, for every $i \in \{1,\ldots, n\}$.

Let  $D=\ConstraintsGraph_{\System(\FilterMatrix^{*})}$. Since $\System$ is
irreducible, it follows by Obs. \ref{obs:reducible_graph_connected} that $D$ is strongly connected.
Let $C=(\Entity_{i_1}, \ldots, \Entity_{i_k})$ be a directed cycle in $D$, i.e.,
$(\Entity_{i_j},\Entity_{i_{j+1}}) \in E(D)$ for every $j \in \{1, \ldots, k\}$ and
$(\Entity_{i_k},\Entity_{i_{1}}) \in E(D)$.
For ease of notation, let $\Entity_{i_k}=\Entity_{i_{-1}}$.
Since $D$ is strongly connected, such a cycle $C$ exists.
By the optimality of $\overline{X}^{*}$ we have that
$$\beta^{*}(\System) \cdot \TotR(\overline{X}^{*}, \System)_{i} ~=~
\TotS(\overline{X}^{*}, \System)_{i}$$
for every $\Entity_i$. Note that by definition
%
%MP: The following equation concerns separately $i=1$. But now we define
%$\Entity_{i_k}=\Entity_{i_{-1}}$ so it is not needed.
%
%\begin{eqnarray*}
%-g(\Affectors_{i_k},\Entity_{i_1}) \cdot X^{*}(\Affectors_{i_k}) ~\leq~
%\TotR(\overline{X}^{*},
%\System)_{i_1};
%\end{eqnarray*}
%and
$|g(\Entity_{i_j},\Affectors_{\IndS(i_{j-1})})| \cdot X^{*}(\Affectors_{\IndS(i_{j-1})}) \leq
\TotR(\overline{X}^{*}, \System)_{i_j}$ for every $j \in \{1, \ldots, k\}$,
and by the graph definition, $\Affectors_{\IndS(i_{j-1})} \in \Repressors_{i_j}$ or
$g(\Entity_{i_j},\Affectors_{\IndS(i_{j-1})})<0$, for every $j \in \{1, \ldots, k\}$.
Combining this with Fact \ref{fc:feasible_tots_totr}, we get that
%Following equation is hidden for same reason as above.
%
%\begin{eqnarray*}
%-\beta^{*}(\System) \cdot g(\Affectors_{i_k},\Entity_{i_1}) \cdot
%X^{*}(\Affectors_{i_k}) ~\leq~
%g(\Affectors_{i_1},\Entity_{i_1}) \cdot X^{*}(\Affectors_{i_1});
%\end{eqnarray*}
%and
\begin{eqnarray*}
\beta^{*}(\System) %\cdot
|g(\Entity_{i_j},\Affectors_{\IndS(i_{j-1})})| %\cdot
X^{*}(\Affectors_{\IndS(i_{j-1})}) ~\leq~
g(\Affectors_{\IndS(i_j)},\Entity_{i_j}) \cdot X^{*}(\Affectors_{\IndS(i_j)})
\end{eqnarray*}
for every $j \in \{1, \ldots, k\}$, and therefore
$$\beta^{*}(\System) ~\leq~ \min_{j \in \{1, \ldots, k\}}\left\{
\frac{g(\Affectors_{\IndS(i_j)},\Entity_{i_j})}{|g(\Affectors_{\IndS(i_{j-1})},\Entity_{i_j})|}
\cdot \frac{X^{*}(\Affectors_{\IndS(i_j)})}{X^{*}(\Affectors_{\IndS(i_{j-1})})}  \right\}.
$$
%for every $j \in [2,k]$ as well as
%$$\beta^{*}(\System) \leq \min\left\{
%\frac{g(\Affectors_{i_{k}},\Entity_{i_1})}{g(\Affectors_{i_1},\Entity_{i_1})}
%\cdot \frac{X^{*}(\Affectors_{i_1})}{X^{*}(\Affectors_{i_{k}})}  \right\}.$$
It is easy to verify that
$\min_{j \in \{1, \ldots, k\}}
\left \{\frac{X^{*}(\Affectors_{\IndS(i_j)})}{X^{*}(\Affectors_{\IndS(i_{j-1})})}\right\} \leq 1$.
Therefore, by Eq. (\ref{eq:gmax})
we get that $\beta^{*}(\System) \leq \MaxGain$, as required.
\QED
\begin{lemma}
\label{cl:phase_1}
Phase 1 of Alg. $\AlgoName$ finds $\beta^{-}$ such that
$\beta^{*}(\System)-\beta^{-} \leq \Delta_{\beta}$.
\end{lemma}
\Proof
By Property (Q5) of Thm. \ref{thm:pf_ext}, $\PFEigenVector(\System)$
is an optimal solution for Program  (\ref{LP:Ext_Perron}) and
$\PFEigenValue(\System)=1/\beta^{*}(\System)$. Therefore
$\fff(\beta,\System)=1$ for every $\beta \in (0,\beta^{*}]$.
Steps 3 and 5(b) in Alg. $\AlgoName$
%(Figure \ref{figure:pseudo code Algorithm ?name?}),
yield
$\fff(\beta^-,\System)=1$. Therefore $\beta^{-} \leq \beta^{*}(\System)$.
By the stopping criterion of step 5, it ends with
$\fff(\beta^{+},\System)=0$, $\fff(\beta^{-},\System)=1$ and
$\beta^{+}-\beta^{-} \leq \Delta_{\beta}$. The first 2 conditions imply that
$\beta^{*} \in [\beta^{-},\beta^{+})$ as required. The claim follows.
\QED
Let $Range_{\beta^{*}}=[\beta^{-},\beta^{+})$.
%Note that at the beginning of phase 2 of Alg. $\AlgoName$, the computed value
%$\beta^{-}$ is at most $\Delta_{\beta}$ apart from $\beta^{*}$.
%We begin by showing the following.

\begin{lemma}
\label{cl:selection_opt}
By the end of phase 2, the selection $\SelectionVec_n$ satisfies
$\PFEigenValue(\System(\SelectionVec_n))=1/\beta^{*}(\System)$.
\end{lemma}
\Proof
Let $\SelectionVec_t$ be the partial selection obtained at step $t$,
$\System_{t}=\System(\SelectionVec_{t})$ be the corresponding system for step $t$
and $\beta_{t}=\beta^{*}(\System_{t})$ the optimal solution of Program
(\ref{LP:Ext_Perron}) for system $\System_{t}$.
We claim that $\SelectionVec_{t}$ satisfies the following properties for each
$t \in \{0, \ldots, n\}$:
\begin{description}
\item{(C1)}
$\SelectionVec_{t}$ is a partial selection vector of length $t$,
such that $\SelectionVec_t \sim \SelectionVec_{t-1}$.
\item{(C2)}
$\System(\SelectionVec_{t})$ is feasible for $\beta^{-}$.
\end{description}
The proof is by induction.
Beginning with $\SelectionVec_{0}=\emptyset$, it is easy to see that (C1) and (C2)
are satisfied (since $\System(\SelectionVec_0)=\System$). Next, assume that (C1) and (C2)
hold for $\SelectionVec_{i}$ for $i \leq t$
and consider $\SelectionVec_{t+1}$. Let $V_{t}\subseteq \EntitySet$ be such that $\SelectionVec_{t}$ is
a partial selection for $V_{t}$ (i.e., $|V_{t}|=|\SelectionVec_{t}|$, and
$|\Supporters_{i}(\System) \cap \SelectionVec_{t}|=1$
for every $\Entity_i \in V_{t}$).
Given that $\SelectionVec_{t}$ is a selection for nodes $\Entity_1, \ldots, \Entity_t$
that satisfies (C1) and (C2), we show that $\SelectionVec_{t+1}$ satisfies
(C1) and (C2) as well.

In particular, it is required to show that there exists at least one supporter
of $\Entity_{t+1}$, namely, $\Affectors_k \in \Supporters_{t+1}(\System)$,
such that $\fff(\beta^-,\System(\SelectionVec_{t} \cup \{\Affectors_k\}))=1$.
This will imply that step 7(a) of the algorithm always succeeds in expanding $\SelectionVec_{t}$.

By Observation \ref{obs:irreducible_selection} and Property (C2) for step $t$,
the system $\System(\SelectionVec_{t})$ is irreducible with $\beta_{t} \geq \beta^{-}$.
In addition, note that
$\FilterMatrixFamily(\System_{t}) \subseteq \FilterMatrixFamily(\System)$
(as every square system of $\System_t$ is also a square system of $\System$).
By Theorem \ref{thm:pf_ext}, there exists a square system
$\System_{t}(\FilterMatrix^{*}_{t})$,
$\FilterMatrix^{*}_{t} \in \FilterMatrix(\System_{t})$, such that
$\PFEigenValue(\System_{t}(\FilterMatrix^{*}_{t}))=1/\beta_{t}$.
In addition, $\PFEigenVector(\System_{t}(\FilterMatrix^{*}_{t}))$ is a feasible
solution for Program (\ref{LP:Ext_Perron_convex}) with the system
$\System_{t}(\FilterMatrix^{*}_{t})$ and $\beta=\beta_{t}$.

By Eq. (\ref{eq:FilterMatrixFamily}), the square system
$\System_{t}(\FilterMatrix^{*}_{t})$ corresponds to a complete selection
$\SelectionVec^{**}$, where $|\SelectionVec^{**}|=n$ and $\SelectionVec_{t} \subseteq \SelectionVec^{**}$,
i.e., $\System_{t}(\FilterMatrix^{*}_{t})=\System(\SelectionVec^{**})$.
Observe that by Property (Q5) of Thm. \ref{thm:pf_ext} for the system
$\System_{t}$,
there exists a $\ZeroStar$ solution for Program (\ref{LP:Ext_Perron_convex})
that achieves $\beta_{t}$. This  $\ZeroStar$ solution is constructed from
$\PFEigenVector(\System_{t}(\SelectionVec^{**}))$,
the PF eigenvector of  $\System_{t}(\SelectionVec^{**})$.

Let $\Affectors_k \in \Supporters_{t+1}(\System_{t}) \cap \SelectionVec^{**}$.
Note that by the choice of $\SelectionVec^{**}$, such an affector
$\Affectors_k$ exists.
We now show that  $\SelectionVec_{t+1}=\SelectionVec_{t} \cup \{\Affectors_k\}$
satisfies Property (C2), thus establishing the existence of
$\Affectors_k \in \Supporters_{t+1}(\System_{t})$ in step 7(a).
We show this by constructing a feasible solution $X^{*}_{\beta^-} \in \R^{m(\SelectionVec_{t+1})}$
for $\System_{t+1}$. By the definition of $\SelectionVec^{**}$,
$\fff(\beta^-,\System(\SelectionVec^{**}))=1$ and therefore there exists
a feasible solution $\overline{X}^{t+1}_{\beta^-} \in \R^{n}$
for $\System(\SelectionVec^{**})$.
Since $\SelectionVec_{t+1}  \subseteq  \SelectionVec^{**}$, it is possible to  extend
$\overline{X}^{t+1}_{\beta^-} \in \R^{n}$ to a feasible solution $X^{*}_{\beta^{-}}$
for system $\System_{t+1}$,
by setting $X^{*}_{\beta^{-}}(\Affectors_q)=X^{t+1}_{\beta^{-}}(\Affectors_q)$ for every
$\Affectors_q \in \SelectionVec^{**}$ and $X^{*}_{\beta^{-}}(\Affectors_q)=0$ otherwise.
It is easy to verify that this is indeed a feasible solution for $\beta^{-}$,
concluding that
$\fff(\beta^-,\System_{t+1})=1$.

So far, we have shown that there exists
an affector $\Affectors_k \in \Supporters_{t+1}(\System_t)$
such that $\fff(\beta^-,\System_{t+1})=1$.
We now claim that for any $\Affectors_k \in \Supporters_{t+1}(\System_t)$
such that $\fff(\beta^-,\System_{t+1})=1$,
Properties (C1) and (C2) are satisfied. This holds trivially,
relying on the criterion for selecting $\Affectors_k$, since $\Supporters_{t+1}(\System_t) \cap \SelectionVec_{t}=\emptyset$.

After $n$ steps, we get that $\SelectionVec_n$ is a complete selection,
$\FilterMatrix(\SelectionVec_n) \in \FilterMatrixFamily(\System_{n-1})$, and therefore
by Property (C1) for steps $t=1,\ldots,n$, it also holds that
$\FilterMatrix(\SelectionVec_n) \in \FilterMatrixFamily(\System)$.
In addition, by Property (C2), $\fff(\beta^-,\System_n)=1$.
Since $\System_n$ is equivalent to
$\System(\SelectionVec_n) \in \SquareSystemFamily$
(obtained by removing the $m-n$ columns corresponding to the affectors not selected
by $\SelectionVec_n$), it is easy to verify that
$\fff(\beta^-,\System(\SelectionVec_n))=1$.
Next, by Thm. \ref{thm:pf} we have that
$1/\PFEigenValue(\System(\SelectionVec_n)) \in Range_{\beta^{*}}$.

It remains to show that
$1/\PFEigenValue(\System(\SelectionVec_n))= \beta^{*}(\System)$.
By Theorem \ref{thm:pf_ext}, there exists a square system
$\System(\FilterMatrix^{*})$, $\FilterMatrix^{*} \in \FilterMatrix(\System)$,
such that
$\PFEigenValue(\System(\FilterMatrix^{*}))=1/\beta^{*}$.
Assume, toward contradiction, that
$1/\PFEigenValue(\System(\SelectionVec_n)) \neq 1/\beta^{*}$.
Obs. \ref{obs:filter_to_square}(b) implies that
$\PFEigenValue(\System(\FilterMatrix^{*})) <
\PFEigenValue(\System(\SelectionVec_n))$.
It therefore follows that $\System(\FilterMatrix^{*})$ and
$\System(\SelectionVec_n)$ are two non-equivalent hidden square systems
of $\System$ such that
$1/\PFEigenValue(\System(\FilterMatrix^{*})),
1/\PFEigenValue(\System(\SelectionVec_n)) \in Range_{\beta^{*}}$, or, that
$1/\PFEigenValue(\System(\SelectionVec_n))-
1/\PFEigenValue(\FilterMatrix^{*}) \leq \Delta_{\beta}$,
in contradiction to Lemma \ref{lem:apart_in_range}.
This completes the proof of Lemma \ref{cl:selection_opt}.
\QED
We are now ready to complete the proof of Thm. \ref{thm:algorithm}.
%and hence also Theorem \ref{thm:algorithm}.
\Proof [Theorem \ref{thm:algorithm}]
We show that Alg. $\AlgoName$ satisfies
the requirements of the theorem.
By Obs. \ref{obs:filter_to_square}(b),
$\min_{\FilterMatrix \in \FilterMatrixFamily} \left\{\PFEigenValue(\System(\FilterMatrix)) \right\}
\geq 1/\beta^{*}(\System)$.
Therefore, since
$\PFEigenValue(\System(\SelectionVec_n))=1/\beta^{*}(\System)$,
the square system $\System(\SelectionVec_n)$ constructed in step 7 of the algorithm indeed yields the Perron value (by Eq. (\ref{eq:general_pf_root})), hence the correctness of the algorithm is established.
\par
Finally, we analyze the runtime of the algorithm.
Note that there are  $O(\log \left(\beta^{*}(\System)/\Delta_{\beta}\right)+n)$
calls for the interior point method (computing $\fff(\beta^-,\System_{i})$),
namely,
$O(\log \left(\beta^{*}(\System)/\Delta_{\beta}\right))$ calls in the first
phase and $n$ calls in the second phase. By plugging Eq. (\ref{eq:gmax})
in Claim \ref{cl:max_beta}, Thm. \ref{thm:algorithm} follows.
\QED

%%%%%%%%%%%%%%%%%%%%%%%%%%%%%%%%%%%%%%%%%%%%
\section{Limitations for the existence of a $\ZeroStar$}
\label{sec:limit}
%%%%%%%%%%%%%%%%%%%%%%%%%%%%%%%%%%%%%%%%%%%%
In this section we provide a characterization of systems in which
a $\ZeroStar$ solution does not exist.
%%%%%%%%%%%%%%%%%%%%%%%%%%%%%%%%%%%%%%%%%%%%
%\subsection{Bounded value systems.}
%\label{subsec:BoundedValue}
\paragraph{Bounded value systems.}
%%%%%%%%%%%%%%%%%%%%%%%%%%%%%%%%%%%%%%%%%%%%
Let $\MaxPower$ be a fixed constant.
For a nonnegative vector $\overline{X}$, let
$$\max(\overline{X}) ~=~
\max\left\{X(j)/X(i) \mid 1 \leq i,j \leq n, X(i)>0\right\}.$$
%\begin{definition}
A system $\System$ is called
a \emph{bounded power system} if $\max(\overline{X}) \leq \MaxPower$.
%\end{definition}
\begin{lemma}
\label{lem:bounded_system}
There exists a bounded power system $\System$ such that no optimal solution
$\overline{X}^{*}$ for $\System$ is a $\ZeroStar$ solution.
\end{lemma}
\Proof
Consider the optimization problem (\ref{LP:Ext_Perron}), and the following
system $\System=\langle \SupportersMatrix,\RepressorsMatrix\rangle$:
\begin{align*}
\SupportersMatrix &=
\left(\begin{array}{cccc}a&a&0&0\\0&0&a&a\end{array}\right),
~~~~~~~~~~~
\RepressorsMatrix =
\left(\begin{array}{cccc}0&0&4c\MaxPower^2&4c\MaxPower^2\\c&c&0&0\end{array}\right),
\end{align*}
for constants $a,c>0$. We first show that it is impossible to attain
the optimal value $\beta^{*}$
if $\overline{X}$ is a $\ZeroStar$ solution. Then, we show that there exists
a non-$\ZeroStar$ solution $\overline{X}$ that attains $\beta^{*}$.
Thus, for a given system, no $\ZeroStar$ solution is optimal.

Assume, by contradiction, that we have a $\ZeroStar$ solution that achieves
$\beta^{*}$ on $\System$. Due to symmetry, every $\ZeroStar$ solution will
yield the same $\beta^{*}$, so without loss of generality assume that $X(2)=0$
and $X(4)=0$, and thus the corresponding square system is
\begin{align*}
\widehat{\SupportersMatrix} &=
\left(\begin{array}{cc}a&0\\0&a\end{array}\right),
~~~~~~~~~~
\widehat{\RepressorsMatrix} =
\left(\begin{array}{cc}0&4c\MaxPower^2\\c&0\end{array}\right).
\end{align*}
By Lemma \ref{lem:strict_equality}, at the optimum value $\beta^{*}$,
the inequality constraints of Eq. (\ref{eq:SR}) holds with equality, namely,
$\left(\widehat{\RepressorsMatrix}-\frac{1}{\beta^{*}}\widehat{\SupportersMatrix}\right)
\cdot \overline{X}=0$. Plugging in the chosen values, we get
\begin{align*}
\left(\begin{array}{cc}
-\frac{a}{\beta}&4c\MaxPower^2
\\
c&-\frac{a}{\beta}
\end{array}\right)
\cdot
\left(\begin{array}{c}
X(1)
\\
X(3)
\end{array}\right)
=0~,
\end{align*}
leading to the equations
%\begin{align*}
%\begin{array}{l}
$-\frac{a}{\beta}X(1)+4c\MaxPower^2X(3)=0$
and
%\\
$cX(1)-\frac{a}{\beta}X(3)=0$.
%\end{array}
%\end{align*}
Rewriting these two equations as
$\frac{X(1)}{X(3)} = 4c\MaxPower^2 / (a/\beta)$
and
$\frac{X(1)}{X(3)} = (a/\beta) / c$~,
%$\frac{X(1)}{X(3)}=\frac{4c\MaxPower^2}{a/\beta}$ and
%$\frac{X(1)}{X(3)}=\frac{a/\beta}{c}$~,
we get that $\left(\frac{X(1)}{X(3)}\right)^2 =
4\MaxPower^2$, or, $\frac{X(1)}{X(3)} =
2\MaxPower$.
But this contradicts the assumption that $\System$ is a bounded value
system, namely, $\max(\overline{X}) \leq \MaxPower$.
It follows that there is no optimal $\ZeroStar$ solution for such a system.

Now we show that there exists a non-$\ZeroStar$ solution $\overline{X}$ for
$\System$ that achieves $\beta^{*}$. Consider some $\overline{X}$ satisfying
$X(2)=0,X(1)>0,X(3)>0$ and $X(4)>0$. Similar to the above steps, we derive that
$\frac{X(1)}{X(3)+X(4)} = \frac{\beta c}{a}$ and
$\frac{X(3)+X(4)}{X(1)} = \frac{4\beta c \MaxPower^2}{a}~$, hence
$\left(\frac{X(3)+X(4)}{X(1)}\right)^2 = 4\MaxPower^2$,
or, $\frac{X(3)+X(4)}{X(1)} =2\MaxPower$.
Clearly, the last equation does not contradict the value boundedness of
$\System$, since $\max(\overline{X}) \leq \MaxPower$ only imposes the
constraint $\frac{X(3)+X(4)}{X(1)}\le2\MaxPower$. It follows that there exists
a non-$\ZeroStar$ solution that attains $\beta^{*}$.
\QED

%%%%%%%%%%%%%%%%%%%%%%%%%%%%%%%%%%%%%%%%%%%%
\paragraph{Second eigenvalue maximization.}
%%%%%%%%%%%%%%%%%%%%%%%%%%%%%%%%%%%%%%%%%%%%
One of the most common applications of the \PFT~is the existence of the
stationary distribution for a transition matrix
(representing a random process).
The stationary distribution is the eigenvector of the largest eigenvalue of
the transition matrix. We remark that if the transition matrix is stochastic,
i.e., the sum of each row is $1$, then the  largest eigenvalue is equal to $1$.
So this case does not give rise to any optimization problem.
However, in many cases we are interested in processes with fast mixing time.
Assuming the process is ergodic, the mixing time is determined by the
difference between the largest eigenvalue and the second largest eigenvalue.
So we can try to solve the following problem. Imagine that there is some rumor
that we are interested in spreading over two or more social networks.
Each node can be a member of several social networks. We would like to merge
all the networks into one large social network in a way that will result in fast
mixing time. This problem looks very similar to the one solved in this paper.
Indeed, one can use similar techniques and get an approximation.
But interestingly, this problem does not have the $\ZeroStar$ solution property,
as illustrated in the following example.

Assume we are given $n$ nodes. Consider the $n!$ different social
networks that arise by taking, for each permutation $\pi \in S(n)$, the path $P_\pi$ corresponding to the permutation $\pi$. Clearly, the best mixing graph
we can get is the complete graph $K_n$. We can get this graph if each node
chooses each permutation with probability $\frac{1}{n!}$.
We remind the reader that the mixing time of the graph $K_n$ is 1.
On the other hand, any $\ZeroStar$ solution have a mixing time $O(n^2)$.
This example shows that in the second largest eigenvalue, the solution
is not always a $\ZeroStar$ solution.

%%%%%%%%%%%%%%%%%%%%%%%%%%%%%%%%%%%%%%
\section{Applications}
\label{short:sec:Applications}
%%%%%%%%%%%%%%%%%%%%%%%%%%%%%%%%%%%%%%%%%%%%
We have considered several applications for our generalized \PFT.
All these examples concern generalizations of well-known applications
of the standard \PFT.
%Section \ref{sec:Applications} illustrates applications
In this section, we illustrate applications
for power control in wireless networks,
and input--output economic model.
(In fact, our initial motivation for the study of generalized \PFT~arose
while studying algorithmic aspects of wireless networks
in the SIR model \cite{Avin2009PODC,KLPP2011STOC, Avin2012SINR}.)
%In the future, change  Avin2009PODC  to  Avin10journal
%
%For each application, we briefly state what are the entities, affectors,
%supporters and repressors and explain the (physical) meaning  of the gain
%and the \PFR~ and \PFV.
%We provide this example in detail in Appendix \ref{subsec:power_control_app}.
%In addition, (A3) is deferred to Section \ref{subsec:Population} as well.
%

%%%%%%%%%%%%%%%%%%%%%%%%%%%%%%%%%%%%%%%%%%%%
\subsection{Power control in wireless networks.}
\label{subsec:power_control_app}
%%%%%%%%%%%%%%%%%%%%%%%%%%%%%%%%%%%%%%%%%%%%
%\def\ApplicationPC{
The rules governing the availability and quality of wireless connections
can be described by {\em physical} or {\em fading channel} models
(cf. \cite{PL95,B96,R96}). Among those, a commonly studied is the
{\em signal-to-interference ratio (SIR)} model\footnote{This is
a special case of the {\em signal-to-interference \& noise ratio (SINR)} model
where the noise is zero.}.
In the SIR model, the energy of a signal fades with the distance
to the power of the {\em path-loss parameter} $\alpha$.
If the signal strength received by a device divided by the interfering
strength of other simultaneous transmissions
%(plus the fixed \emph{background
%noise} \(\Noise\)\footnote{In the sequel we assume $N=0$.})
is above some \emph{reception threshold} $\beta$, then the
receiver successfully receives the message, otherwise it does not. Formally,
let $\dist{p,q}$ be the Euclidean distance between $p$ and $q$,
and assume that each transmitter $t_i$ transmits with power $X_i$.
At an arbitrary point $p$, the transmission of station $t_i$ is
correctly received if
%\begin{align}\label{eq:sinr}
%\frac{X_i \cdot \dist{p, t_i}^{-\alpha}}
%{N + \sum_{j \neq i} X_j \cdot \dist{p, t_j}^{-\alpha}}
%~ \geq ~ \beta ~ .
%\end{align}
\begin{align}\label{eq:sinr}
\frac{X_i \cdot \dist{p, t_i}^{-\alpha}}
{\sum_{j \neq i} X_j \cdot \dist{p, t_j}^{-\alpha}}
~ \geq ~ \beta ~ .
\end{align}
In the basic setting, known as the SISO (Single Input, Single Output) model,
we are given a network of $n$ receivers $\{r_i\}$ and transmitters $\{t_i\}$
embedded in $\mathbb{R}^d$ where each transmitter is assigned to a single
receiver. The main question is then is to find the optimal
(i.e., largest) $\beta^*$ and the power assignment $\overline{X^{*}}$ that
achieves it when we consider Eq. (\ref{eq:sinr}) at each receiver $r_i$.
The larger $\beta$, the simpler (and cheaper) is the hardware implementation
required to decode messages in a wireless device.  In a seminal and elegant
work, Zander \cite{Zander92b} showed how to compute  $\beta^*$ and
$\overline{X}^{*}$, which are essentially the \PFR~ and \PFV,
if we generate a square
matrix $A$ that captures the signal and interference for each station.

The motivation for the general \PFT\ appears when we consider Multiple Input
Single Output (MISO) systems.
In the MISO setting, a set of multiple synchronized transmitters,
located at different places, can transmit at the same time to the same receiver.
Formally, for each receiver $r_i$ we have a set of $k_i$ transmitters,
to a total of $m$  transmitters. Translating this to the generalized \PFT,
the $n$ receivers are the entities and the $m$ transmitters are affectors.
For each receiver, its supporter set consists of its $k_i$ transmitters and
its repressor set contains all other transmitters.
The SIR equation at receiver $r_i$ is then:
\begin{align}\label{eq:sinr2}
\frac{\sum_{\ell \in \Supporters_i} X_{\ell} \cdot \dist{r_i, t_{\ell}}^{-\alpha}}
{\sum_{\ell \in \Repressors_i} X_{\ell} \cdot \dist{r_i, t_{\ell}}^{-\alpha}}
~ \geq ~ \beta ~,
\end{align}
where $\Supporters_i$ and $\Repressors_i$ are the sets of supporters and
repressors of $r_i$, respectively.
As before, the gain $g(i,j)$ is proportional to $1/\dist{r_i, t_j}^{-\alpha}$
(where the sign depends on whether $t_j$ is a supporter or repressor of $r_i$).
Using the generalized \PFT~we can again find the optimal reception threshold
$\beta^*$ and the power assignment $\overline{X}^{*}$ that achieves it.

% Perhaps place Figure \ref{figure:sinr_diag}} here?

An interesting observation is that since our optimal power assignment is
a $\ZeroStar$ solution using several transmitters at once for a receiver
is not necessary, and will not help to improve $\beta^*$, i.e.,
only the ``best" transmitter of each receiver needs to transmit
(where ``best" is with respect to the entire set of receivers).

%%%%%%%%%%%%%%%%%%%%%%%%%%%%%%%%%%%%%%%%%%%%
%\subsubsection{Related Work on MISO Power Control}
%\label{sec:related_work_miso_power}
\paragraph{Related work on MISO power control.}
%%%%%%%%%%%%%%%%%%%%%%%%%%%%%%%%%%%%%%%%%%%%
We next highlight the differences between our proposed
MISO power-control algorithm and the existing approaches to this problem.
The vast literature on power control in MISO and MIMO systems
considers mostly the joint optimization of power control with beamforming
(which is represented by a precoding and shaping matrix).
In the commonly studied {\em downlink scenario}, a single transmitter
with $m$ antennae sends independent information signals to $n$ decentralized
receivers. With this formulation, the goal is to find an optimal power vector
of length $n$ and a $n \times m$ beamforming matrix.
The standard heuristic applied to this problem is an iterative strategy
that alternatively repeats a {\em beamforming} step
(i.e., optimizing the beamforming matrix while fixing the powers)
and a {\em power control} step
(i.e., optimizing powers while fixing the beamforming matrix)
till convergence \cite{CaiQT11,Cai2011,Chiang2007,Schu2004,Chee11}.
In \cite{CaiQT11}, the geometric convergence of such scheme has been
established.  In addition, \cite{WieselES06} formalizes the problem
as a conic optimization program that can be solved numerically.
In summary, the current algorithms for MIMO power-control (with beamforming)
are of numeric and iterative flavor, though with good convergence guarantees.
In contrast, the current work considers the simplest MISO setting
(without coding techniques) and aims at \emph{characterizing} the mathematical
\emph{structure} of the optimum solution. In particular, we establish the fact
that the optimal max-min SIR value is an algebraic number
(i.e., the root of a characteristic polynomial) and the optimum power vector
is a $\ZeroStar$ solution. Equipped with this structure, we design
an efficient algorithm which is more accurate than off-the-shelf
numeric optimization packages that were usually applied in this context.
Needless to say, the structural properties of the optimum solution are of
theoretical interest in addition to their applicability.

We note that our results are (somewhat) in contradiction to the well-established
fact that MISO and MIMO (Multiple Input Multiple Output) systems, where
transmitters transmit in parallel, do improve the capacity of wireless networks,
which corresponds to increasing $\beta^*$ \cite{Foschini98onlimits}.
There are several reasons for this apparent dichotomy, but they are all
related to the simplicity of our SIR model. For example, if the ratio
between the maximal power to the minimum power is bounded, then our result
does not hold any more (as discussed in Section \ref{sec:limit}).
In addition, our model does not capture random noise and small scale fading
and scattering \cite{Foschini98onlimits}, which are essential for the benefits
of a MIMO system to manifest themselves.

%%%%%%%%%%%%%%%%%%%%%%
\subsection{Input--output economic model.}
%\paragraph{Input--Output Economic Model (Leontief's Model)}
%%%%%%%%%%%%%%%%%%%%%%
Consider a group of $n$ industries that each produce (output) one type of commodity,
but requires inputs from other industries
\cite{meyer2000matrix,pillai2005pft}.
Let $a_{ij}$ represent the number of $j$th industry commodity units that need to be purchased
by the $i$th industry to operate its factory for one time unit divided by the number of
commodity units produced by the $i$th industry in one time unit, where $a_{ij} \ge 0$.

Let $X_j$ represent a unit price of the $i$th commodity to be determined
by the solution.
In the following profit model (variant of Leontief's Model \cite{pillai2005pft}), the percentage
profit margin of an industry for a time unit is:
$$\beta_i ~=~ \text{Profit} ~=~ \text{Total income}/\text{Total expenses}.$$
%\\
%&= X_i - \sum_{j=1}^n a_{ij}X_j
That is, $\beta_i = X_i /\left(\sum_{j=1}^n a_{ij}X_j\right)$.
Maximizing the the profit of each industry can be solved
via Program (\ref{LP:Stand_Perron}), where $\beta^*$ is the minimum profit and  $\overline{X}^{*}$ is the optimal pricing.

Consider now a similar model where the $i$th industry can produce $k_i$
alternative commodities in a time unit and requires inputs from other commodities of industries.
The industries are then the entities in the generalized Perron--Frobenius
setting, and for each industry, its own commodities are the supporters and
input commodities are optional repressors.

The repression gain $\RepressorsMatrix(i,j)$ of industry $i$ and commodity $j$ (produced by some other industry $i'$), is the number of $j$th commodity units that are required by the $i$th industry to produce (i.e., operate) for a one unit of time. Thus, $(\RepressorsMatrix \cdot \overline{X})_i$ is the total expenses of industry $i$ in one time unit.

The supporter gain $\SupportersMatrix(i,j)$ of industry $i$ to its commodity $j$ is the number of units it can produce in one time unit.
Thus, $(\SupportersMatrix \cdot \overline{X})_i$ is the total income of industry $i$ in one time unit.
Now, similar to the basic case, $\beta^*$ is the best minimum percentage profit for an industry and  $\overline{X}^{*}$ is the optimal pricing for the commodities. The existence of a $\ZeroStar$ solution implies that it is sufficient for each industry to charge a nonzero cost for only \emph{one} of its commodities and produce the rest for free.
%%%%%%%%%%%%%%%%%%%%%%%%
\section{Discussion and open problems}
%%%%%%%%%%%%%%%%%%%%%%%%
Our results concern the generalized eigenpair of a nonsquare system
of dimension $n \times m$, for $m \geq n$.
We provide a definition, as well as a geometric and a graph theoretic characterization
of this eigenpair, and present centralized algorithm for computing it.
A natural question for future study is whether there exists an iterative
method with a good convergence guarantees for this task, as exists for
(the maximal eigenpair of) a square system.
In addition, another research direction involves studying the other eigenpairs
of a nonsquare irreducible system. In particular, what might be the meaning
of the 2nd eigenvalue of this spectrum?
Yet another interesting question involves studying the relation of our
spectral definitions with existing spectral theories for nonsquare matrices.
Specifically, it would be of interest to characterize the relation between
the generalized eigenpairs of irreducible systems according to our definition
and the eigenpair provided by the SVD approach. Finally, we note that
a setting in which $n < m$ might also be of practical use (e.g., for the
power control problem in \emph{Single Input Multiple Output} systems), and therefore deserves exploration.

\clearpage

\renewcommand{\baselinestretch}{1}

\end{document}